\documentclass[11pt]{memoir}
\usepackage{amsmath,amsthm,graphicx}
\usepackage{amssymb}

\makeatletter
\def\gacsparsep{6pt}

\newif\if@pfLongIndent %  true if indent proof to width of label.
\@pfLongIndentfalse
\newcommand{\pflongindent}{\@pfLongIndenttrue}
\newcommand{\pfshortindent}{\@pfLongIndentfalse}

\newlength{\pfindent}      
\newlength{\pftopsep}  
\newlength{\pfbotsep}  
\newlength{\pfsep}  
\newlength{\stepsep}        
\newbox{\pfbox}
\newcounter{pfhidelevel}
\newcommand{\pfshortnumbers}[1]{\pfshortNumberLevel=#1\relax}

\newcount\pfshortNumberLevel \pfshortNumberLevel=0

\newcommand{\pfstepnumber}[3]{%
   \ifnum \pfLevelCount < \pfshortNumberLevel
       #3%
     \else $\langle#1\rangle#2$%
   \fi}

\newcommand{\pflevelnumber}[2]{%
   \ifnum \pfLevelCount < \pfshortNumberLevel
       #2%
     \else $\langle#1\rangle$%
   \fi}

\newif\if@pfSideNumbers %  true if putting long numbers at margin
\@pfSideNumbersfalse
\newcounter{pf@sidenumberdepth}
\newlength{\pf@sidenumberoutdent}
\newcommand{\pfsidenumbers}[2]{\@pfSideNumberstrue
     \setcounter{pf@sidenumberdepth}{#1}%
     \addtocounter{pf@sidenumberdepth}{-1}%
     \setlength{\pf@sidenumberoutdent}{#2}}
\newcommand{\pfnosidenumbers}{\@pfSideNumbersfalse}

\newcount\pfLevelCount  \pfLevelCount=0 % current level number
\newcount\pfStepCount   \pfStepCount=0  % current step number
\newcommand{\pfLongLevel}{} % The long version of current level--e.g., 2.7.5
\newcommand{\pfLongStep}{}  % The long version of current step--e.g., 2.7.5.2
\newcommand{\pfStepName}{}  % {current level name}{current step name} 

\newcommand{\pfSetName}{%
  \edef\pfLongStep{%
    \ifnum\pfLevelCount>\@ne
      \pfLongLevel\pfdot\the\pfStepCount
     \else\the\pfStepCount\fi}%
  \edef\pfStepName{%
   \ifnum \pfLevelCount < \pfshortNumberLevel
    {\pfLongStep}{\pfLongStep}%
    \else
     {$\langle\the\pfLevelCount\rangle$}%
     {$\langle\the\pfLevelCount\rangle\the\pfStepCount$}%
   \fi}}

\newcommand{\pfSetRef}[1]{%
  \@ifundefined{pf@#1}%
    {\expandafter
      \edef\csname pf@#1\endcsname{\pfStepName}}%
    {\typeout{WARNING:
         proof step "#1" (<\the\pfLevelCount>\the\pfStepCount) 
        already defined}}%
    }

\newcommand{\pfPrintStepNumber}[2]{#2}
\newcommand{\pfPrintLevelNumber}[2]{#1}

\newcommand{\stepref}[1]{\@ifundefined{pf@#1}{{\bf ??}\typeout{WARNING:
 proof step "#1" undefined}}{\expandafter\expandafter\expandafter
 \pfPrintStepNumber\csname pf@#1\endcsname}}

\newcommand{\levelref}[1]{\@ifundefined{pf@#1}{{\bf ??}\typeout{WARNING:
 proof step #1 undefined}}%
  {\expandafter\expandafter\expandafter
 \pfPrintLevelNumber\csname pf@#1\endcsname}}

\newlength{\pf@outdent}
\newcommand{\pfSideNumber}{%
  \if@pfSideNumbers 
   \ifnum\pfLevelCount>\value{pf@sidenumberdepth}%
    \hspace*{-\pf@outdent}%
    \makebox[0pt][l]{\footnotesize\pfLongStep}%
    \hspace*{\pf@outdent}%
    \else\fi
   \else\fi}

 \newcounter{branch}%

 \newcounter{branchi}[branch]

 \newcounter{branchii}[branchi]

 \newcounter{branchiii}[branchii]

 \newcounter{branchiv}[branchiii]

 \newcounter{branchv}[branchiv]

 \newcounter{branchvi}[branchv]

 \def\CurBranch{{branch\romannumeral\pfLevelCount}}%

 \def\TheCurBranch{\csname thebranch\romannumeral\pfLevelCount\endcsname}%

 \def\RefStepCounter#1{\refstepcounter{#1}} % for expandafter

 \def\SetTheBranch#1#2{%
  \expandafter\def\csname thebranch#2\endcsname{%
    \csname thebranch#1\endcsname\pfdot\arabic{branch#2}%
   }}

 \def\OpenLevel{%
  \edef\ParentBranch{\CurBranch}
  \edef\romanPrevLevel{{\romannumeral \pfLevelCount}}
  \advance\pfLevelCount\@ne%
  \edef\romanLevel{{\romannumeral\pfLevelCount}}\@ifundefined
{c@branch\romannumeral\pfLevelCount}{\typeout
{WARNING: too many prooof levels}}{}\expandafter\setcounter\CurBranch{0}
}

\def\pflist#1#2{\ifnum \@listdepth >15\relax \@toodeep 
     \else \global\advance\@listdepth\@ne \fi
  \rightmargin \z@ \listparindent\z@ \itemindent\z@
  \csname @list\romannumeral\the\@listdepth\endcsname 
  \def\@itemlabel{#1}\let\makelabel\@mklab \@nmbrlistfalse #2\relax
  \@trivlist
  \parskip\parsep \parindent\listparindent
  \advance\linewidth -\rightmargin \advance\linewidth -\leftmargin
  \advance\@totalleftmargin \leftmargin
  \parshape \@ne \@totalleftmargin \linewidth 
  \ignorespaces}

\def\endpflist{\global\advance\@listdepth\m@ne
    \endtrivlist}

\let\@listvii=\@listv
\let\@listviii=\@listv
\let\@listix=\@listv
\let\@listx=\@listv
\let\@listxi=\@listv
\let\@listxii=\@listv
\let\@listxiii=\@listv
\let\@listxiv=\@listv
\let\@listxv=\@listv
\let\@listxvi=\@listv

\newenvironment{prooof}{%       BEGIN:
  \edef\pfLongLevel          % LongLevel := 
   {\ifnum\pfLevelCount>\z@  %    IF LevelCount > 0
     \ifnum\pfLevelCount>\@ne%     THEN IF LevelCount > 1
      \pfLongLevel\pfdot\else\fi  %            THEN LongLevel * "." FI
      \the\pfStepCount       %          * StepCount
     \else\fi}%              %     ELSE FI
  \OpenLevel                 % Gacs replacement for the above line
  \@tempcnta=\value{pfhidelevel}%
  \ifnum\@tempcnta<\@ne
    \setcounter{pfhidelevel}{1}%
    \typeout{WARNING: pfhidelevel < 1, setting to 1}%
    \@tempcnta=\@ne\fi
  \advance\@tempcnta\@ne
  \if@qedstep
    \advance\@tempcnta\@ne
    \ifnum\pfLevelCount 
        = \@tempcnta 
     \sbox{\pfbox}\bgroup\begin{minipage}{\textwidth}\fi
  \else     
   \ifnum\pfLevelCount 
        = \@tempcnta 
     \sbox{\pfbox}\bgroup\begin{minipage}{\textwidth}\fi\fi
  \pfStepCount=\z@           %   StepCount  := 0
  \ifnum\pfLevelCount>\@ne   %   IF LevelCount > 1
   \begin{pflist}{}{%          %     THEN Begin List with
    \topsep=\pfsep\relax     %           \topsep     := \pfsep
    \itemsep=\z@             %           \itemsep    := 0
     \parsep=\gacsparsep             % Gacs
    \partopsep=\z@           %           \partopsep  := 0
    \if@pfLongIndent         %          IF LongIndent
     \settowidth{\leftmargin}%%           THEN
      {\expandafter          %              \leftmargin := width of step name
       \pfPrintStepNumber    %                             + \labelsep
       \pfStepName.}%
     \advance\leftmargin
     \labelsep 
    \else                    %           ELSE
    \leftmargin=\pfindent    %            \leftmargin := \pfindent
    \fi\relax        
   }   \item[]
   \else \par                %     ELSE \par
    \addvspace{\pftopsep}%   %          addvspace of \pftopsep
    \parindent=\z@           %          \parindent := 0
    \parskip = \z@           %          \parskip := 0
    \@ifundefined{mathindent}{%
    }{%
      \mathindent=1em}%
 \fi
\@qedstepfalse
    }%                       % END:
  {\ifnum\pfLevelCount>\@ne  %   IF LevelCount > 1
    \end{pflist}%              %     THEN End List
   \else \par                %     ELSE \par
     \addvspace{\pfbotsep}%  %          addvspace of \pfbotsep
   \fi
  \@tempcnta=\value{pfhidelevel}%
  \advance\@tempcnta\@ne
  \if@qedstep
    \advance\@tempcnta\@ne
    \ifnum\pfLevelCount 
        = \@tempcnta 
     \end{minipage}\egroup\sbox{\pfbox}{}\@qedstepfalse\fi
  \else 
      \ifnum\pfLevelCount 
        =\@tempcnta 
        \end{minipage} 
        \egroup\sbox{\pfbox}{}\fi
 \fi}

  {\par}                     %  \par

\newcommand{\step}[2]{\begin{step+}{#1}#2\end{step+}}
\newif\if@qedstep %  true right after a \qedstep if in the scope of an 
\@qedstepfalse

\newif\if@unhideqedstep % set true by \unhideqprooof, false by \hideqprooof. 

\@unhideqedstepfalse
\newcommand{\unhideqedprooof}{\@unhideqedsteptrue}
\newcommand{\hideqedprooof}{\@unhideqedstepfalse}

\newcommand{\qedstep}{\step{qedstep\the\pfLevelCount}{Q.E.D.}
\if@unhideqedstep\@qedsteptrue\fi}

\newcommand{\nostep}[1]{% 
  \expandafter\RefStepCounter\CurBranch    % Gacs addition
  \label{#1}                               % Gacs addition
  \advance\pfStepCount\@ne       % StepCount := StepCount + 1
  \pfSetName                     % SetName
  \pfSetRef{#1}%                 % SetRef
  }

\newenvironment{step+}[1]%
 {\endgroup                      % Move Outside environment scope
  \expandafter\RefStepCounter\CurBranch    % Gacs addition
  \label{#1}                               % Gacs addition
  \advance\pfStepCount\@ne       % StepCount := StepCount + 1
  \pfSetName                     % SetName
  \pfSetRef{#1}%                 % SetRef
 \begingroup\@endpefalse         % Move inside environment scope
   \def\@currenvir{step+}%       %   by simulating \begin{step+}
 \begin{pflist}{}{%                % Begin list environment with
   \setlength
     {\pf@outdent}{\textwidth}%  %   \pf@outdent = outdent
   \addtolength                  %                 for side numbers
    {\pf@outdent}{-\linewidth}%  
   \addtolength
    {\pf@outdent}%
    {\pf@sidenumberoutdent}%
   \topsep=\stepsep\relax        %   \topsep     := \stepsep
   \itemsep=\z@                  %   \itemsep    := 0
   \parsep=\gacsparsep                   %    % Gacs
   \partopsep=\z@                %   \partopsep  := 0
   \settowidth{\labelwidth}%     %   \labelwidth := width of step name.
     {\expandafter
      \pfPrintStepNumber
      \pfStepName.}%
   \leftmargin=\labelwidth\relax %   \leftmargin := \labelwidth
   \advance\leftmargin\labelsep %                   + \labelsep
   \relax}%
 \item[\pfSideNumber%               % \item[ \pfSideNumber
  \expandafter                       %        StepNumber]
  \pfPrintStepNumber\pfStepName.]}%
 {\end{pflist}}

\newenvironment{assume+}{%
 \begin{pflist}{}{%                % Begin list environment with
   \topsep=\z@                   %   \topsep     := 0
   \itemsep=\z@                  %   \itemsep    := 0
   \parsep=\gacsparsep                   %   Gacs
   \partopsep=\z@                %   \partopsep  := 0
   \settowidth{\labelwidth}%     %   \labelwidth := width "Assume:"
     {{\kwfont Assume\/}:}%
   \leftmargin=\labelwidth\relax %   \leftmargin := \labelwidth
   \advance\leftmargin\labelsep  %                   + \labelsep
   \relax}%
 \item[{\kwfont Assume\/}:]}%    % \item[Assume:]]
 {\end{pflist}}

\newenvironment{prove+}{%
 \begin{pflist}{}{%                % Begin list environment with
   \topsep=\z@                   %   \topsep     := 0
   \itemsep=\z@                  %   \itemsep    := 0
   \parsep=\gacsparsep   % Gacs
   \partopsep=\z@                %   \partopsep  := 0
   \settowidth{\labelwidth}%     %   \labelwidth := width "Assume:"
     {{\kwfont Assume\/}:}%
   \leftmargin=\labelwidth\relax %   \leftmargin := \labelwidth
   \advance\leftmargin\labelsep  %                   + \labelsep
   \relax}%
 \item[{\kwfont 
          Prove\/}:\hfill]}%     % \item[Prove:]]
 {\end{pflist}}

\newenvironment{pflet+}{%
 \begin{pflist}{}{%                % Begin list environment with
   \topsep=\z@                   %   \topsep     := 0
   \itemsep=\z@                  %   \itemsep    := 0
   \parsep=\gacsparsep % Gacs
   \partopsep=\z@                %   \partopsep  := 0
   \settowidth{\labelwidth}%     %   \labelwidth := width "Assume:"
     {{\kwfont Let\/}:}%
   \leftmargin=\labelwidth\relax %   \leftmargin := \labelwidth
   \advance\leftmargin\labelsep  %                   + \labelsep
   \relax}%
 \item[{\kwfont 
          Let\/}:\hfill]}%     % \item[Let:]]
 {\end{pflist}}

\newcommand{\case}[1]{%
 \begin{pflist}{}{%                % Begin list environment with
   \topsep=\z@                   %   \topsep     := 0
   \itemsep=\z@                  %   \itemsep    := 0
   \parsep=\gacsparsep % Gacs
   \partopsep=\z@                %   \partopsep  := 0
   \settowidth{\labelwidth}%     %   \labelwidth := width "Case:"
     {{\kwfont Case\/}:}%
   \leftmargin=\labelwidth\relax %   \leftmargin := \labelwidth
   \advance\leftmargin\labelsep  %                   + \labelsep
   \relax}%
 \item[{\kwfont
           Case\/}:]             % \item[Case:]]
 #1
 \end{pflist}}

\newcommand{\pf}{{\kwfont Proof\/}.} % Gacs, to conform with amslatex

\let\kwfont=\scshape
\def\pfdot{.}

\def\@push#1#2{{\let\@nil\relax\let\@elt\relax\xdef#1{#2\@elt#1}}}
\def\@pop#1#2{{\let\@nil\relax\let\@elt\relax
              \xdef#2{\expandafter\@innerhead#1}
              \xdef#1{\expandafter\@innerpop#1}}}
\def\@innerpop#1\@elt#2\@nil{#2\@nil}
\def\@head#1#2{{\let\@elt\relax\xdef#2{\expandafter\@innerhead#1}}}
\def\@innerhead#1\@elt#2\@nil{#1}
\def\newstack#1{\gdef#1{\@nil}}

\newcounter{pf@conjCounter}   % counter for conj* and disj*
\newstack\pf@conj             % stack of counters

\newenvironment{conj*}{%
 \@push\pf@conj{\the\value{pf@conjCounter}}%
 \setcounter{pf@conjCounter}{0}%
 \begin{array}[t]{@{\addtocounter{pf@conjCounter}{1}%
   \mbox{\protect\small\protect\arabic{pf@conjCounter}.}
   \land\;}l@{}}%
 }{%
 \end{array}%
 \@pop\pf@conj\pf@temp
  \setcounter{pf@conjCounter}{\pf@temp}}

\newenvironment{disj*}{%
 \@push\pf@conj{\the\value{pf@conjCounter}}%
 \setcounter{pf@conjCounter}{0}%
 \begin{array}[t]{@{\addtocounter{pf@conjCounter}{1}%
   \mbox{\protect\small\protect\alph{pf@conjCounter}.}
   \lor\;}l@{}}%
 }{%
 \end{array}%
 \@pop\pf@conj\pf@temp
  \setcounter{pf@conjCounter}{\pf@temp}}

\newcounter{pfenum}
\newcounter{pfenumdepth}
\newlength{\enumindent}
\@definecounter{pfenumi}
\@definecounter{pfenumii}
\@definecounter{pfenumiii}
\def\thepfenumi{\arabic{pfenumi}}

\def\thepfenumii{\alph{pfenumii}}
\def\p@pfenumii{\thepfenumi}

\def\p@pfenumiii{\thepfenumi\thepfenumii}

\newcommand{\pf@setEnumWidth}[1]{%
  \settowidth{#1}{\setcounter{\@pfenumctr}{2}%
  \csname the\@pfenumctr\endcsname.%
  \setcounter{\@pfenumctr}{0}}}

\newcommand{\pf@enumLabel}{%
  \hfill\makebox[0pt][r]{\csname the\@pfenumctr\endcsname.}}

\newenvironment{pfenum*}{%    % BEGIN
  \ifnum \value{pfenumdepth}>2%  %  IF pfenumdepth > 2
    \relax\@toodeep \else        %    THEN error
  \addtocounter{pfenumdepth}{1}% %    ELSE  pfenumdepth := pfenumdepth + 1
  \edef\@pfenumctr{pfenum%           %          @pfenumctr  := pfenumN
    \romannumeral\the            %            where N = Roman(pfenumdepth)
    \value{pfenumdepth}}%        %
   \fi                           %  FI
  \begin{pflist}%                  %   Begin list environment with
  {\pf@enumLabel}{%              %    Default label = pf@enumlabel
   \labelsep=                    %   \labelsep =
    \ifcase\value{pfenumdepth}   %     CASE OF pfenumdepth
      \labelsep                  %        
     \or .67\labelsep            %       1 -> .67\labelsep
     \or .67\labelsep            %       2 -> .67\labelsep
     \else \labelsep             %      >2 -> \labelsep
     \fi                         %
   \topsep=\z@                   %    \topsep     := 0
   \itemsep=\z@                  %    \itemsep    := 0
   \parsep=\gacsparsep % Gacs
   \partopsep=\z@                %    \partopsep  := 0
   \pf@setEnumWidth\labelwidth   %    set \labelwidth with setEnumWidth
   \leftmargin=\labelwidth\relax %    \leftmargin := \labelwidth
   \advance\leftmargin\labelsep  %                   + \labelsep
   \advance\leftmargin\enumindent%                   + \enumindent
   \relax
   \usecounter{\@pfenumctr}%       %    counter @pfenumctr
 }%
 }{%                             % END
   \end{pflist}%                   %   End list
   \addtocounter{pfenumdepth}{-1}%   pfenumdepth := pfenumdepth - 1
 }

\def\makefcn#1#2#3{{\let\or\relax
     \gdef\fcn@temp{}%
     \gdef\fcn@tempb{#2}%
     \@tfor\foo:=#3\do{\xdef\fcn@temp{\fcn@temp\or\foo}\xdef\fcn@tempb{\foo}}%
     \let\ifcase\relax\let\else\relax\let\fi\relax
     \xdef#1##1{\ifcase ##1 #2\fcn@temp\else\fcn@tempb\fi}}}

\newcommand{\pfmathdef}[1]{\relax\ifmmode #1\else $#1$\fi}

\newenvironment{pproof}{\begin{prooof}\pf\ }{\end{prooof}}
\newenvironment{Proof}[1][\proofname]{\vspace{0pt}\begin{proof}[#1]\refstepcounter{branch}\begin{prooof}}%
  {~\qed\end{prooof}\renewcommand{\qed}{}\end{proof}} 

\makeatother

\pfshortnumbers{4} % Use short step numbers for all levels >= 3
\let\kwfont\itshape % cheatpf defines it as \scshape

\usepackage
[backref,hyperindex,colorlinks,citecolor=blue,pdftex,pdfpagelabels,debug]{hyperref}
\usepackage[all]{hypcap} % After hyperref, to anchor floats correctly.

\setlrmarginsandblock{*}{1.5in}{1}
\checkandfixthelayout

\counterwithout{section}{chapter}
\counterwithout{figure}{chapter}
\firmlists

\abstractrunin
\newcommand{\mathboldit}[1]{\boldsymbol{#1}}
\newcommand{\myTimes}{\usepackage{mathptmx}
\DeclareSymbolFont{letters}{OML}{txmi}{m}{it}
\DeclareMathAlphabet{\mathbit}{OT1}{ptm}{bx}{it}
\renewcommand{\mathboldit}[1]{\mathbit{##1}}
}

\myTimes

\makeatletter
\renewcommand \thefigure {\@arabic\c@figure}

 \newcommand{\df}[1]{{\rmfamily\itshape\mdseries#1}}

\numberwithin{equation}{section} % in amsmath

\newcommand{\customqed}[1]{{\renewcommand{\qedsymbol}{#1}\qed}}
\newcommand{\varqed}{\customqed{\hbox{$\lrcorner$}}}

 \theoremstyle{plain}

 \newtheorem{lemma}{Lemma}[section]

 \newtheorem{theorem}{Theorem}

 \newtheorem{corollary}[lemma]{Corollary}

 \theoremstyle{definition}

 \newtheorem{Definition}{Definition}[section]
 \newenvironment{definition}{%
   \begin{Definition}}{\varqed\end{Definition}}

 \newtheorem{Defstep}{Step}

 \newenvironment{defstep}{%
   \begin{Defstep}}{\varqed\end{Defstep}}

 \theoremstyle{remark}

 \newtheorem{Condition}[lemma]{Condition}

 \newenvironment{condition}{%
   \begin{Condition}}{\varqed\end{Condition}}

 \newtheorem{Example}[lemma]{Example}
 \newtheorem{Examples}[lemma]{Examples}
 
 \newenvironment{example}{%
   \begin{Example}}{\varqed\end{Example}}

 \newtheorem{Remark}[lemma]{Remark}
 \newtheorem{Remarks}[lemma]{Remarks}
 
 \newenvironment{remark}{%
   \begin{Remark}}{\varqed\end{Remark}}
 \newenvironment{remarks}{%
   \begin{Remarks}}{\varqed\end{Remarks}}

\newenvironment{alphenum} % needs package enumerate, or memoir.
 {\begin{enumerate}[{\upshape (a)}]}{\end{enumerate}}

 \newcommand \ag{\alpha}

 \newcommand \eps{\varepsilon}

 \renewcommand \lg{\lambda} % it conflicts with lg (this is not serious).

 \newcommand{\cA}{\mathcal{A}}
 \newcommand{\cB}{\mathcal{B}}
 \newcommand{\cC}{\mathcal{C}}
 
 \newcommand{\cE}{\mathcal{E}}
 \newcommand{\cF}{\mathcal{F}}
 \newcommand{\cG}{\mathcal{G}}

 \newcommand{\cL}{\mathcal{L}}
 \newcommand{\cM}{\mathcal{M}}

 \newcommand{\cQ}{\mathcal{Q}}
 
 \newcommand{\cS}{\mathcal{S}}
 \newcommand{\cT}{\mathcal{T}}
 
 \newcommand{\cV}{\mathcal{V}}
 \newcommand{\cW}{\mathcal{W}}

 \newcommand{\bU}{\mathbf{U}}

 \newcommand{\bbM}{\mathbb{M}}

 \newcommand{\bbR}{\mathbb{R}}
 \newcommand{\bbZ}{\mathbb{Z}}

 \newcommand\set[1]{\mathopen\{#1\mathclose\}}
 \newcommand\setof[1]{\mathopen\{\,#1\,\mathclose\}}

 \newcommand\ang[1]{{\langle #1\rangle}}
 \newcommand\bigparen[1]{{\bigl(\,#1\,\bigr)}}

 \newcommand\Paren[1]{{\left( #1\right)}}
 
 \newcommand\cei[1]{{\lceil #1\rceil}}
 \newcommand\flo[1]{{\lfloor #1\rfloor}}
\newcommand\Cei[1]{{\left\lceil #1\right\rceil}}
 \newcommand\Flo[1]{{\left\lfloor #1\right\rfloor}}
 \newcommand\Pbof[1]{\mathbf{P}\mskip 1mu\mathopen[\,#1\,\mathclose]}
 \newcommand\Prob{\mathbf{P}}

  \renewcommand{\le}{\leqslant}
  \renewcommand{\ge}{\geqslant}
 
 \newcommand {\sbsq}{\subseteq}
 
 \newcommand {\spsq}{\supseteq}

 \newcommand{\txt}[1]{\text{\rmfamily\mdseries\upshape{#1}}}

\makeatother

 \newcommand{\pair}[2]{\mathopen(#1,#2\mathclose)}
 \newcommand{\tup}[1]{(#1)}

 \newcommand{\clint}[2]{\mathopen[#1,#2\mathclose]}
\newcommand{\lint}[2]{\mathopen[#1, #2\mathclose)}
\newcommand{\rint}[2]{(#1,#2]}

 \newcommand{\aux}{c}

 \renewcommand{\d}{d} % what was it?
 \newcommand{\D}{D}
 \newcommand{\f}{f}
 \newcommand{\fxp}{\varphi}
 \newcommand{\g}{g}
 \newcommand{\gxp}{\gamma}
 \newcommand{\h}{h}
 
 \newcommand{\tubh}{\overline w}
 \newcommand{\hxp}{\chi}
 \renewcommand{\l}{l}
 \newcommand{\m}{m}
 \newcommand{\p}{p}
 \newcommand{\ncln}{q}
 \newcommand{\pub}{{\overline p}}
 \renewcommand{\r}{r} % what was it?
 \newcommand{\s}{s}
 \newcommand{\tub}{w}
 \newcommand{\tubxp}{\omega}
 \newcommand{\txp}{\tau}
 \newcommand{\txpub}{\overline\tau}
 \newcommand{\R}{R}

 \newcommand{\T}{T}

 \newcommand{\Rect}{\text{\rmfamily\mdseries\upshape Rect}}

 \newcommand{\Wvalues}{\text{\rmfamily\mdseries\upshape Wvalues}}

 \newcommand{\slope}{\text{\rmfamily\mdseries\upshape slope}}
 \newcommand{\minslope}{\text{\rmfamily\mdseries\upshape minslope}}
 
 \newcommand{\Body}{\text{\rmfamily\mdseries\upshape Body}}

 \newcommand{\bub}{{\Delta}}
 \newcommand{\bubxp}{\delta}

 \newcommand{\slb}{{\sigma}}
 \newcommand{\slopeincr}{\Lambda}

\newenvironment{flushdescription}{\leftmargini=1em\begin{description}}
{\end{description}}

\begin{document}

  \title{Clairvoyant scheduling of random walks}

  \author{Peter G\'acs
\\Boston University
\\ gacs@bu.edu
}
\maketitle

\begin{abstract}
  Two infinite walks on the same finite graph are called \df{compatible} if
it is possible to introduce delays into them in such a way that they never
collide.  
Years ago, Peter Winkler asked the question: for which
graphs are two independent random walks compatible with positive
probability.  
Up to now, no such graphs were found.  
We show in this paper that large complete graphs have this property.  
The question is equivalent to a
certain dependent percolation with a power-law behavior: the probability
that the origin is blocked at distance $n$ but not closer decreases only
polynomially fast and not, as usual, exponentially.
 \end{abstract}

% \keywords{Percolation, random walk, renormalization.}
% \subjclass{82B43; 60K35; 68Q85}

 \section{Introduction}

 \subsection{The model}

Let us call any strictly increasing sequence $t = \tup{t(0)=0, t(1),\dotsc}$ of
integers a \df{delay sequence}.
For an infinite sequence $z=\tup{z(0), z(1),\dotsc}$, the 
delay sequence $t$ introduces a timing arrangement in which 
the value $z(n)$ occurs at time $t(n)$.
Given infinite sequences $z_{\d}$ and delay
sequences $t_{\d}$, for $\d = 0, 1$, let $a,b\in\{0,1\}$.
We say that there is a \df{collision}
at $\tup{a,n,k}$ if $t_{a}(n) \le t_{b}(k) < t_{a}(n+1)$ and
$z_{b}(k)=z_{a}(n)$, for $a\ne b$.
We call the two sequences $z_{0}, z_{1}$
\df{compatible} if there are delay sequences for them that avoid
collisions.

% Miert nem 0-tol?
For a finite undirected graph, a Markov chain
$Z(1),Z(2),\dots$ with values that are vertices in this graph
is called a \df{random walk} over this graph
if it moves, going from $Z(n)$ to $Z(n+1)$, 
from any vertex with equal probability to any one of its neighbors.

\begin{figure}
% \if\picsize\SlideSize
%     \psset{unit=0.5cm}
% \else\if\picsize\PaperSize
%     \psset{unit=0.8cm}
% \fi\fi
% \if\palette\DefaultColors
%     \definecolor{my-orchid}{cmyk}{0.00, 0.40, 0.02, 0.00}
% \else\if\palette\SlideColors
%     \definecolor{my-orchid}{cmyk}{0.00, 0.20, 0.01, 0.00}
% \else\if\palette\GrayColors
%     \definecolor{my-orchid} {gray}{0.8}
% \fi\fi\fi
 \[
  \includegraphics{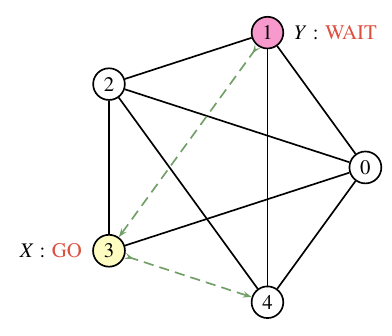}
 \]
 \caption{The clairvoyant demon problem.
$X,Y$ are ``tokens'' performing independent random walks on the same
graph: here the complete graph $K_{5}$.
A ``demon'' decides every time, whose turn it is.
She is clairvoyant and wants to prevent collision.
}
\label{fig:rwalk} 
 \end{figure}

Take two infinite random sequences $Z_{\d}$ for $\d=0,1$
independent from each other,
both of which are random walks on the same finite undirected graph.
Here, the delay sequence $t_{\d}$ can be viewed as causing the sequence
$Z_{\d}$ to stay in state $z_{\d}(n)$ between times $t_{\d}(n)$ and
$t_{\d}(n+1)$.
(See the example on the graph $K_{5}$ in Figure~\ref{fig:rwalk}.)
A collision occurs when the two delayed walks enter the
same point of the graph.
Our question is: are $Z_{0}$ and $Z_{1}$ compatible with positive
probability?
The question depends, of course, on the graph.
Up to the present paper, no graph was known with an affirmative
answer.
% Temporarily out:
% Miert 0-4?

Consider the case when the graph is the complete graph $K_{\m}$ of size
$\m$. 
It is known that if $\m \le 3$ then the two sequences are compatible
only with zero probability.
Simulations suggest strongly that the walks do not collide if $\m \ge 5$, and,
somewhat less strongly, even for $\m = 4$.
The present paper proves the following theorem.

 \begin{theorem}[Main]\label{thm:main}
If $\m$ is sufficiently large then on the graph $K_{\m}$, 
the independent random walks $Z_{0}$, $Z_{1}$ 
are compatible with positive probability.
 \end{theorem}

The upper bound computable for $\m$ from the proof is very bad.
In what follows we will also use the simpler notation
 \[
   X = Z_{0},\quad Y = Z_{1}.
 \]
The problem, called the clairvoyant demon problem, 
arose first in distributed computing.
The original problem was to find a leader among a finite number of
processes that form the nodes of a communication graph.
There is a proposed algorithm: at start, let each process have a ``token''.
The processes pass the tokens around in such a way that each token performs
a random walk.
However, when two tokens collide they merge.
Eventually, only one token will remain and whichever process has it becomes
the leader.
The paper~\cite{CopperTetWink93} examined the algorithm in the traditional
setting of distributed computing, when the timing of this procedure is
controlled by an adversary.
Under the (reasonable) assumption that
the adversary does not see the future sequence of
moves to be made by the tokens, the work~\cite{CopperTetWink93} 
gave a very good upper bound on the expected time by which a leader will be found.
It considered then the question whether a clairvoyant adversary (a
``demon'' who sees far ahead into the future token moves) can, 
by controlling the timing alone, with positive probability,
prevent two distinct tokens from ever colliding.
The present paper solves Conjecture 3 of~\cite{CopperTetWink93},
which says that this is the case when the communication graph is a
large complete graph.

The proof is long and complex, but its construction is based on some simple
principles presented below in Section~\ref{sec:outline}, 
after first transforming the problem into a percolation problem in Section~\ref{subsec:perc}.
The rest of the paper is devoted mainly to proving Lemma~\ref{lem:main}, stated at the end
of Section~\ref{sec:outline}.
The main ideas can be summarized as follows.
\begin{enumerate}[1.]
\item Reformulate into a percolation problem in the upper quarter plane,
where closed sites $(i,j)$ are those with $X(i)=Y(j)$ (Section~\ref{subsec:perc}).
\item Let us call these closed sites \df{obstacles}.
They have long-range dependencies, making it hard to
handle them directly: therefore for ``bad events'' (like when too many obstacles
crowd together), we introduce
a hierarchy of new kinds of obstacles in Section~\ref{sec:mazery}.
These fall into two categories.
Rectangular \df{traps}, and infinite horizontal or
vertical \df{walls} of various widths.
A wall can only be penetrated at certain places.

\item Traps at the higher levels of the hierarchy are larger, walls are 
wider and denser.
But these higher-level objects have smaller probability of
occurrence, therefore on level $k$ it will be possible to assume the local absence of
level $k+1$ objects with high probability.
That percolation is possible under this assumption for each $k$, is proved in
section~\ref{sec:approx-proof}.

\item
The structure of random traps, walls, and so on,
whose versions appear on all levels, is defined in Section~\ref{sec:mazery}. 
The higher-level objects are defined by a recursive procedure in
Section~\ref{subsec:plan.constr}, though some parameters will only be fixed in 
Section~\ref{sec:params}.
Their combinatorial properties are proved in Section~\ref{subsec:plan.prop}.

\item Those probabilistic estimates for the higher-level objects 
that can be proved without fixing all parameters are given in
Sections~\ref{sec:bounds}.
The rest is proved in Section~\ref{sec:after-scale-up}.

\end{enumerate}

\subsection{Related synchronization problems}\label{sec:related}

Let us define a notion of collision somewhat different from the previous
section.
For two infinite 0-1-sequences $z_{\d}$
($\d = 0, 1$) and corresponding delay sequences $t_{\d}$ we say that there
is a \df{collision} at $(a,n)$ if $z_{a}(n) = 1$, and there is no $k$
such that $z_{b}(k)=0$ and $t_{a}(n)=t_{b}(k)$, for $a\ne b$.
We say that the sequences $z_{\d}$ are \df{compatible} if there is a pair
of delay sequences $t_{\d}$ without collisions.
It is easy to see that this is equivalent to saying that 0's can be
deleted from both sequences in such a way that the resulting sequences have
no collisions in the sense that they never have a 1 in the same position.

Suppose that for $\d = 0, 1$, $Z_{\d} = \tup{Z_{\d}(0), Z_{\d}(1),\dotsc}$
are two independent infinite sequences of independent random variables
where $Z_{\d}(j) = 1$ with probability $\p$ and 0 with probability $1-\p$.
Our question is: are $Z_{0}$ and $Z_{1}$ compatible with positive
probability?
The question depends, of course, on the value of $\p$: intuitively, it
seems that they are compatible if $\p$ is small.

Peter Winkler and Harry Kesten~\cite{KestenWinkler}, 
independently of each other, found an upper
bound smaller than 1/2 on the values $\p$ for which $Z_{0}, Z_{1}$ are
compatible.
Computer simulations by John Tromp
suggest that when $\p < 0.3$, with positive probability the sequences are
compatible.
The paper~\cite{GacsChat04} proves  that if
$\p$ is sufficiently small then with positive probability, $Z_{0}$ and
$Z_{1}$ are compatible.
The threshold for $\p$ obtained from the proof is only $10^{-400}$, so
there is lots of room for improvement between this number and the
experimental $0.3$.

The author recommends the reader to refer to paper~\cite{GacsChat04} while
reading the present one.
Its technique prefigures the one used here:
the main architecture is similar, but many details are simpler.

\subsection{A percolation}\label{subsec:perc}

The clairvoyant demon
problem has a natural translation into a percolation problem.
Consider the lattice $\bbZ_{+}^{2}$, and a directed graph obtained from it
in which each point is connected to its right and upper neighbor.
For each $i,j$, let us ``color'' the $i$th vertical line 
by the state $X(i)$, and the $j$th horizontal line by the state
$Y(j)$.
The ingoing edges of a point $\pair{i}{j}$ will be deleted from the graph if
$X(i)=Y(j)$, if its horizontal and vertical colors coincide.
We will also say that point $\pair{i}{j}$ is \df{closed}; otherwise, it will be
called \df{open}.
(It is convenient to still keep the closed point $\pair{i}{j}$ in the graph,
even though it became unreachable from the origin.)
The question is whether, with positive probability, an infinite oriented
path starting from $\pair{0}{0}$ exists in the remaining random graph
 \[
   \cG = (\cV, \cE).
 \]
In~\cite{GacsChat04}, we proposed
to call this sort of percolation, where two infinite random
sequences $X,Y$ are given on the two coordinate axes and the openness of a
point or edge at position $\pair{i}{j}$ depends on the pair $\tup{X(i),Y(j)}$, a
\df{Winkler percolation}.
 \begin{figure}
\[
\includegraphics{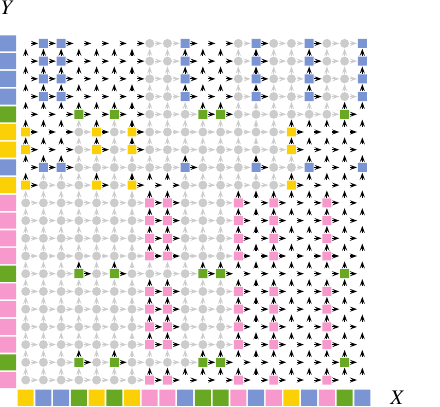}
\]
\[
\mbox{
\includegraphics{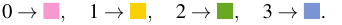}
}
\]
 \caption{
Percolation for the clairvoyant demon problem, for random
walks on the complete graph $K_{4}$.
Round light-grey dots mark the reachable points.}
 \end{figure}
This problem permits an interesting variation: undirected percolation,
where the the whole lattice $\bbZ^{2}$ is present, and the edges are
undirected. 
This variation has been solved, independently, in~\cite{WinklerCperc99} and
\cite{BalBollobStacCperc99}.
On the other hand, the paper~\cite{GacsClairv99} shows that the
directed problem has a different nature, since 
it has power-law convergence (the undirected percolations have
the usual exponential convergence).

 \section{Outline of the proof}\label{sec:outline}
The proof introduces a whole army of auxiliary concepts, which are cumbersome to
keep track of.
The reader will find it helpful occasionally to refer to 
the glossary and notation index provided at the end of the paper.

 \subsection{Renormalization}\label{subsec:renorm}
The proof method used is \df{renormalization} (scale-up).
An example of the
ordinary renormalization method would be when, say, in an
Ising model, the space is partitioned into blocks, spins in
each block are summed into a sort of ``superspin'', and it is shown that
the system of super-spins exhibits a behavior that is in some sense similar
to the original system.
We will also map our model
repeatedly into a series of higher-order models similar to each other.
However the definition of the new models is more complex than just
taking the sums of some quantity over blocks.
The model which will scale up properly may contain a number of new objects,
and restrictions more combinatorial than computational in character.

The method is messy, laborious, and rather crude (rarely leading
to the computation of exact constants).
However, it is robust and well-suited to ``error-correction'' situations.
Here is a rough first outline.

 \begin{enumerate}[\upshape 1.]

  \item
Fix an appropriate sequence 
$\bub_{1} < \bub_{2} < \dotsb$, of scale parameters for which
$\bub_{k+1} > 3\bub_{k}$ holds\footnote{In the present
  paper, the actual quotient between these
  parameters will be more than $10^{7}$.}.
(The actual values of our
parameters will only be fixed in Section~\ref{sec:params}.) 
Let $\cF_{k}$ be the event that point $\pair{0}{0}$ is blocked in the square 
$\clint{0}{\bub_{k}}^{2}$.
(In other applications, it could be some other \df{ultimate bad event}.)
Throughout the proof, we will denote the probability of an event $E$ by 
$\Prob(E)$.
We want to prove 
 \[
   \Prob\bigparen{\bigcup_{k} \cF_{k}} < 1.
 \]
This will be sufficient: if $\pair{0}{0}$ is not blocked in any
finite square then by compactness (or by what is sometimes called K\"onig's
Lemma), there is an infinite oriented path starting at $\pair{0}{0}$.

  \item\label{i:Fprime}
Identify some events that we call
\df{bad events} and some others called \df{very bad events},
where the latter are much less probable.
 
Define a series $\cM^{1},\cM^{2},\dotsc$ of models
similar to each other (in our case each based on the same directed lattice graph on
$\bbZ_{+}^{2}$)
where the very bad events of $\cM^{k}$ become
the bad events of $\cM^{k+1}$.
(The structure of these
models will be introduced in Section~\ref{sec:mazery}, the actual
recursive definition is carried out in Section~\ref{subsec:plan.constr}, with the
fixing of parameters left to Section~\ref{sec:params}.)

Let $\cF'_{k}$
hold iff some bad event of $\cM^{k}$ happens in the square
$\clint{0}{\bub_{k+1}}^{2}$.
  \item
Prove
 \begin{equation}\label{eq:bd-non-reach-by-wall-all}
  \cF_{k} \sbsq \bigcup_{i \le k} \cF'_{i}.
 \end{equation}
(This will follow from the structure of definitions in
Section~\ref{subsec:plan.constr}.)

  \item
Prove $\sum_{k} \Prob(\cF'_{k}) < 1$.
(This results from the estimates in Section~\ref{sec:params}.)
  
  \end{enumerate}
In later discussions, we will
frequently delete the index $k$ from $\cM^{k}$ as well as from other
quantities defined for $\cM^{k}$.
In this context, we will refer to $\cM^{k+1}$ as $\cM^{*}$.

 \subsection{Application to our case}\label{subsec:application}

The role of the
``bad events'' of Subsection~\ref{subsec:renorm} will be played by \df{traps}
and \df{walls}.
\begin{figure}
  \centering
 \includegraphics[scale=0.3]{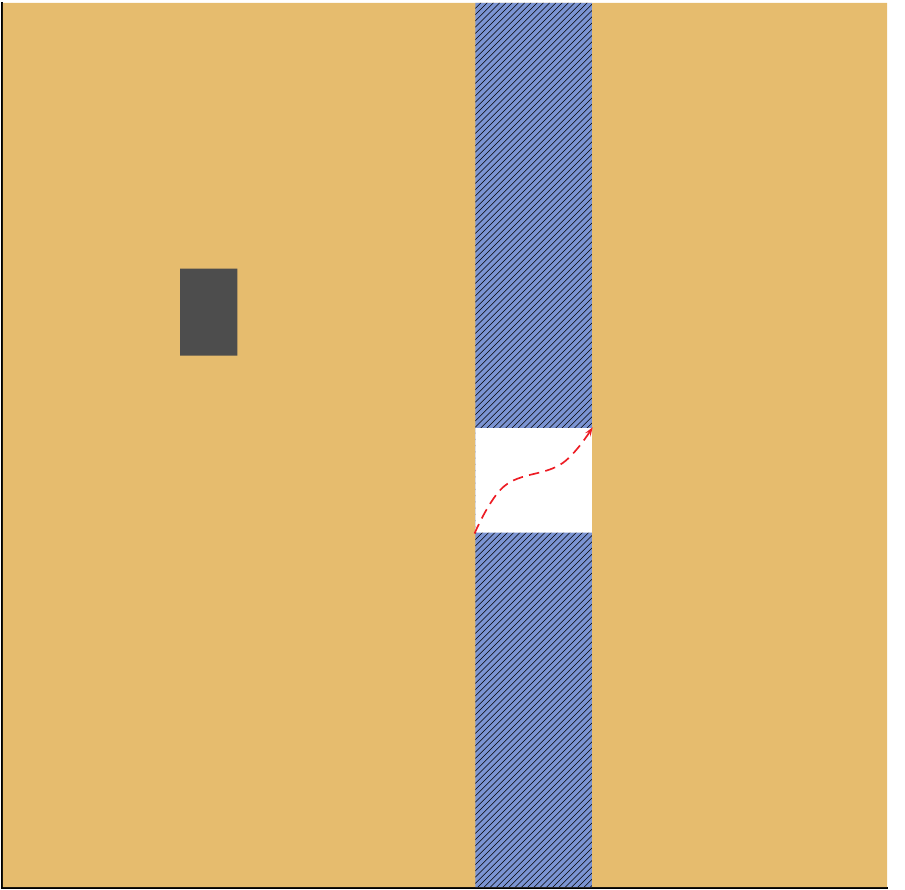}  
  \caption{A trap and a wall with a hole.}
  \label{fig:trap-wall}
\end{figure}
The simplest kind of trap is a point $\pair{i}{j}$ in the plane such that
$X(i) = Y(j)$; in other words, a closed point.
More generally, traps will be certain rectangles in the plane. 
We want to view the occurrence of two traps close to each other as a very
bad event; however, this is justified only if this is indeed very
improbable.
Consider the events 
 \[
       \cA_{5} = \set{X(1)=X(2)=X(3)=Y(5)},
\quad \cA_{13} = \set{X(1)=X(2)=X(3)=Y(13)}.
 \]
(For simplicity, this example assumes that the random walk has the option
of staying at the same point, that is loops have been added to the
graph $K_{m}$.)
The event $\cA_{5}$ makes the rectangle
$\clint{1}{3}\times \{5\}$ a trap of size 3, and has
probability $\m^{-3}$.
Similarly for the event $\cA_{13}$ and the rectangle
$\clint{1}{3}\times \{13\}$.
However, these two events are not independent: 
the probability of $\cA_{5} \cap \cA_{13}$ is only $\m^{-4}$, not
$\m^{-6}$.
The reason is that the event $\cE=\set{X(1)=X(2)=X(3)}$ 
significantly increases the conditional
probability that, say, the rectangle $\clint{1}{3}\times \{5\}$ becomes a
trap.
In such a case, we will want to say that event $\cE$ creates a
\df{vertical wall} on the segment $\rint{0}{3}$.

Though our study only concerns the integer lattice, it is convenient to use
the notations of the real line and Euclidean plane.
In particular, walls will be right-closed intervals.
(Even though, of course,
$\rint{a}{b}\cap\bbZ =\clint{a+1}{b}\cap\bbZ$, but we will 
not consider the interval $\rint{4}{9}$ to be contained in $\clint{5}{10}$.)
We will say that a certain rectangle \df{contains} a wall if the
corresponding projection contains it, and that the same
rectangle \df{intersects} a wall if the
corresponding projection intersects it.

Traps will have low probability.
If there are not too many traps, it is possible to get around them.
On the other hand, to get through walls, one also needs extra luck:
such lucky events will be called \df{holes} (see Figure~\ref{fig:trap-wall}).
Our proof systematizes the above ideas by
introducing an abstract notion of traps, walls and holes.
We will have walls of many different types.
To each (say, vertical) wall of a given type, the probability that a
(horizontal) hole goes through it
at a given point will be much higher than the probability that a (horizontal)
wall of this type occurred at that point.
Thus, the ``luck'' needed to go through some wall type is still smaller
than the ``bad luck'' needed to create a wall of this type.

This model will be called a \df{mazery} $\cM$ (a system for creating
mazes). 
In any mazery, whenever it happens that
walls and traps are well separated from each other and
holes are not missing, then paths can pass through.
(Formally, this claim will be called Lemma~\ref{lem:approx} (Approximation)---as
the main combinatorial tool in a sequence of successive approximations.)
Sometimes, however, unlucky events arise.
These unlucky events can be classified in the types listed below.
For any mazery $\cM$, we will define a mazery $\cM^{*}$ whose walls
and traps correspond (essentially) to these typical unlucky events.
 \begin{enumerate}[--]

  \item
A minimal rectangle enclosing two traps very close to each other, both of
whose projections are disjoint, is an \df{uncorrelated compound trap}
(see Figure~\ref{fig:compound-traps}).

  \item
For both directions $\d = 0, 1$, a (essentially) 
minimal rectangle enclosing 4
traps very close to each other, whose $\d$ projections are disjoint,
is a \df{correlated compound trap} (see Figure~\ref{fig:compound-traps}).

\begin{figure}
  \centering
\includegraphics[scale=0.6]{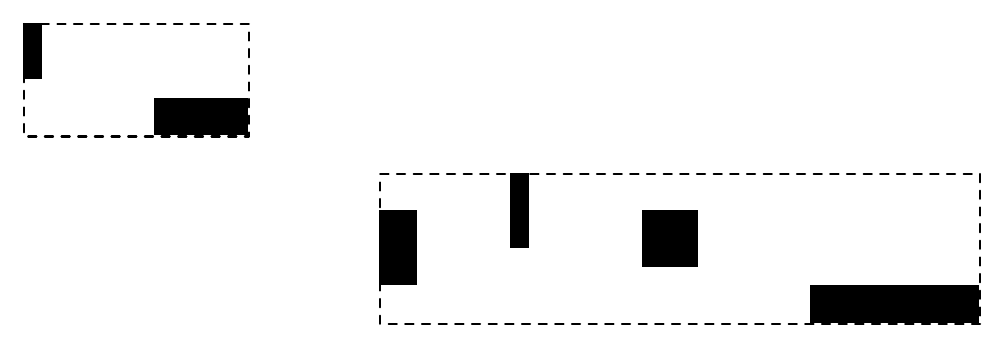}  
  \caption{An uncorrelated and a correlated compound trap}
  \label{fig:compound-traps}
\end{figure}

  \item
Whenever a certain horizontal wall $W$ appears and at the same
time there is a large interval without a vertical hole of $\cM$ 
through $W$, this situation gives rise to a
trap of $\cM^{*}$ of the \df{missing-hole} kind
(see Figure~\ref{fig:no-hole}). 

\begin{figure}
  \centering
  \includegraphics[scale=0.6]{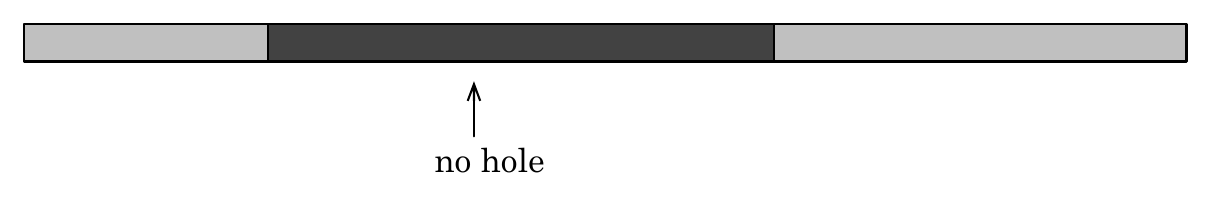}
  \caption{A missing-hole trap}
  \label{fig:no-hole}
\end{figure}

  \item
A pair of very close walls of $\cM$ gives rise to a wall of
$\cM^{*}$ called a \df{compound wall} (see Figure~\ref{fig:compound-wall}).

\begin{figure}
  \centering
  \includegraphics[scale=0.4]{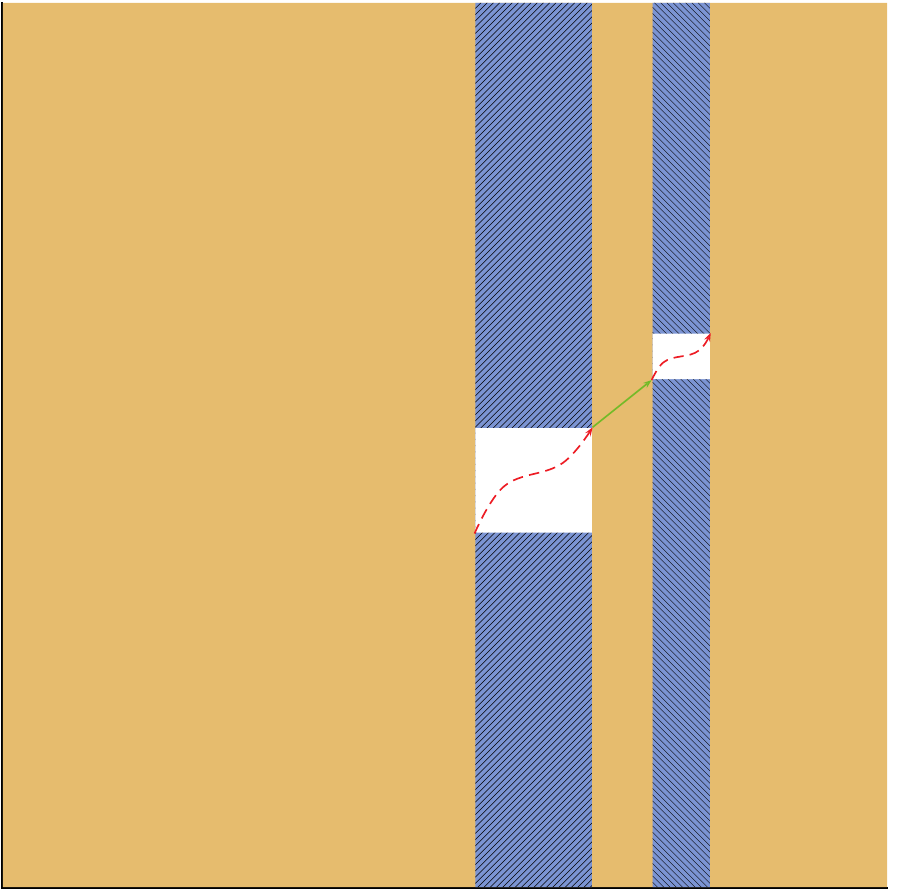}
  \caption{A compound wall, with a hole through it.}
  \label{fig:compound-wall}
\end{figure}

  \item
A segment of the $X$ or $Y$ sequence such that conditioning on it,
a correlated trap or a trap of the missing-hole kind
occurs with too high conditional probability,
is a new kind of wall called an \df{emerging} wall.
(These are the walls that, indirectly, give rise to all other walls.)

 \end{enumerate}
(The exact definition of these objects involves some extra technical
conditions: here, we are just trying to give the general idea.)

At this point, it would be hard for the reader to appreciate that the set of
kinds of objects (emerging traps and walls) is \emph{complete}: that
there is no need in any other ones.
An informal effort to try to prove Lemma~\ref{lem:approx} (Approximation)
should give some
insight: in other words, the reader should try to convince herself that
as long as the kind of ``very bad events'' covered by the emerging traps and walls do
not occur, percolation can occur.
 
There will be a constant 
 \begin{equation}\label{eq:hxp}
   \hxp = 0.015.
 \end{equation}
with the property that if a wall has probability $\p$ then
a hole getting through it has probability lower bound $\p^{\hxp}$.
Thus, the ``bad events'' of the outline in
Subsection~\ref{subsec:renorm} are the traps and walls of $\cM$, 
the ``very bad events'' are (modulo some details that are not important
now) the new traps  and walls of $\cM^{*}$.
Let $\cF, \cF'$ be the events $\cF_{k},\cF'_{k}$ formulated in
Subsection~\ref{subsec:renorm}.
Thus, $\cF'$ says that in $\cM$ a wall or a
trap is contained in the square $\clint{0}{\bub^{*}}^{2}$.

\begin{remark}
  The paper uses a number of constants:
 \begin{align*}
\aux_{1},\aux_{2},\aux_{3},H,\lg,\slopeincr,\R_{0},\gxp,\hxp,\txp,\txp',\tubxp.
 \end{align*}
It would be more confusing to read (and a nightmare to debug)
if we substituted numerical values for them everywhere.
This is debatable in some cases, like $\lg=\sqrt 2$, $\tubxp=4.5$, but 
even here, the symbol $\lg$ emphasizes that there is no
compelling reason for using the exact value $\sqrt 2$.
There is a notation index allowing the reader to look up the definition of each
constant, whenever needed.
\end{remark}

We do not want to see all the details of $\cM$ once we are
on the level of $\cM^{*}$: this was the reason for creating $\cM^{*}$ in
the first place.
The walls and traps of $\cM$ will indeed become transparent; however,
some restrictions will be inherited from them: these are distilled in the
concepts of a \df{clean point} and of a \df{slope constraint}.
Actually, we distinguish the concept of \df{lower left clean} and
\df{upper right clean}.
Let 
 \[
   \cQ
 \]
be the event that point $\pair{0}{0}$ is not upper right clean in $\cM$.

We would like to say that in a mazery, if points $\pair{u_{0}}{u_{1}}$,
$\pair{v_{0}}{v_{1}}$ are such that for $\d=0,1$ we have
$u_{\d} < v_{\d}$ and there are no walls and traps
in the rectangle $\clint{u_{0}}{v_{0}} \times\clint{u_{1}}{v_{1}}$,
then $\pair{v_{0}}{v_{1}}$  is reachable from $\pair{u_{0}}{u_{1}}$.
However, this will only hold with some restrictions.
What we will have is the following, with an appropriate parameter
 \[
   0 \le \slb < 0.5.
 \]
 \begin{condition}\label{cond:reachable-hop}
Suppose that points $u = \pair{u_{0}}{u_{1}}$, $v = \pair{v_{0}}{v_{1}}$ are 
such that for $\d=0,1$ we have
$u_{\d} < v_{\d}$ and there are no traps or walls contained
in the rectangle between $u$ and $v$.
If $u$ is upper right clean, $v$ is lower left clean and
these points also satisfy the slope-constraint
 \begin{equation*}
  \slb \le \frac{v_{1}-u_{1}}{v_{0}-u_{0}} \le 1/\slb
 \end{equation*}
then $v$  is reachable from $u$.
 \end{condition}

We will also need sufficiently many clean points:

  \begin{condition}\label{cond:dense}
For every square $\pair{a}{b}+\rint{0}{3\bub}^{2}$ that does not contain walls or
traps, there is a lower left clean point in its middle third
$\pair{a}{b}+\rint{\bub}{2\bub}^{2}$.
  \end{condition}

 \begin{lemma}\label{lem:imply-reachable}
We have $\cF \sbsq \cF' \cup \cQ$.
 \end{lemma}
 \begin{proof}
(Please, refer to Figure~\ref{fig:main-lemma}.) 
Suppose that $\cQ$ does not hold, then $\pair{0}{0}$ is upper right clean.

Suppose also that $\cF'$ does not hold: then by Condition~\ref{cond:dense}, 
there is a point $u = \pair{u_{0}}{u_{1}}$ in the square
$\clint{\bub}{2\bub}^{2}$ that is lower left clean in $\cM$. 
This $u$ also satisfies the slope condition
$1/2 \le u_{1}/u_{0} \le 2$ and is hence, by
Condition~\ref{cond:reachable-hop}, reachable from $\pair{0}{0}$.
 \end{proof}
 
We will define a sequence of mazeries $\cM^{1},\cM^{2},\dotsc$ with
$\cM^{k+1}=(\cM^{k})^{*}$, with $\bub_{k} \to \infty$.
All these mazeries are on a common probability space, since $\cM^{k+1}$ is
a function of $\cM^{k}$.
All components of the mazeries will be indexed correspondingly: for 
example, the event $\cQ_{k}$ that $\pair{0}{0}$ is not upper right clean in
$\cM^{k}$ plays the role of $\cQ$ for the mazery$\cM^{k}$.
We will have the following property:

\begin{condition}\label{cond:Qbd} We have
 $\cQ_{k} \sbsq \bigcup_{i<k} \cF'_{i}$.
 \end{condition}
In other words, if there are no traps or walls of any order $i<k$
in the $\bub_{i+1}$-neighborhood of the origin, then the origin is upper-right
clean on level $k$.

This, along with Lemma~\ref{lem:imply-reachable} implies
$\cF_{k} \sbsq \bigcup_{i \le k} \cF'_{i}$, which is
inequality~\eqref{eq:bd-non-reach-by-wall-all}.
Hence the theorem is implied by the following lemma, which will be
proved after all the details are given:
 \begin{figure}
   \begin{center}
  \includegraphics{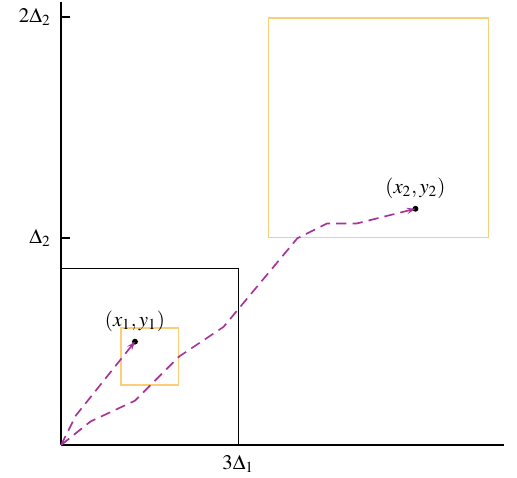}
  \caption{Illustration of  
Lemma~\protect\ref{lem:imply-reachable} and Theorem~\protect\ref{thm:main}}
\label{fig:main-lemma}
   \end{center}
 \end{figure}
% Miert (x,y) az abraban?

For the following lemma, recall that $\m$ is the
the number of ``colors'' in the percolation setting, and therefore $1/m$
upper-bounds 
the conditional probability of a trap $\{X(i)=Y(j)\}$ (under, for
example, fixing $X(i)$).

 \begin{lemma}[Main]\label{lem:main}
If $\m$ is sufficiently large then the sequence $\cM^{k}$ can be
constructed, in such a way that it satisfies all the above conditions and
also 
 \begin{equation}\label{eq:main}
   \sum_{k} \Prob(\cF'_{k}) < 1.
 \end{equation}
 \end{lemma}

\subsection{The rest of the paper}

The proof structure is quite similar to~\cite{GacsChat04}.
That paper is not simple, but it is still simpler than the present one, and
we recommend very much looking at it in order to see some of the ideas
going into the present paper in their simpler, original setting.
Walls and holes, the general form of the definition of a mazery and the
scale-up operation are similar.
There are, of course, differences: traps are new.

Section~\ref{sec:mazery} defines the random structures called mazeries.
This is probably the hardest to absorb, since the structure has a large number
of ingredients and required properties called Conditions.
The proof in the sections that follow will clarify, however, the role of each
concept and condition.

Section~\ref{sec:plan} defines the scale-up operation $\cM^{k}\mapsto\cM^{k+1}$.
It also proves that scale-up preserves almost all \emph{combinatorial}
properties, that is those that do not involve probability bounds.
The one exception is the reachability property, formulated by
Lemma~\ref{lem:approx} (Approximation): its proof is more complex, and
is postponed to Section~\ref{sec:approx-proof}.

Section~\ref{sec:bounds} estimates how the probability bounds are transformed by
the scale-up operation.
Section~\ref{sec:params} specifies the parameters in such a way that guarantees
that the probability conditions are also preserved by scale-up.
Section~\ref{sec:after-scale-up} carries out the computations leading to the
proof of those conditions.

Section~\ref{sec:main-proof}
ties up all threads into the proof of Lemma~\ref{lem:main} and the
proof of the theorem.

\section{Mazeries}\label{sec:mazery}

This section consists almost exclusively of definitions.

 \subsection{Notation}\label{subsec:notation}

The notation $\pair{a}{b}$ for real numbers $a,b$ will generally mean for us
the pair, and not the open interval.
Occasional exceptions would be pointed out separately.
We will use
 \[
   a \land b = \min(a,b),\quad a \lor b = \max(a,b).
 \]
As mentioned earlier, we will use intervals on the real line 
and rectangles over the Euclidean plane, even though we are really only
interested in the lattice $\bbZ_{+}^{2}$.
To capture all of $\bbZ_{+}$ this way, 
for our right-closed intervals $\rint{a}{b}$, we allow the
left end $a$ to range over all the values $-1,0,1,2,\dots$.
For an interval $I=\rint{a}{b}$, we will denote 
 \[
   X(I) = \tup{X(a+1),\dots,X(b)}.
 \]
The \df{size} of an interval $I$ with endpoints $a,b$ (whether it is open,
closed or half-closed), is denoted by $|I| = b-a$.
By the \df{distance} of two points 
$a = \pair{a_{0}}{a_{1}}$, $b=\pair{b_{0}}{b_{1}}$ of the plane, we mean 
 \[
   |b_{0}-a_{0}| \lor |b_{1}-a_{1}|.
 \]
The \df{size} of a rectangle
 \[
   \Rect(a, b) 
    = \clint{a_{0}}{b_{0}} \times \clint{a_{1}}{b_{1}}
 \]
in the plane is defined to be equal to the distance between $a$ and $b$.
For two different points $u = \pair{u_{0}}{u_{1}}$, $v = \pair{v_{0}}{v_{1}}$
in the plane, when $u_{0} \le v_{0}$, $u_{1} \le v_{1}$:
 \begin{alignat*}{2}
      &\slope(u, v) &&= \frac{v_{1}-u_{1}}{v_{0}-u_{0}}, 
\\ &\minslope(u, v) &&= 
   \min\bigparen{\slope(u, v), 1/\slope(u, v)}.
  \end{alignat*}
We introduce the following partially open rectangles
 \begin{equation}\label{eq:rectangles}
 \begin{alignedat}{2}
   &\Rect^{\rightarrow}(a, b)
                      &&= \rint{a_{0}}{b_{0}} \times \clint{a_{1}}{b_{1}},
\\ &\Rect^{\uparrow}(a, b)
                      &&= \clint{a_{0}}{b_{0}} \times \rint{a_{1}}{b_{1}}.
 \end{alignedat}
 \end{equation}
The relation
 \[
   u \leadsto v
 \]
says that point $v$ is reachable from point $u$ (the underlying graph
will always be clear from the context).
For two sets $A, B$ in the plane or on the line,
 \[
   A + B = \setof{a + b : a \in A,\; b \in B}.
 \]

\subsection{The structure}

\subsubsection{The tuple}

All our structures defined below refer to ``percolations'' over the same lattice
graph, in $\bbZ_{+}^{2}$,
defined by the pair of sequences of random variables
 \[
  Z = \tup{X, Y} = \tup{Z_{0}, Z_{1}},
 \]
where
$Z_{\d}=\tup{Z_{\d}(0),Z_{\d}(1),\dots}$ with $Z_{\d}(t)\in \{1,\dots,m\}$
independent random walks on the
set $\{1,\dots,\m\}$ of nodes of the graph $K_{\m}$ for some fixed $\m$.

A \df{mazery}
 \begin{equation}\label{eq:mazery}
 \bbM = (\cM, \bub, \slb, \R, \tub, \ncln) 
 \end{equation}
consists of a random process $\cM$, the parameters $\bub>0$, $\slb\ge 0$, 
$\R>0$, and the probability bounds $\tub>0$, $\ncln$,
all of which will be detailed below, along with conditions that 
they must satisfy.
Let us describe the random process
 \[
   \cM  = (Z, \cT, \cW, \cB, \cC, \cS).
 \]
We have the random objects
 \begin{equation*}
 \cT, \quad
 \cW = \tup{\cW_{0},\cW_{1}},\quad
 \cB = \tup{\cB_{0},\cB_{1}},\quad
 \cC = \tup{\cC_{0}, \cC_{1}},\quad
 \cS = \tup{\cS_{0}, \cS_{1}, \cS_{2}}.
 \end{equation*}
all of which are functions of $Z$.
The set $\cT$ of random \df{traps}
is a set of some closed rectangles of size $\le \bub$.
For trap $\Rect(a, b)$, we will say that it
\df{starts} at its lower left corner $a$.

 \begin{definition}[Wall values]\label{def:walls}
To describe the process $\cW$, we introduce the concept of a
 \df{wall value} $E = (B, \r)$.  % Jeloles?
Here $B$ is the \df{body} which is a right-closed interval,\footnote{This
is different from the definition in
the paper~\cite{GacsChat04}, where walls were open
intervals.} and \df{rank} 
 \[
  \r \ge \R > 0.
 \]
We write $\Body(E) = B$, $|E| = |B|$.
We will sometimes denote the body also by $E$.
Let $\Wvalues$ denote the set of all possible wall values.
 \end{definition}

Walls will arise in a variety of ways, but the properties 
we are interested in will only depend on the body and the rank.
Walls of higher rank have a smaller probability upper bound of 
occurrence, and smaller probability lower bound of holes through them.
They arise at higher levels in the hierarchical construction,
but the rank depends on more details of how the wall arose than just the level
of the mazery it is in.

Let
 \begin{align*}
  \bbZ_{+}^{(2)}
 \end{align*}
denote the set of pairs $\pair{u}{v}$ with $u<v$, $u,v\in\bbZ_{+}$.
The random objects
 \begin{equation*}
  \begin{split}
     \cW_{\d}\sbsq\cB_{\d}  &\sbsq \Wvalues,
\\   \cS_{\d}\subseteq\cC_{\d} &\sbsq \bbZ_{+}^{(2)} \times \{-1, 1\}
                                           \txt{ for } \d = 0, 1,
\\   \cS_{2} &\sbsq \bbZ_{+}^{(2)}\times\bbZ_{+}^{(2)} \times \{-1, 1\} \times \{0,1,2\}
  \end{split}
 \end{equation*}
are also functions of $Z$.
(Note that we do not have any $\cC_{2}$.)

 \begin{definition}[Barriers and walls]\label{def:barrier}
The elements of $\cW_{\d}$ and $\cB_{\d}$ are called \df{walls} and
\df{barriers} of $Z_{\d}$ respectively, where the sets
$\cW_{\d},\cB_{\d}$ are functions of $Z_{\d}$.
In particular, elements of $\cW_{0}$ are called \df{vertical walls}, and
elements of $\cW_{1}$ are called \df{horizontal walls}.
Similarly for barriers.
When we say that a certain interval contains a wall or barrier we mean that
it contains its body.

A right-closed interval is called \df{external} if it intersects no walls.
A wall is called \df{dominant} if it is
surrounded by external intervals each of which is either of size $\ge\bub$ or is
at the beginning of $\bbZ_{+}$.
Note that if a wall is dominant then it contains every wall intersecting it.
 \end{definition}

The set of barriers is a random subset of the set of all possible
wall values, and the set of walls is a random subset of the set of barriers.
In the scale-up operation, we first define barriers, and then
we select the walls from among them.
The form of the definition of barriers implies their simple dependency
properties required in Condition~\ref{cond:distr}.\ref{i:distr.indep}, which
then make simple probability upper bounds possible.
These then hold for walls as well, since walls are barriers.
On the other hand, walls have the nicer combinatorial properties
we need to prove eventually reachability (percolation).

The following definition uses the fact following from
Condition~\ref{cond:distr}.\ref{i:distr.indep.barrier} that whether an interval $B$
is a barrier of the process $X$ depends only $X(B)$.

\begin{definition}[Potential wall]\label{def:potential-wall}
For a vertical wall value $E=(B,\r)$ and a value of $X(B)$ making $E$ a
barrier of rank $\r$
we will say that $E$ is a \df{potential vertical wall} of rank $\r$
if there is an extension of $X(B)$ to a complete sequence $X$ that makes
$E$ a vertical wall of rank $\r$.
Similarly for horizontal walls.
\end{definition}

 \begin{remarks}
 \begin{enumerate}[\upshape 1.]

  \item
We will see below that,
for any rectangle with projections $I \times J$, the event that it is a
trap is a function of the pair $X(I), Y(J)$.
Also for any interval $I$, the event that it is a
(say, vertical) barrier depends only on $X(I)$, but the same is not true of
walls.

  \item In the definition of the mazery
$\cM^{k+1}$ from mazery $\cM^{k}$, we will drop \df{low rank} walls of
$\cM^{k}$, (those with $\le\R_{k+1}$).
These walls will have high probability
of holes through them, so reachability will be conserved.

To control the proliferation of walls, a pair of close
walls of $\cM^{k}$ will give rise to a \df{compound wall} of $\cM^{k+1}$
only if at least one of the components has low rank.

%   \item
% In what follows we will refer to $\cM$ by itself 
% also as a \df{mazery}, and will mention $\bbM$ only rarely.
% This should not cause confusion; though $\cM$ is a component of $\bbM$,
% it relies implicitly on all the other components.

 \end{enumerate}
  \end{remarks}

The following condition holds for the parts discussed above.

 \begin{condition}
The parameter $\bub$ is an upper bound on the size of every wall and trap.
 \end{condition}

\subsubsection{Cleanness}

Intuitively, a point $x$ is clean in $\cM^{k}$
when none of the mazeries $\cM^{i}$ for $i < k$ has any bad events 
near $x$.
This interpretation will become precise by the rescaling
operation; at this point, we treat cleanness as a primitive, just like
walls. 
Several kinds of
cleanness are needed, depending on the direction in which the absence of
lower-order bad events will be guaranteed.

The set $\cC_{\d}$ is a function of the process $Z_{\d}$, and is used to
formalize (encode) the notions of cleanness given descriptive names below.

 \begin{definition}[One-dimensional cleanness]\label{def:1dim-clean}
For an interval $I = \rint{a}{b}$ or $I=\clint{a}{b}$,
if $(a, b, -1) \in \cC_{\d}$ then we say that point $b$ of
$\bbZ_{+}$ is \df{clean} in $I$ for the sequence $Z_{\d}$.
If $(a,b,1) \in \cC_{\d}$ then we say that point $a$ is clean in $I$.
From now on, whenever we talk about cleanness of an element of $\bbZ_{+}$,
it is always understood with respect to one of the sequences
$Z_{\d}$ for $\d=0,1$ (that is either for the sequence $X$ or for $Y$).

Let us still fix a direction $\d$ and talk about cleanness, and so on, with
respect to the sequence $Z_{\d}$.
A point $x \in \bbZ_{+}$ is called left-clean (right-clean)
if it is clean in all intervals of the form $\rint{a}{x}$, $\clint{a}{x}$ (all intervals
of the form $\rint{x}{b}$, $\clint{x}{b}$).
It is \df{clean} if it is both left- and right-clean.
If both ends of an interval $I$ are clean in $I$ 
then we say $I$ is \df{inner clean}.

To every notion of one-dimensional
cleanness there is a corresponding notion of \df{strong
cleanness}, defined with the help of the process $\cS$ in place of the
process $\cC$.
 \end{definition}

 \begin{figure}
   \centering
\includegraphics[scale=0.6]{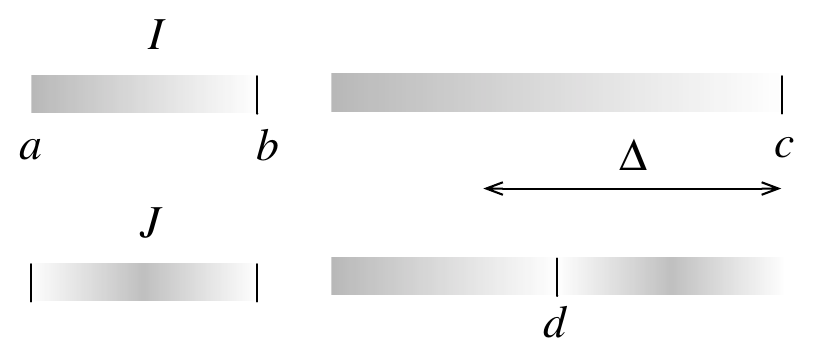}   
   \caption{One-dimensional notions of cleanness.  Point $a$ is not clean in
     interval $I$, but point $b$ is.  Point $c$ is left-clean.  Interval $J$ is
     inner clean.  Point $d$ is clean.}
   \label{fig:one-dim-clean}
 \end{figure}

The relation of strong cleanness to cleanness is dual to the relation of walls
to barriers: every strongly clean point is clean but not vice versa, and every
wall is a barrier, but not vice versa.
This duality is not accidental, since the scale-up operation will define strong
cleanness recursively requiring the absence of nearby barriers, and cleanness
requiring the absence of nearby walls.

 \begin{definition}[Trap-cleanness]\label{def:trap-clean}
For points $u=\pair{u_{0}}{u_{1}}$, $v=\pair{v_{0}}{v_{1}}$, 
$Q = \Rect^{\eps}(u,v)$
where $\eps = \rightarrow$ or $\uparrow$ or nothing,
we say that point $u$ is \df{trap-clean in} $Q$
(with respect to the pair of sequences $\tup{X,Y}$)
if $(u, v, 1, \eps') \in \cS_{2}$, where $\eps' = 0,1,2$ depending on
whether $\eps = \rightarrow$ or $\uparrow$ or nothing.
Similarly, point $v$ is \df{trap-clean in} $Q$
if $(u, v, -1, \eps') \in \cS_{2}$.
 \end{definition}

It is not seen here, but the scale-up operation will define trap-cleanness 
recursively requiring the absence of nearby traps. 

\begin{figure}
  \centering
  \includegraphics[scale=0.6]{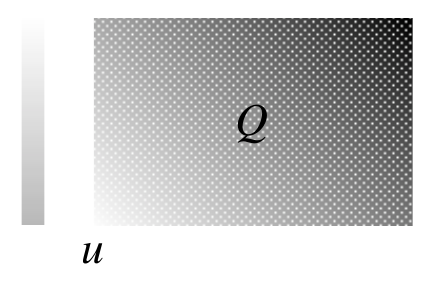}
  \caption{Cleanness in a rectangle.  Point $u$ is trap-clean in rectangle $Q$,
    but is not clean in it, since its projection is not clean in the corresponding
    projection of $Q$. } 
  \label{fig:2d-clean}
\end{figure}

 \begin{definition}[Complex two-dimensional sorts of cleanness]\label{def:H-clean}
We say that point $u$ is \df{clean} in $Q$ when
it is trap-clean in $Q$ and its projections are clean in the
corresponding projections of $Q$.

If $u$ is clean in all such left-open rectangles
then it is called \df{upper right rightward-clean}.
We delete the ``rightward'' qualifier here
if we have closed rectangles in the definition here instead of
left-open ones.
(It is hopeless to illustrate visually the difference made by the ``rightward''
qualifier, but the distinction seems to matter in the proof.)
Cleanness with qualifier ``upward'' is defined similarly.
Cleanness of $v$ in $Q$ and lower left cleanness of $v$ are
defined similarly, using $(u,v, -1, \eps')$,
except that the qualifier is unnecessary: all our rectangles are upper
right closed.

A point is called \df{clean} if it is upper right clean and lower left clean.
If both the lower left and upper right 
points of a rectangle $Q$ are clean in $Q$ then $Q$ is called \df{inner clean}.
If the lower left endpoint is lower left clean and the upper
right endpoint is upper right rightward-clean 
then $Q$ is called \df{outer rightward-clean}.
Similarly for \df{outer upward-clean} and \df{outer-clean}.

We will also use a \df{partial} version of cleanness.
If point $u$ is trap-clean in $Q$ and its projection on the $x$ axis 
is \emph{strongly} clean in the same projection of $Q$
then we will say that $u$ is H-clean in $Q$.
Clearly, if $u$ is H-clean in $Q$ and its projection on the $y$ axis is
clean in (the projection of) $Q$ then it is clean in $Q$.
We will call rectangle $Q$ \df{inner} H-clean if both its lower left and upper
right corners are H-clean in it.

The notion V-clean is defined similarly when we interchange horizontal and
vertical.
 \end{definition}

\subsubsection{Hops}
Hops are intervals and rectangles for which we will be able to
give some guarantees that they can be passed.

\begin{definition}[Hops]\label{def:hop}
A right-closed horizontal interval $I$ is called a \df{hop}
if it is inner clean and contains no vertical wall.
A closed interval $\clint{a}{b}$ is a hop if $\rint{a}{b}$ is a hop.
Vertical hops are defined similarly.

We call a rectangle $I\times J$
a \df{hop} if it is inner clean and contains no trap or wall.
\end{definition}
  
 \begin{remark}
An interval or rectangle that is a hop can be empty: this is the case if
the interval is $\rint{a}{a}$, or the rectangle is, say,
$\Rect^{\rightarrow}(u,u)$.
 \end{remark}

 \begin{definition}[Sequences of walls]\label{def:neighbor-seq}
Two disjoint walls are called \df{neighbors} if the interval
between them is a hop.
A sequence $W_{i} \in \cW$ of walls
$i=1,2,\dotsc$ along with the intervals $I_{1},\dots,I_{n-1}$ between them 
is called a \df{sequence of neighbor walls}
if for all $i > 1$, $W_{i}$ is a right neighbor of $W_{i-1}$.
We say that an interval $I$ is \df{spanned} by the sequence of neighbor
walls $W_{1},W_{2},\dots,W_{n}$ if 
$I=W_{1}\cup I_{1}\cup W_{2}\cup\dots\cup W_{n}$.
We will also say that $I$ is spanned by the sequence $\tup{W_{1},W_{2},\dots}$
if both $I$ and the sequence are infinite and $I=W_{1}\cup I_{1}\cup
W_{2}\cup\dots$.
If there is a hop $I_{0}$ adjacent on the left to $W_{1}$ and a hop
$I_{n}$ adjacent on the right to $W_{n}$ (or the sequence $W_{i}$ is
infinite) then this system is called an
\df{extended sequence of neighbor walls}.
We say that an interval $I$ is \df{spanned} by this extended sequence
if $I=I_{0}\cup W_{1}\cup I_{1}\cup \dots\cup I_{n}$ (and correspondingly
for the infinite case).
 \end{definition}

\subsubsection{Holes}\label{sss.holes}

 \begin{definition}[Reachability]\label{def:reachability}
To each mazery $\cM$ belongs a random graph 
 \begin{equation*}
       \cV = \bbZ_{+}^{2},
\quad  \cG = (\cV, \cE)
 \end{equation*}
where $\cE$ is determined by the  above random processes as in 
Subsection~\ref{subsec:perc}.
We say that point $v$ is \df{reachable} from point $u$ in $\cM$ 
(and write $u \leadsto v$) if it is reachable in $\cG$.
 \end{definition}
  
 \begin{remark}
According to our definitions in Subsection~\ref{subsec:perc}, point $u$ itself
may be closed even if $v$ is reachable from $u$.
 \end{remark}

Intuitively, a hole is a place at which we can pass through a wall.
We will also need some guarantees of being able to reach the hole and being
able to leave it.

 \begin{definition}[Holes]\label{def:holes}
Let $a=\pair{a_{0}}{a_{1}}$, $b=\pair{b_{0}}{b_{1}}$, and
let the interval $I = \rint{a_{1}}{b_{1}}$ be the body of
a horizontal barrier $B$.
For an interval $J = \rint{a_{0}}{b_{0}}$ with 
$|J| \le |I|$  we say that $J$
is a vertical \df{hole passing through} $B$, or \df{fitting} $B$, if
$a\leadsto b$ within the rectangle $J\times\clint{a_{1}}{b_{1}}$.

The above hole is called \df{lower left clean}, upper right clean, and so on,
if the rectangle is.
Consider a point $\pair{u_{0}}{u_{1}}$ with $u_{i}\le a_{i}$, $i=0,1$.
The hole $J$ is called \df{good as seen from} point $u$ if $a$ is H-clean in
$\Rect^{\to}(u,a)$, and $b$ is upper-right rightward
H-clean (recall Definition~\ref{def:H-clean}). 
It is \df{good} if it is good as seen from any such point $u$.
Note that this way the horizontal cleanness is required to be strong, but no
vertical cleanness is required (since the barrier $B$ was not required to
be outer clean).

Horizontal holes are defined similarly.\footnote{
The notion of hole in the present paper is different
from that in~\cite{GacsChat04}.
Holes are not primitives; rather, they are defined with the help of
reachability.}
 \end{definition}

 \begin{remark}
Note that the condition of passing through a wall depends on an interval
slightly larger than the wall itself: it also depends on the left end of the
left-open interval that is the body of the wall. 
 \end{remark}

\subsection{Conditions on the random process}

Most of our conditions on the distribution of process $\cM$
are fairly natural; however, the need for 
some of them will be seen only later.
For example, for Condition~\ref{cond:distr}.\ref{i:distr.hole-lb},
only its special case (in Remark~\ref{rem:distr}.\ref{i:rem.hole-lb}) is well
motivated now: it says
that through every wall there is a hole with sufficiently large probability.
The general case will be used in
the inductive proof showing that the hole lower bound also holds
on compound walls after renormalization (going from $\cM^{k}$ to
$\cM^{k+1}$).

It is a fair question to ask at this point,
why all these conditions are necessary, and whether they are sufficient.
Unfortunately, at this point I can only answer that each condition will
be used in the proof, suggesting their necessity (at least in this proof).
On the other hand, the proof that the scale-up operation conserves the
conditions, shows their sufficiency.

The combinatorial conditions derive ultimately from the necessity of proving
Lemma~\ref{lem:approx} (Approximation).
The dependency conditions and probability estimates derive ultimately from the
necessity of proving inequality~\eqref{eq:main}.
But some conditions are introduced just in order to help the proof of the
conservation in the scale-up.

 \begin{definition}
The function 
 \begin{equation}\label{eq:p(,)}
  \p(\r,l)
 \end{equation}
is defined as the supremum of probabilities (over all points $t$)
that any barrier with rank $\r$ and size $l$ starts at $t$
conditional over all possible conditions of the form $Z_{\d}(t)=k$
for $k\in\{1,\dots,\m\}$.
 \end{definition}

The constant $\hxp$ has been introduced in~\eqref{eq:hxp}.
Its choice, as well as the choice of some of the other expressions we are about
to introduce, will be motivated in Section~\ref{sec:params}.
We will use three additional constants, $\aux_{1},\aux_{2},\aux_{3}$.
Constant $\aux_{1}$ will be chosen at the end of the proof of
Lemma~\ref{lem:compound-contrib}, 
$\aux_{2}$ in the proof or Lemma~\ref{lem:pub},
while $\aux_{3}$ will be chosen
at the end of the proof of Lemma~\ref{lem:all-compound-hole-lb}. 

\begin{definition}\label{def:lg-def}
We will make use of constant
 \begin{equation}\label{eq:lg-def}
   \lg = 2^{1/2}.
 \end{equation}
\end{definition}

Let us define a function that will serve as an
an upper bound on $\sum_{l} \p(\r, l)$.

 \begin{definition}[Barrier probability upper bound]
Let
 \begin{equation}\label{eq:wall-prob}
  \p(\r) = \aux_{2} \r^{-\aux_{1}} \lg^{-\r}.
 \end{equation}
 \end{definition}

The term $\aux_{2}\r^{-\aux_{1}}$
serves for absorbing some lower-order factors that arise in estimates
like~\eqref{eq:compound-hole-lb}.

 \begin{definition}[Hole probability lower bound]
Let
 \begin{align}\label{eq:h-def}
         \h(\r)  &= \aux_{3}\lg^{-\hxp\r}.
 \end{align}
 \end{definition}

Here come the conditions.

 \begin{condition}\label{cond:distr}\
 \begin{enumerate}[\upshape 1.]

  \item\label{i:distr.indep}(Dependencies)
  \begin{enumerate}[\upshape a.]

    \item\label{i:distr.indep.trap}
For any rectangle $I \times J$, the event that it is a
trap is a function of the pair $X(I), Y(J)$.

    \item\label{i:distr.indep.barrier}
For a vertical wall value $E$ the event $\setof{E\in\cB}$ (that is the event
that it is a vertical barrier) is a function of $X(\Body(E))$.

Similarly for horizontal barriers.

    \item\label{i:distr.indep.left-clean} 
For integers $a < b$, and $\d=0,1$, the events defining strong cleanness, that
is $\setof{(a,b,-1) \in \cS_{\d}}$ and $\setof{(a,b,1) \in \cS_{\d}}$,
are functions of $Z_{\d}(\rint{a}{b})$.

When $Z$ is fixed, then for a fixed $a$, the
(strong and not strong) cleanness of $a$ in $\rint{a}{b}$ 
is decreasing as a function of $b-a$, 
and for a fixed $b$, the (strong and not strong) cleanness 
of $b$ in $\rint{a}{b}$ is decreasing as a function of $b-a$. 
These functions reach their minimum at $b-a=\bub$: thus, if
$x$ is (strongly or not strongly) 
left clean in $\rint{x-\bub}{x}$ then it is (strongly or not strongly) 
left clean.

    \item\label{i:distr.indep.ur-clean}
For any rectangle $Q=I \times J$,
the event that its lower left corner is trap-clean in $Q$,
is a function of the pair $X(I), Y(J)$.

Among rectangles with a fixed lower left corner, the
event that this corner is
trap-clean in $Q$ is a decreasing function of $Q$
(in the set of rectangles partially ordered by containment).
In particular, the trap-cleanness of $u$ in $\Rect(u,v)$ implies its 
trap-cleanness
in $\Rect^{\rightarrow}(u,v)$ and in $\Rect^{\uparrow}(u,v)$.
If $u$ is upper right trap-clean in the left-open or bottom-open or
closed square of size $\bub$,
then it is upper right trap-clean in all rectangles $Q$ of the same
type. 
Similar statements hold if we replace upper right with lower left.

  \end{enumerate}

Whether a certain wall value $E=(B,\r)$ is a vertical barrier depends only on
$X(B)$.
Whether it is a vertical 
wall depends also on only on $X$, but may depend on the values of $X$
outside $B$.
Similarly, whether a certain horizontal interval is inner clean depend only the
sequence $X$ but may depend on its elements outside it,
but whether it is strongly inner clean depends only on $X$ inside the interval.

Similar remarks apply to horizontal wall values and vertical cleanness.

  \item\label{i:distr.combinat}(Combinatorial requirements)
   \begin{enumerate}[\upshape a.]
  
    \item\label{i:distr.inner-clean}
A maximal external interval  (see
Definition~\ref{def:barrier}) of size $\ge\bub$ or one starting at $-1$
is inner clean.

    \item\label{i:distr.cover}
An interval $I$ that is surrounded by maximal external intervals of size
$\ge\bub$ is spanned by a
sequence of neighbor walls (see Definition~\ref{def:neighbor-seq}).
This is true even in the case when $I$ starts at 0 and even if it is infinite.
To accomodate these cases, we require the following, which is somewhat harder to
parse:
Suppose that interval $I$ is adjacent on the left to a maximal external interval
that either starts at $-1$ or has size $\ge \bub$.
Suppose also that it is either adjacent on the right to a similar
interval or is infinite.
Then it is spanned by a (finite or infinite) sequence of neighbor walls.
In particular, the whole line is spanned by an extended sequence of
neighbor walls.

    \item\label{i:distr.clean.1}
If a (not necessarily integer aligned)
right-closed interval of size $\ge 3 \bub$ contains no wall,
then its middle third contains a clean point.

    \item\label{i:distr.clean.2}
Suppose that a rectangle $I\times J$ with (not necessarily integer aligned)
right-closed $I,J$ with $|I|, |J| \ge 3 \bub$ contains
no horizontal wall and no trap, and $a$ is a right 
clean point in the middle third of $I$.
There is an integer $b$ in the middle third of $J$ such that the point
$\pair{a}{b}$ is upper right clean.
A similar statement holds if we replace upper right with
lower left (and right with left).
Also, if $a$ is clean then we can find a point $b$ 
in the middle third of $J$ such that $\pair{a}{b}$ is clean.

There is also a similar set of statements if we vary $a$ instead of $b$.

   \end{enumerate}

   \item\label{i:distr.bounds}(Probability bounds)
  \begin{enumerate}[\upshape a.]

   \item\label{i:distr.trap-ub}
Given a string $x = \tup{x(0), x(1), \dotsc}$, a point $\pair{a}{b}$,
let $\cF$ be the event that a trap starts at $\pair{a}{b}$.
Let $\s\in\{1,\dots,\m\}$, then
 \begin{equation*}
     \Pbof{\cF \mid X = x,\; Y(b-1)=\s} \le \tub.
 \end{equation*}
The same is required if we exchange $X$ and $Y$.

   \item\label{i:distr.wall-ub}
We have $\p(\r) \ge \sum_{l} \p(\r,l)$.

  \item\label{i:distr.ncln}
We require $\ncln < 0.1$, and that for all
$k\in\{1,\dots,m\}$, for all $a < b$ and all $u=\pair{u_{0}}{u_{1}}$,
$v=\pair{v_{0}}{v_{1}}$, for all sequences $y$, the following quantities are all
$\le\ncln/2$:
\begin{align}\label{eq:ncln.1dim}
   &\Pbof{ a \txt{ (resp. $b$) } 
  \txt{ is not strongly clean in } \rint{a}{b} \mid X(a)=k},
\\\label{eq:ncln.lower}
 &\Pbof{ u \txt{ (resp. $v$) }
  \txt{ is not trap-clean in } 
     \Rect^{\rightarrow}(u,v) \mid X(u_{0})=k, Y=y},
\\\label{eq:ncln.upper}
   &\Pbof{u \txt{ (resp. $v$) }
  \txt{ is not trap-clean in }
     \Rect(u,v) \mid X(u_{0}-1)=k, Y=y},
 \end{align}
and similarly with $X$ and $Y$ reversed.
 \begin{figure}
\begin{center}
   \includegraphics{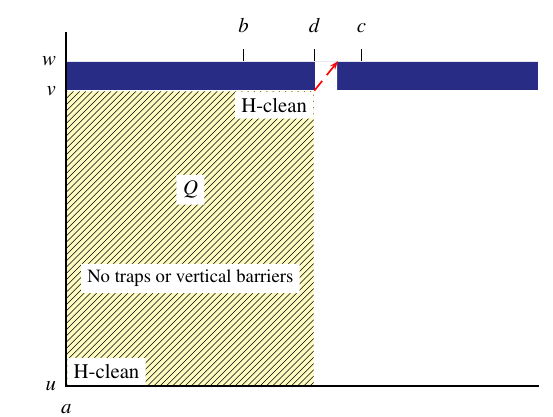}
 \caption{Condition~\protect\ref{cond:distr}.\ref{i:distr.hole-lb}.
The hatched rectangle is inner H-clean and contains no traps or vertical barriers.}
  \end{center}
 \end{figure}

   \item\label{i:distr.hole-lb}
Let $u \le v < w$, and $a$ be given with $v-u \le 12\bub$, and define
 \begin{equation*}
                 b = a + \cei{(v-u)/2},
  \qquad         c = a + (v-u) + 1.
 \end{equation*}
Assume that $Y=y$ is fixed in such a way that
$B$ is a horizontal wall of rank $\r$ with body $\rint{v}{w}$.
For a $d \in \clint{b}{c-1}$ let $Q(d)=\Rect^{\rightarrow}(\pair{a}{u}, \pair{d}{v})$.
Let $E = E(u, v, w;\, a)$ be the event (a function of $X$)
that there is a $d$ such that
 \begin{enumerate}[\upshape (i)]
  \item A vertical hole fitting $B$ starts at $d$.
  \item Rectangle $Q(d)$ contains no traps or vertical barriers, and is inner H-clean.
 \end{enumerate}
Let $k\in\{1,\dots,\m\}$.
Then
 \begin{equation*}
   \Pbof{E \mid X(a)=k,\; Y = y} \ge (c-b)^{\hxp} \h(\r).
 \end{equation*}
The same is required if we exchange horizontal and vertical.

  \end{enumerate}
 \end{enumerate}
 \end{condition}

\begin{remarks}\label{rem:distr}\  % hyperref complains without this
 \begin{enumerate}[\upshape 1.]

  \item\label{i:middle-third}
 Conditions~\ref{cond:distr}.\ref{i:distr.clean.1} 
and~\ref{cond:distr}.\ref{i:distr.clean.2} imply the following.
Suppose that a right-upper closed square $Q$ of size $3 \bub$ contains
no wall or trap.
Then its middle third contains a clean point.
In particular, this implies Condition~\ref{cond:dense}.

  \item\label{i:rem.hole-lb}
The most important special case of
Condition~\ref{cond:distr}.\ref{i:distr.hole-lb} is
$v = u$, implying $b = a$, $c = b + 1$: then it says 
that for any horizontal wall $B$ of rank $\r$, at any point $a$,
the probability that there is a vertical hole passing through $B$ at point
$a$ is at least $\h(\r)$.

 \end{enumerate}
\end{remarks}
  
The graph $\cG$ defined in Definition~\ref{def:reachability}
is required to satisfy the following conditions.

 \begin{condition}[Reachability]\label{cond:reachable}
We require $0 \le \slb < 0.5$.
Let $u, v$ be points with $\minslope(u, v) \ge \slb$.
If they are the starting and endpoint of a rectangle 
that is a hop, then $u \leadsto v$.
The rectangle in question is allowed to be
bottom-open or left-open, but not both.
 \end{condition}

 \begin{example}[Base case]\label{xmp:base}
The clairvoyant demon problem can be seen as a special case of a mazery.
Let us choose the scale parameter $\bub$ for any valuel $\ge 1$, and $\slb = 0$,
that is there is no 
lower bound on the minimum slope for the reachability condition.
The parameters $\R>0$ and $0<\ncln<0.1$ are chosen arbitrarily,
and we choose 
 \begin{align}\label{eq:base-tub}
 1>\tub\ge\frac{1}{\m-1},   
 \end{align}
where $\m$  is the size of the complete graph on which
 the random walks are performed.

Let $\cT=\setof{\pair{i}{j}: X(i)=Y(j)}$, that is
traps are points $\pair{i}{j}$ with $X(i) = Y(j)$.
We set $\cB_{d}=\cW_{d}=\emptyset$, that is 
there are no barriers (and then, of course, no walls).

Let $\cS_{d}=\cC_{d}=\bbZ_{+}^{{(2)}}\times\{-1,1\}$  for $d=0,1$.
In other words, every point is strongly clean in all one-dimensional senses.
Also $\cS_{2}=\bbZ_{+}^{(2)}\times\bbZ_{+}^{(2)}\times\{-1,1\}\times\{0,1,2\}$,
that is every point is trap-clean in all senses.

All combinatorial and dependency conditions are satisfied trivially.
Of the probability bounds, 
Condition~\ref{cond:distr}.\ref{i:distr.trap-ub}
is satisfied by our requirement~\eqref{eq:base-tub}.
Since there are no walls and every point is clean in all possible ways, the
other probability bounds are satisfied trivially.

Since now in a trap-free rectangle nothing blocks reachability, the reachability
condition also holds trivially.
Note that it is violated in the bottom-left open 
rectangle $\rint{0}{1} \times \rint{0}{1}$ if $X(0)=1$, $X(1)=2$, 
$Y(0)=2$, $Y(1)=1$.
Indeed, the traps $\pair{0}{1}$, $\pair{1}{0}$ are not
part of the rectangle where they are prohibited, but
they block point $\pair{0}{0}$ from point $\pair{1}{1}$.
 \end{example}

 \section{The scaled-up structure}\label{sec:plan}

We will use Example~\ref{xmp:base} as the mazery $\cM^{1}$ in our sequence of
mazeries $\cM^{1},\cM^{2},\dots$ where $\cM^{k+1}=(\cM^{k})^{*}$.
In this section, we will define the scaling-up operation 
$\cM \mapsto \cM^{*}$: we still postpone to Section~\ref{sec:params}
the definition of several parameters and probability bounds for $\cM^{*}$.

Let us recall the meaning of the scale-up operation from
Section~\ref{sec:outline}.
Our final goal is to prove reachability of points far away from the origin, with
large probability.
In our model $\cM$, reachability is guaranteed in a rectangle $Q=\Rect(u,v)$ from
$u$ to $v$ if $u,v$ are inner clean in $Q$, and there are no traps or walls in
$Q$.
The absence of traps cannot be guaranteed in our base model when the rectangle
is large.
It should be sufficient for traps to be far from each other, but
even this condition will fail occasionally.

The idea of the scale-up strategy is to define new kinds of ``obstacles''
on which we can blame such failures. 
If these obstacles have sufficiently small probability, then they can be
regarded as the traps and walls of a model $\cM^{*}$.
The crucial combinatorial test of this procedure is the proof of the
reachability condition in model $\cM^{*}$, which is Lemma~\ref{lem:approx}.
This may not be so complicated in going up just one step from the base model
$\cM^{1}$ to $(\cM^{1})^{*}=\cM^{2}$, since the base model has no walls.
But other models $\cM^{k}$ do have them: in this case, reachability
in $\cM^{k+1}$ in a rectangle $Q$ without walls means that $Q$ may have
walls and traps of $\cM^{k}$, just no walls or traps of $\cM^{k+1}$.
So, the walls and traps of $\cM^{*}$ impersonate the difficulties of passing
through the walls and around the traps of $\cM$.
Walls of $\cM^{*}$ will also be used as scapegoats for some excessive
correlations among traps of $\cM$.

 \subsection{The scale-up construction}\label{subsec:plan.constr}

Some of the following parameters will be given values only later, but they
are introduced by name here.

 \begin{definition}\label{def:f-g}
Let $\slopeincr$ be a constant and let parameters $\f,\g$ satisfy
 \begin{equation}\label{eq:bub-g-f}
 \begin{split}
   \slopeincr &= 500,
\\  \bub / \g &\le \g / \f < (0.5 - \slb) / (2\slopeincr).
 \end{split}
 \end{equation}
 \end{definition}

The parameters $\bub\ll\g\ll\f$ will be different for each level of the
construction.
The scale parameter $\bub$ is part of the definition of a mazery.
Here is the approximate meaning of $\f$ and $\g$:
We try not to permit walls closer than $\f$ to each other,
and we try not to permit intervals larger than $\g$ without holes.

 \begin{definition}\label{def:new-slb}
Let $\slb^{*} = \slb + \slopeincr \g / \f$.
 \end{definition}

The value $\bub^{*}$ will be defined later, but we will guarantee the inequality
 \begin{equation}\label{eq:bub-star}
  3\f \le \bub^{*}.
 \end{equation}
After defining the mazery $\cM^{*}$, eventually we will have to prove the 
required properties.
To be able to prove Condition~\ref{cond:reachable}
for $\cM^{*}$, we will introduce some new walls
and traps in $\cM^{*}$ whenever some larger-scale obstacles prevent
reachability.
There will be two kinds of new walls, so-called \df{emerging} walls, and
\df{compound} walls. 
A pair of traps too close to each other will define,
under certain conditions, a compound trap, which becomes part of
$\cM^{*}$.
A new kind of trap, called a trap of the missing-hole kind will arise when
some long stretch of a low-rank wall is without a hole.

For the new value of $\R$ we require
 \begin{equation}\label{eq:R-cond}
   \R^{*} \le 2 \R - \log_{\lg} \f.
 \end{equation}

 \begin{definition}[Light and heavy]
Barriers and walls of rank lower than $\R^{*}$ are called \df{light}, the
other ones are called \df{heavy}.
 \end{definition}

Heavy walls of $\cM$ will also be walls of $\cM^{*}$ (with some exceptions given
below).
We will define walls only for either $X$ or $Y$, but it is understood
that they are also defined when the roles of $X$ and $Y$ are reversed.

The rest of the scale-up construction will be given in the following steps.

 \begin{defstep}[Cleanness]\label{defstep:clean}
For an interval $I$, its right endpoint $x$ will be called
clean in $I$ for $\cM^{*}$ if 
 \begin{enumerate}[--]
  \item It is clean in $I$ for $\cM$.
  \item The interval $I$ contains no wall
of $\cM$ whose right end is closer to $x$ than $\f/3$.
 \end{enumerate}
We will say that a point is strongly clean in $I$ for $\cM^{*}$ if
it is strongly clean in $I$ for $\cM$ and $I$ contains no barrier
of $\cM$ whose right end is closer to it than $\f/3$.
Cleanness and strong cleanness of the left endpoint is defined similarly.

Let a point $u$ be a starting point or endpoint of a rectangle $Q$.
It will be called trap-clean in $Q$ for $\cM^{*}$ if
 \begin{enumerate}[--]
  \item It is trap-clean in $Q$ for $\cM$.
  \item Any trap contained in $Q$ is at a distance $\ge \g$ from $u$.
 \end{enumerate}
 \end{defstep}

For the next definitions, Figure~\ref{fig:compound-traps} may help.

 \begin{defstep}[Uncorrelated traps]\label{defstep:uncorrel}
A rectangle $Q$ is called an \df{uncorrelated compound trap} if
it contains two traps with disjoint projections, with a distance of their
starting points at most $\f$, and if it is minimal among the rectangles
containing these traps.
 \end{defstep}
Clearly, the size of an uncorrelated trap is bounded by $\bub+\f$.

 \begin{defstep}[Correlated trap]\label{defstep:correl}
  Let 
 \begin{equation}\label{eq:g-prime}
  \g' = 2.2 \g,\; \l_{1}=7\bub,\; \l_{2} = \g'.
 \end{equation}
(Choice motivated by the proof of Lemmas~\ref{lem:clean} and~\ref{lem:approx}.)
    \begin{figure}
      \centering
         \includegraphics[scale=0.5]{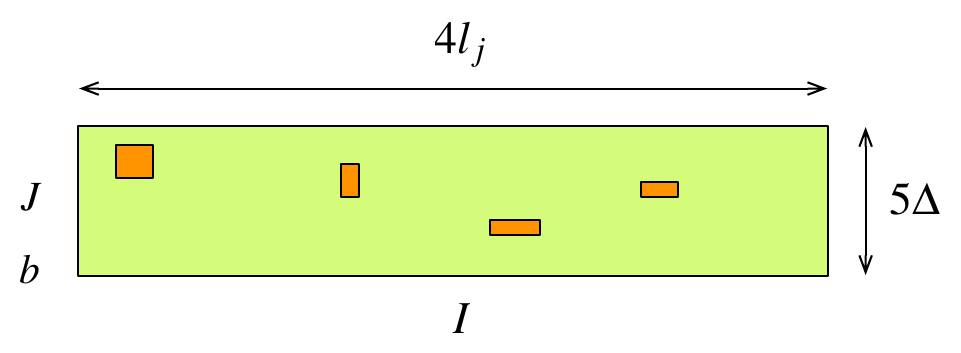}
 \caption{A correlated trap}
    \end{figure}
For a $j\in\{1,2\}$ let $I$ be a closed interval with length $|I| = 4 \l_{j}$,
and $b\in\bbZ_{+}$, with $J=\clint{b}{b+5 \bub}$.
Let $x(I),y(J)$ be given.
We say that event 
  \[
   \cL_{j}(x,y,I,b)
  \]
  holds if for all right-closed intervals $\hat I \sbsq I$ of size $\l_{j}$,
the rectangle $\hat I \times J$ contains a trap.
We will say that $I \times J$ is a \df{horizontal correlated trap} 
of type $j$ if $\cL_{j}(x,y,I,b)$ holds and for all $\s$ in $\{1,\dots,\m\}$,
\begin{equation}\label{eq:correl-trap-ineq}
  \Pbof{\cL_{j}(x,Y,I,b) \mid X(I) = x(I),\, Y(b-1) = \s} \le \tub^{2}.
\end{equation}
Vertical correlated traps are defined analogously\footnote{
  The smallness of the conditional probability in the other direction
will be proved in Lemma~\ref{lem:correl}, without
having to require it.}.
 \end{defstep}

 \begin{defstep}[Traps of the missing-hole kind]\label{defstep:missing-hole}
Let $I$ be a closed interval of size $\g$, let $b$ be
a site with $J = \clint{b}{b + 3\bub}$.
   \begin{figure}
     \centering
       \includegraphics[scale=0.5]{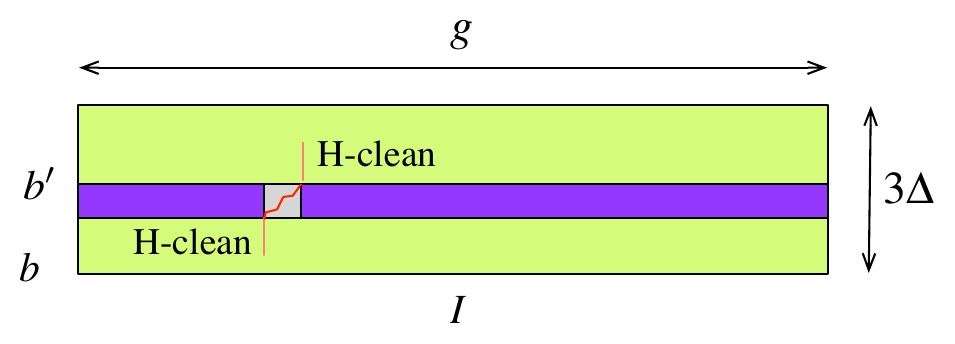}
     \caption{A trap of the missing-hole kind: good holes (see
       Definition~\ref{def:holes}, the type of hole shown in the picture with
       H-cleanness on the lower-left and upper-right) are missing}
   \end{figure}
Let $x(I),y(J)$ be fixed.
We say that event 
 \[
  \cL_{3}(x,y,I,b)
 \]
holds if, with there is a $b'>b+\bub$ such that
$\rint{b+\bub}{b'}$ is the body of
a light horizontal potential wall $W$, and no good vertical hole
(in the sense of Definition~\ref{def:holes}) 
$\rint{a_{1}}{a_{2}}$ with $\rint{a_{1}-\bub}{a_{2}+\bub} \sbsq I$ 
passes through $W$.

We say that $I \times J$ is a \df{horizontal trap of the missing-hole kind}  
if event $\cL_{3}(x,y,I,b)$ holds 
and for all $\s\in\{1,\dots,\m\}$ we have
 \begin{equation}\label{eq:missing-hole-trap-ineq}
 \Pbof{\cL_{3}(x,Y,I,b) \mid X(I) = x(I),\, Y(b-1) = \s} \le \tub^{2}.
 \end{equation}
 \end{defstep}

Inequalities~\eqref{eq:bub-g-f} and~\eqref{eq:bub-star} bound
the size of all new traps by $\bub^{*}$.

 \begin{defstep}[Emerging walls]\label{defstep:emerg}
It is in this definition where the difference between barriers and walls first
appears in the paper constructively.
   \begin{figure}
     \centering
       \includegraphics[scale=0.5]{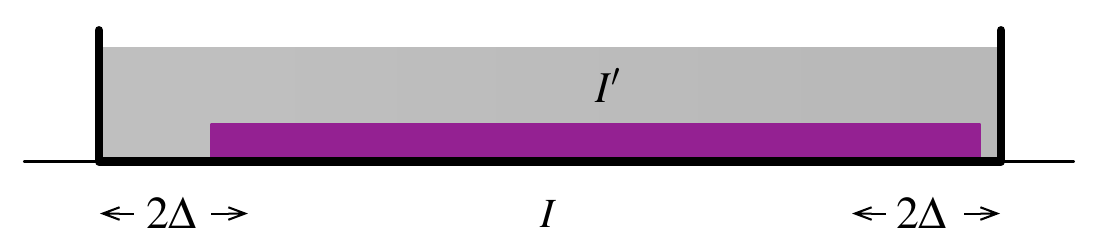}
     \caption{An emerging barrier}
   \end{figure}
We define some objects as barriers, and then designate some of the barriers (but
not all) as walls.

A vertical emerging barrier is, essentially, a horizontal interval over which
the conditional probability of a bad event $\cL_{j}$ is not
small (thus preventing a new trap).
But in order to find enough barriers, the ends are allowed to be sligthly extended.
Let $x$ be a particular value of the sequence $X$ over an 
interval $I=\rint{u}{v}$.
For any $u' \in \rint{u}{u+2\bub}$, $v' \in \rint{v-2\bub}{v}$,
let us define the interval $I'=\clint{u'}{v'}$.
We say that interval $I$ 
is the body of a vertical \df{barrier} of the \df{emerging kind}, of type 
$j\in\{1,2,3\}$ if the following inequality holds:
 \begin{equation}\label{eq:emerg-barrier}
 \sup_{I',k} \Pbof{\cL_{j}(x,Y,I',1) \mid X(I') = x(I'),\, Y(0) = k} 
  > \tub^{2}.  
 \end{equation}
To make it more explicit, for example interval $I$ is an emerging barrier of
type 1 for the process $X$
if it has a closed subinterval $I'$ of size $4\l_{1}$ within $2\bub$ of its
two ends, such that conditionally over the value of $X(I')$ and $Y(0)$, with
probability $>\tub^{2}$, for 
all right-closed subintervals $\hat I$ of $I'$, the rectangle
$\hat I\times \clint{b}{b+5 \bub}$ contains a trap.
More simply, the value $X(I')$ makes not too improbable (in terms of a randomly
chosen $Y$) for a sequence of closely placed traps to exist reaching horizontally across
$I'\times \clint{b}{b+5 \bub}$.

Note that emerging barriers of type 1 are smallest, and those of type 2
are largest.
More precisely, let 
 \[
  L_{1}=4\l_{1},\quad L_{2} = 4\l_{2}, \quad L_{3}=\g.
 \]
Then emerging barriers of type $j$ have length in 
$L_{j}+\clint{0}{4\bub}$.

We will designate some of the emerging barriers as walls.
We will say that $I$ is a \df{pre-wall} of the emerging kind if also the
following properties hold:
 \begin{alphenum} 

  \item\label{i:emerg.hop}
Either $I$ is an external hop of $\cM$ or it is the union of a dominant
light wall and one
or two external hops of $\cM$, of size $\ge\bub$, surrounding it.

  \item\label{i:emerg.outside}
Each end of $I$ is adjacent to either an external hop of size $\ge\bub$
or a wall of $\cM$.

 \end{alphenum}
 \begin{figure}
   \centering
     \includegraphics[scale=0.6]{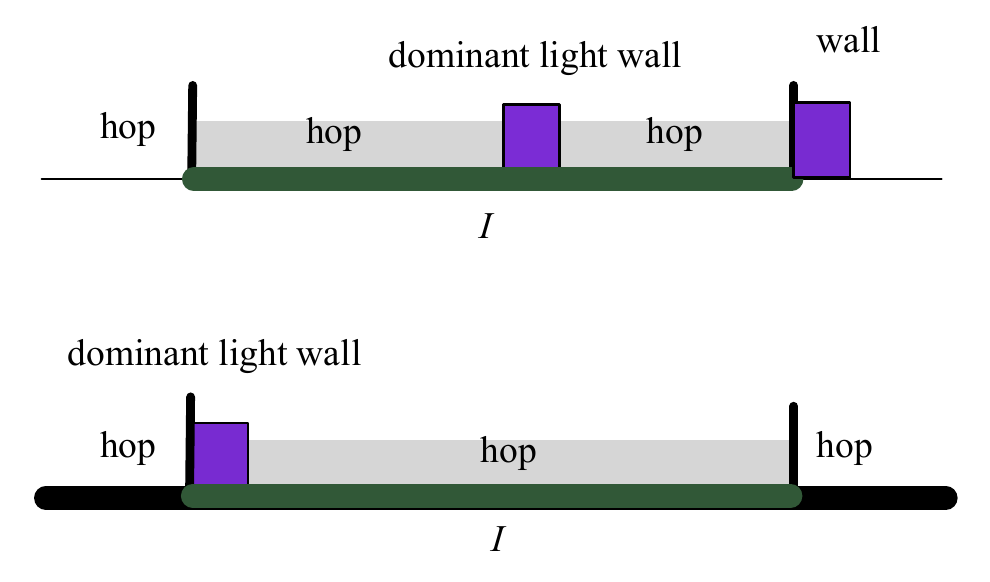}
   \caption{Some possibilities for a pre-wall}
 \end{figure}

Now, for $j=1,2,3$, list all emerging pre-walls of type $j$ in a sequence
$\tup{B_{j1},B_{j2},\dots}$. 
First process pre-walls $B_{11},B_{12},\dots$ one-by-one.
Designate $B_{1n}$ a wall if and only if it
is disjoint of all emerging pre-walls designated as walls earlier.
Next process the sequence $\tup{B_{31},B_{32},\dots}$.
Designate $B_{3n}$ a wall if and only if it
is disjoint of all emerging pre-walls designated as walls earlier.
Finally process the sequence $\tup{B_{21},B_{22},\dots}$.
Designate $B_{2n}$ a wall if and only if it
is disjoint of all emerging pre-walls designated as walls earlier.

To emerging barriers and walls, we assign rank 
 \begin{equation}\label{eq:emerg-rank}
  \hat\R > \R^{*}
 \end{equation}
to be determined later.

 \end{defstep}

 \begin{defstep}[Compound walls]\label{defstep:compound}
We make use of a certain sequence of integers: 
 \begin{equation}\label{eq:jump-types}
   d_{i} = \begin{cases}
               i &\txt{ if $i=0,1$},
\\ \cei{\lg^{i}} &\txt{ if $i \ge 2$}.
           \end{cases}
 \end{equation}
A \df{compound barrier} occurs in $\cM^{*}$ for $X$
wherever barriers $W_{1}, W_{2}$ occur (in this order) for $X$
at a distance $d \in \lint{d_{i}}{d_{i+1}}$, $d \le \f$,
and $W_{1}$ is light.
We will call this barrier a wall if $W_{1},W_{2}$ are neighbor walls
(that is, they are walls separated by a hop).
We denote the new compound wall or barrier by
 \[
  W_{1} + W_{2}.
 \]
Its body is the smallest right-closed interval containing the bodies of
$W_{j}$. 
For $\r_{j}$ the rank of $W_{j}$, we will say that the compound wall or
barrier in question has \df{type} 
 \[
   \ang{\r_{1},\r_{2},i}.
 \]
Its rank is defined as
 \begin{equation}\label{eq:compound-rank}
  \r  = \r_{1} + \r_{2} - i.
 \end{equation}
Thus, a shorter distance gives higher rank.
This definition gives 
  \[
    \r_{1}+\r_{2}-\log_{\lg}\f \le \r \le \r_{1}+\r_{2}.
 \]
 Inequality~\eqref{eq:R-cond} will make sure that the rank of 
the compund walls is lower-bounded by $\R^{*}$.

Now we repeat the whole compounding step, introducing compound walls and
barriers in which now $W_{2}$ is required to be light.
The barrier $W_{1}$ can be any barrier introduced until now, also a compound
barrier introduced in the first compounding step.
\begin{figure}
  \centering
\includegraphics[scale=0.7]{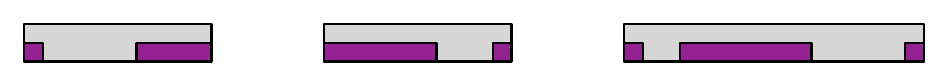}
    \caption{Three (overlapping) 
types of compound barrier obtained: light-any, any-light,
light-any-light.
Here, ``any'' can also be a recently defined emerging barrier.}
\end{figure}
 \end{defstep}

The walls that will occur as a result of the compounding operation are 
of the type $L$-$*$, $*$-$L$, or $L$-$*$-$L$, where $L$ is a light wall
of $\cM$ and $*$ is any wall of $\cM$ or an emerging wall of
$\cM^{*}$.
Thus, the maximum size of a compound wall is
 \[
   \bub + \f + (4\g'+4\bub) + \f + \bub < \bub^{*},
 \]
where we used~\eqref{eq:bub-g-f} and~\eqref{eq:bub-star}.

 \begin{defstep}[Finish]\label{defstep:finish}
The graph $\cG$ does not change in the scale-up: $\cG^{*}=\cG$.
Remove all traps of $\cM$.

Remove all light walls and barriers.
If the removed light wall was dominant, remove also all other
walls of $\cM$ (even if not light) contained in it.
 \end{defstep}

The reader may miss an explanation for why
we introduced exactly these higher-order traps and walls, and no others.
The only explanation available is, however, that these are the only objects
whose absence is used in the proof of the reachability property in $\cM^{*}$
(Lemma~\ref{lem:approx}).

 \subsection{Combinatorial properties}\label{subsec:plan.prop}

Let us prove some properties of $\cM^{*}$ that can already be established.
Note first that Condition~\ref{cond:reachable-hop} follows from
Condition~\ref{cond:reachable}, Condition~\ref{cond:dense} follows from
Conditions~\ref{cond:distr}.\ref{i:distr.clean.1}-\ref{i:distr.clean.2},
and~\ref{cond:Qbd} for $\cM^{*}$
follows from the definition of cleanness in $\cM^{*}$
given in the present section.

 \begin{lemma}\label{lem:distr.indep}
The new mazery $\cM^{*}$ satisfies
Condition~\ref{cond:distr}.\ref{i:distr.indep}.
 \end{lemma}
 \begin{proof}
We will see that all the properties in the condition follow essentially
from the form of our definitions.

Condition~\ref{cond:distr}.\ref{i:distr.indep.trap} says that
for any rectangle $I \times J$, the event that it is a trap
is a function of the pair $X(I), Y(J)$.
To check this, consider all possible traps of $\cM^{*}$.
We have the following kinds:
 \begin{enumerate}[--]

  \item Uncorrelated and correlated compound trap.
The form of the definition shows that this event only depends on
$X(I), Y(J)$.

  \item Trap of the missing-hole kind.
Since the definition of good holes uses H-cleanness, this
depends in $\cM$ only on a $\bub$-neighborhood of a point.
Therefore event $\cL_{3}(x,y,I,b)$ depends only on the $x(I)$ part of $x$.
Since $W$ is required to be a potential wall, the event only depends
on the $y(J)$ part of $y$.
The conditional probability inequality also depends only on $x(I)$.

 \end{enumerate}

Condition~\ref{cond:distr}.\ref{i:distr.indep.barrier} says that, say,
for a vertical wall value $E$ the event $\setof{E\in\cB_{0}}$ (that is the event
that it is a vertical barrier) is a function of $X(\Body(E))$.
There are two kinds of vertical barriers in $\cM^{*}$: emerging and
compound barriers.
The definition of both of these refers only to $X(\Body(E))$.

Condition~\ref{cond:distr}.\ref{i:distr.indep.left-clean} says first that
for every interval $I=\rint{a}{b}$, the strong cleanness of $a$ or $b$ in $I$
are functions of $Z_{\d}(I)$.
The property that $a$ or $b$ is strongly clean in interval $I$ in $\cM^{*}$ 
is defined in terms of strong cleanness in $\cM$ and the
absence of barriers contained in $I$.
Therefore strong cleanness of $a$ or $b$ in $I$ for 
$\cM^{*}$ is a function of $Z_{\d}(I)$.

Since (strong) cleanness in $I$ for $\cM$ is a decreasing
function of $I$, and the property stating the absence of walls
(barriers) is a decreasing function of $I$, (strong) cleanness for
$\cM^{*}$ is also a decreasing function of $I$.
The inequality $\f/3 + \bub < \bub^{*}$, implies that
these functions reach their minimum for $|I|=\bub^{*}$.

Condition~\ref{cond:distr}.\ref{i:distr.indep.ur-clean} says first that
for any rectangle $Q=I \times J$,
the event that its lower left corner is trap-clean in $Q$,
is a function of the pair $X(I), Y(J)$.
If $u$ is this point then, our definition of its trap-cleanness for
$\cM^{*}$ in rectangle $Q$ required the following:
 \begin{enumerate}[--]
  \item It is trap-clean in $Q$ for $\cM$;
  \item The starting point of any trap in $Q$ is at a 
distance $\ge \g$ from $u$.
 \end{enumerate}
All these requirements refer only to the projections of $Q$ and depend
therefore only on the pair $X(I), Y(J)$.

It can also be seen that, among rectangles with a fixed lower left corner,
the event that this corner is trap-clean
for $\cM^{*}$ in $Q$ is a decreasing function of $Q$ 
(in the set of rectangles partially ordered by containment).
And, since $\g+\bub < \bub^{*}$, if point
$\pair{x}{y}$ is upper right trap-clean in a square of size $\bub^{*}$,
then it is upper right trap-clean.

 \end{proof}

 \begin{lemma}\label{lem:cover}
The mazery $\cM^{*}$ satisfies 
conditions~\ref{cond:distr}.\ref{i:distr.inner-clean}
and~\ref{cond:distr}.\ref{i:distr.cover}.
 \end{lemma}
 \begin{Proof}
We will prove the statement only for vertical walls; it is proved for
horizontal walls the same way.
In what follows, ``wall'', ``hop'', and so on, mean vertical wall, horizontal
hop, and so on.

Let $\tup{U_{1},U_{2},\dots}$ be a (finite or infinite)
sequence of disjoint walls of $\cM$ and
$\cM^{*}$, and let $I_{0},I_{1},\dots$ be the (possibly empty) 
intervals separating them
(interval $I_{0}$ is the interval preceding $U_{1}$).
This sequence will be called \df{pure} if 
 \begin{alphenum}

  \item The intervals $I_{j}$ are hops of $\cM$.

  \item\label{i:pure.external}
$I_{0}$ is an external interval of $\cM$ starting at $-1$, while
$I_{j}$ for $j>0$ is external if its size is $\ge 3\bub$.

 \end{alphenum}

 \begin{step+}{step:cover.initial-pure}
Let us build an initial pure sequence of $\cM$ which has also an additional property:
\emph{every dominant light wall of $\cM$ belongs to it}.
 \end{step+}
 \begin{prooof}
   \begin{figure}
     \centering
         \includegraphics[scale=0.5]{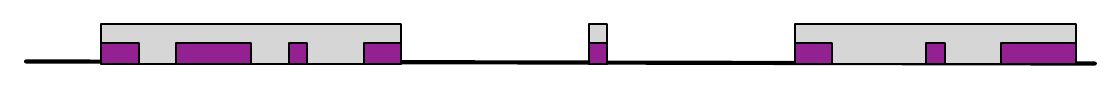}
     \caption{An initial pure sequence.  
The light rectangles show the intervals $K_{j}$ separated by maximal external
intervals.
The dark rectangles form the sequences of neighbor walls $W_{jk}$ spanning the
intervals $K_{j}$.}
   \end{figure}
First we will use only elements of $\cM$; however,
later, walls of $\cM^{*}$ will be added.
Let $\tup{E_{1},E_{2},\dots}$ be the (finite or infinite) sequence of 
maximal external intervals of size $\ge \bub$, and let us add to it the
maximal external interval starting at $-1$.
Let $K_{1},K_{2},\dots$ be the intervals between them (or possibly after
them, if there are only finitely many $E_{i}$).
Clearly each dominant wall has one of the $K_{j}$ as body.
If there is both a dominant light wall and a dominant heavy wall
with the same body, then we will take the light one as part of the sequence.

By Condition~\ref{cond:distr}.\ref{i:distr.cover}
of $\cM$, each $K_{j}$ that is not a wall
can be spanned by a sequence of neighbors $W_{jk}$.
Each pair of these neighbors will be closer than $3\bub$ to each other.
Indeed, each point of the hop between them belongs either to a 
wall intersecting one of the neighbors,
or to a maximal external interval of size $<\bub$,
so the distance between the neighbors is $< 2\bub + \bub = 3\bub$.
The union of these sequences is a single infinite pure 
sequence of neighbor walls 
  \begin{align}\label{eq:M-seq}
    \bU = \tup{U_{1},U_{2},\dots},  &&\Body(U_{j})=\rint{a_{j}}{b_{j}}.
  \end{align}
Every wall of $\cM$ intersects an element of $\bU$.

A light wall in this sequence is called
\df{isolated} if its distance from other elements of the sequence is
greater than $\f$.
By our construction, all isolated light walls of the sequence $\bU$ are
dominant.
 \end{prooof}

Let us change the sequence $\bU$ using the sequence 
$\tup{W_{1},W_{2},\dots}$ of all emerging walls (disjoint by definition)
as follows.
For $n=1,2,\dots$, add $W_{n}$ to $\bU$.
If $W_{n}$ intersects an element $U_{i}$ then delete $U_{i}$.
\begin{figure}
  \centering
    \includegraphics[scale=0.5]{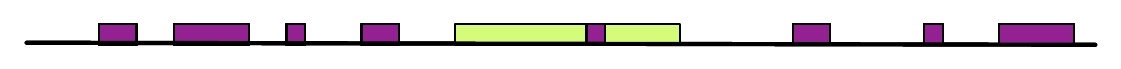}
  \caption{Adding an emerging wall to the pure sequence}
\end{figure}

 \begin{step+}{step:cover.emerg}
 \begin{alphenum}
  \item\label{i:cover.emerg.pure}
 The result is a pure sequence $\bU$ containing all the emerging walls.
  \item\label{i:cover.emerg.adding}
 When adding $W_{n}$, if $W_{n}$ intersects an element $U_{i}$
then $U_{i}$ is a dominant wall of $\cM$ contained in $W_{n}$, and
$W_{n}$ intersects no other element $U_{j}$.
 \end{alphenum}
 \end{step+}
 \begin{pproof}
The proof is by induction.
Suppose that we have already processed $W_{1},\dots,W_{n-1}$, and we
are about to process $W=W_{n}$.
The sequence will be called $\bU$ before processing $W$ and
$\bU'$ after it.

Let us show~\eqref{i:cover.emerg.adding} first.
By the requirement~\eqref{i:emerg.hop} on emerging walls,
either $W$ is an external hop of $\cM$ or it is the union of a dominant
light wall and one
or two external hops of $\cM$, of size $\ge\bub$, surrounding it.
If $W$ is an external hop then it intersects no elements of $\bU$.
Otherwise, the dominant light wall inside it can only be one of the $U_{i}$.

Let us show now~\eqref{i:cover.emerg.pure}, namely that if $\bU$ is pure
then so is $\bU'$.
Property~\eqref{i:pure.external} of the definition of purity
follows immediately, since the intervals
between elements of $\bU'$ are subintervals of the ones
between elements of $\bU$.
For the same reason, these intervals do not contain walls of $\cM$.
It remains to show that if $I'_{j-1}=\rint{b'_{j-1}}{a'_{j}}$ and 
$I'_{j}=\rint{b'_{j}}{a'_{j+1}}$ are
the intervals around $W$ in $\bU'$ then
$a'_{j}$ is clean  in $I'_{j-1}$ and $b'_{j}$ is clean in $I'_{j}$.
Let us show that, for example, $a'_{j}$ is clean in $I'_{j-1}$.

By the requirement~\eqref{i:emerg.outside} on emerging walls,
$a'_{j}$ is adjacent from the right to either an external hop of size $\ge\bub$
or a wall $W'$ of $\cM$.
If the former case, it is left clean and therefore clean in $I'_{j-1}$.
By the definition of emerging walls, $a'_{j}$ is adjacent from the left to
either a dominant wall of $\cM$ or an
external interval $J$ of size $\ge\bub$.
The former is impossible now, since the definition of dominance
excludes the presence of the adjacent wall $W'$ on the left of $a'_{j}$.
The existence of the external interval $J$ along with $W'$
implies the existence of a wall $W''$
in the original sequence $\bU$ whose right end is $a'_{j}$, and this shrinks
interval $I'_{j-1}$ to nothing.
 \end{pproof}

 \begin{figure}
   \centering
       \includegraphics[scale=0.5]{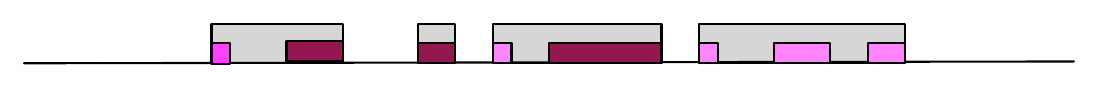}
   \caption{Forming compound walls}
 \end{figure}

 \begin{step+}{step:candid.compound}
Let us break up the pure sequence $\bU$ containing all the emerging walls
into subsequences separated by its intervals $I_{j}$ of size $>\f$.
Consider one of these (possibly infinite) 
sequences, call it $W_{1},\dots,W_{n}$, which is not
just a single isolated light wall.
 \end{step+}
 \begin{prooof}
We will create a sequence of consecutive neighbor walls $W'_{i}$
of $\cM^{*}$ spanning the same interval as $W_{1},\dots,W_{n}$.
In the process, all non-isolated light walls of the sequence 
will be incorporated into a compound wall.

Assume that $W_{i}$ for $i < j$ have been processed already,
and a sequence of neighbors $W'_{i}$ for $i < j'$ has been created
in such a way that 
 \[
   \bigcup_{i<j} W_{i} \sbsq \bigcup_{i<j'} W'_{i},
 \]
and $W_{j}$ is not a light wall which is the last in the series.
(This condition is satisfied when $j=1$ since we assumed that our sequence
is not an isolated light wall.)
We show how to create $W'_{j'}$.

If $W_{j}$ is the last element of the series then it is heavy, and we set
$W'_{j'}=W_{j}$.
Suppose now that $W_{j}$ is not last.

Suppose that it is heavy.
If $W_{j+1}$ is also heavy, or light but not last then $W'_{j'}=W_{j}$.
Else $W'_{j'} = W_{j}+W_{j+1}$, and $W'_{j}$ replaces
$W_{j},W_{j+1}$ in the sequence.
In each later operation also, the introduced new wall will replace
its components in the sequence.

Suppose now that $W_{j}$ is light: then it is not last.
If $W_{j+1}$ is last or $W_{j+2}$ is heavy then
$W'_{j'} = W_{j} + W_{j+1}$.

Suppose that $W_{j+2}$ is light.
If it is last then $W'_{j'} = (W_{j} + W_{j+1}) + W_{j+2}$;
otherwise, $W'_{j'} = W_{j} + W_{j+1}$.
 \end{prooof}

Remove all isolated light walls from $\bU$ and 
combine all the subsequences created in part~\ref{step:candid.compound} above
into a single infinite sequence $\bU$ again.
Consider an interval $I$ between or before its elements.
Then $I$ is inner clean for $\cM$, and the only walls of $\cM$ in $I$ are covered by
some isolated dominant light walls at distance at least $\f/3$ from the endpoints.
Thus, $I$ is inner clean in $\cM^{*}$.
It does not contain any compound walls either (other than possibly those inside
some dominant light wall that was removed), and by definition
it does not contain emerging walls.
Therefore it is a hop of $\cM^{*}$.

 \begin{step+}{step:candid.inner-clean}
Condition~\ref{cond:distr}.\ref{i:distr.inner-clean} holds for $\cM^{*}$.
 \end{step+}
 \begin{pproof}
Let $J$ be a maximal external interval $J$ of $\cM^{*}$, of size
$\ge\bub^{*}$ or starting at $-1$.
Since $J$ has size $\ge\bub^{*}>\f$ it
is an interval separating two elements of $\bU$ and as such is a hop of $\cM^{*}$.
Otherwise it is the interval $I_{0}$.
We have seen that this is also a hop of $\cM^{*}$.
 \end{pproof}

 \begin{step+}{step:candid.cover}
Condition~\ref{cond:distr}.\ref{i:distr.cover} holds for $\cM^{*}$.
 \end{step+}
 \begin{pproof}
By our construction, a maximal external interval of size $\ge\bub^{*}>\f$
is an interval separating two elements of $\bU$.
The segment between two such intervals (or one such and $I_{0}$)
is spanned by elements of $\bU$, separated by hops of $\cM^{*}$.
 \end{pproof}
 \end{Proof}

 \begin{lemma}\label{lem:new-hop}
Suppose that interval $I$ 
contains no walls of $\cM^{*}$, and no wall of $\cM$ closer to its ends than
$\f/3$ (these conditions are satisfied if it is a hop of $\cM^{*}$).
Then it either contains no walls of $\cM$ or all walls of $\cM$ in it are
covered by a sequence $W_{1},\dots,W_{n}$ of dominant 
light neighbor walls of $\cM$ separated from each other 
by external hops of $\cM$ of size $>\f$.

If $I$ is a hop of $\cM^{*}$ then either it is also a hop of $\cM$ or the above
end intervals are hops of $\cM$.
 \end{lemma}
\begin{proof}
If $I$ contains no walls of $\cM$ then there is nothing to prove.
Otherwise, let $U$ be the union of all walls of $\cM$ in $I$.
By assumption, it is not closer than $\f/3$ to the ends of $I$.
Let intervals $J,K$ separate $U$ from both ends, then there are
maximal external intervals $J'$ and $K'$ of size $\ge\bub$ of $\cM$, 
adjacent to $U$ on the left and right.
Let $\tup{E_{1},E_{2},\dots,E_{n}}$ be the sequence of 
maximal external intervals of size $\ge \bub$ in $U$,
and let $E_{0}=J'$, $E_{n+1}=K'$.
Let $F_{i}$ be the interval between $E_{i-1}$ and $E_{i}$.

We claim that all intervals $F_{i}$ are dominant light walls separated by a
distance greater than $\f$.
Note that if $F_{i}$ intersects a dominant light wall $L$ then $F_{i}=L$.
Indeed, 
$L$ is surrounded by maximal external intervals of size $\ge \bub$ which 
then must coincide with $E_{i-1}$ and $E_{i}$.
There is a sequence of neighbor walls spanning $F_{i}$.
If any one is heavy then, since $I$ contains no balls of $\cM^{*}$, it
must be contained in a dominant light wall and then by the observation above,
$F_{i}$ is a dominant light wall.
Assume therefore that  all elements of the sequence are light.
They cannot be farther than $\f$ from each other, since then 
the interval between them would contain an external interval of size $\ge\bub$.
If there is more than one then any two neighbors form a compound wall, which the
absence of walls of $\cM^{*}$ forces to be part of a dominant light wall.

The hops $E_{j}$ have size $>\f$ since
otherwise two neighbors would combine into a compound wall (which could not be
covered by a dominant light wall).

Now suppose $I$ is a hop of $\cM^{*}$.
If it contains no walls of $\cM$ then it is clearly a hop of $\cM$.
Otherwise, look at an end interval, say $J$.
Its right end is also the right end of a maximal external subinterval $J'$, hence it
is clean in $J$.
Since  $I$ is a hop of $\cM^{*}$, the left end of $J$ is also clean in $J$.
So, $J$ is a hop of $\cM$, and the same holds for $K$.
 \end{proof}

The following lemma shows that an emerging barrier
in a ``nice neighborhood'' implies an emerging wall there.

 \begin{lemma}\label{lem:emerg-wall}
Let us be given intervals $I' \subset I$, and also $x(I)$,
with the following properties for some $j\in\{1,2,3\}$.
 \begin{alphenum}
  \item All walls of $\cM$ in $I$ are covered by a sequence $W_{1},\dots,W_{n}$
of dominant light neighbor walls of $\cM$ such that the $W_{i}$ are at a
distance $>\f$ from each other and at a distance $\ge\f/3$ from the ends of
$I$.
  \item  $I'$ is an emerging barrier of type $j$.
  \item\label{i:emerg-wall.away}
$I'$ is at a distance $\ge L_{j}+7\bub$ from the ends of $I$.
 \end{alphenum}
  Then $I$ contains an emerging wall.
 \end{lemma}
 \begin{figure}
   \centering
 \includegraphics[scale=0.5]{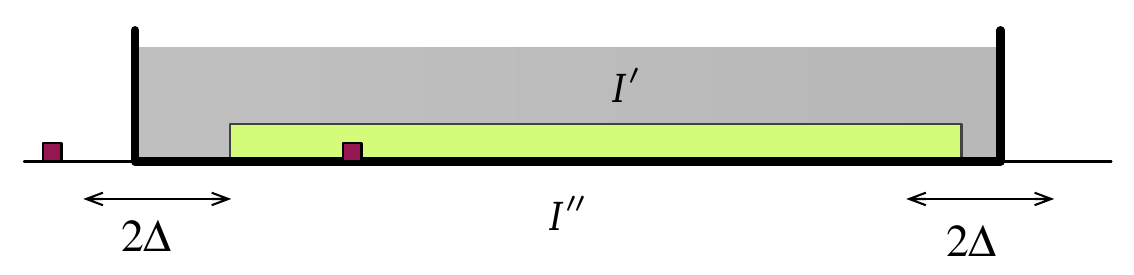}  
      \caption{To Lemma~\protect\ref{lem:emerg-wall}: finding an emerging wall}
\end{figure}
 \begin{proof}
By the definition of emerging barriers, $I'$ contains an emerging barrier of 
size exactly $L_{j}$.
From now on, we assume $I'$ has this size.

Let $I=\rint{a}{b}$, $I'=\rint{u'}{v'}$.
We will define an emerging wall $I''=\rint{u''}{v''}$.
The assumptions imply that the intervals between the walls $W_{i}$ are external
hops. 
However, the interval $\rint{a}{c}$ between the left end of $I$ and $W_{1}$ may
not be one. 
Let $\rint{\hat a}{c}$ be a maximal external subinterval of $\rint{a}{c}$ ending at $c$.
Then $\hat a-a\le\bub$.
Let us define $\hat b$ similarly on the right end of $I$, 
and let $\hat I=\rint{\hat a}{\hat b}$.
We will find an emerging wall in $\hat I$, so let us simply redefine $I$ to be $\hat I$.
We now have the property that the interval between the left end of $I$ and
$W_{1}$ is an external hop of size $\ge\f/3-\bub$, and similarly at the right end.
Also, $I'$ is at a distance $\ge L_{j}+6\bub$ from the ends of $I$.

Assume first that $I$ is a hop of $\cM$ (by the assumption, an external one).
Let us define the interval $I''$ as follows.  
Assumption~\eqref{i:emerg-wall.away} implies $u' \ge a + 2\bub$.
Then, since no wall is contained in $\rint{u'-2\bub}{u'+\bub}$, 
by Condition~\ref{cond:distr}.\ref{i:distr.clean.1}, there is a point 
$u''\in \rint{u'-\bub}{u'}$ clean in $\cM$.
(Since $|I'| > \bub$, there is no problem with walls on the right of $v'$
when finding clean points on the left of $u'$.)
Similarly, $b-v' \ge 2\bub$, and there is a point 
$v''\in \rint{v'}{v'+\bub}$ clean in $\cM$.

Assume now that $I$ is not a hop of $\cM$:
then $I$ is spanned by a nonempty extended sequence $W_{1},\dots,W_{n}$
of neighbor walls of $\cM$ such that the $W_{i}$ are at a distance
 $>\f$ from each other and at a distance $>\f/3-\bub$ from the ends of $I$.

Assume that $I'$ falls into one of the hops, let us call this hop $\rint{a'}{b'}$.
If $u' \ge a' + 2\bub$ then, just as in the paragraph above in interval $I=\rint{a}{b}$,
there is a point $u''\in \rint{u'-\bub}{u'}$ clean in $\cM$.
Otherwise, set $u''=a'$.
Similarly, if $b'-v' \ge 2\bub$, then there is a point 
$v''\in \rint{v'}{v'+\bub}$ clean in $\cM$, else we set $v''=b'$.

Assume now that $I'=\rint{u'}{v'}$ intersects one of these walls, say
$W_{i}=\rint{c}{d}$. 
Now, if $c \le u' < d$ then take $u''=c$.
If  $u' < c$ then there are no walls in the interval $\rint{u'-3\bub}{u'}$,
since it is in the hop on the left of $W_{i}$.
Find a point $u''$ clean in $\cM$ in the middle $\rint{u'-2\bub}{u'-\bub}$ of this interval.
The point $v''$ is defined similarly.

By this definition, interval $I''$
satisfies both requirements~\eqref{i:emerg.hop} and~\eqref{i:emerg.outside}
of emerging pre-walls, and is at a distance $\ge 4\bub+L_{j}$ from the ends of $I$.

If $I$ contains no emerging walls then, in particular, it contains no
walls of type $i$ with $L_{i}\le L_{j}$.
Since $I''$ is at a distance $\ge 4\bub+L_{j}$ (the bound on the size of
emerging walls of type $j$) from the ends of $I$, 
it follows therefore that no wall of such type $i$ intersects it.
But then the process of designating walls in Step~\ref{defstep:emerg} of the
scale-up construction would designate $I''$, or some other interval intersecting it,
a wall, contrary to the assumption that $I$ contains no emerging walls.
 \end{proof}

 \begin{lemma}\label{lem:if-not-missing-hole}
Let the rectangle $Q$ with $X$ projection $I$ contain no traps or vertical walls
of $\cM^{*}$, 
and no vertical wall of $\cM$ closer than $\f/3$ to its sides.
Let $I' = \clint{a}{a + \g}$, $J = \clint{b}{b + 3\bub}$
with $I' \times J \sbsq Q$ be such that $I'$ is 
at a distance $\ge\g+7\bub$ from the ends of $I$.
Suppose that a light horizontal wall $W$ starts at position $b + \bub$.
Then $\clint{a+\bub}{a+\g-\bub}$ 
contains a vertical hole passing through $W$ that is good in the
sense of Definition~\ref{def:holes}.
The same holds if we interchange horizontal and vertical.
 \end{lemma}
 \begin{proof}
Suppose that this is not the case.
Then event $\cL_{3}(X,Y,I',b)$ holds, as defined in the introduction of missing-hole 
traps in Step~\ref{defstep:missing-hole} of the scale-up construction.
Now, if inequality~\eqref{eq:missing-hole-trap-ineq} holds then $I'\times J$
is a trap of the missing-hole kind; but this was excluded, since $Q$
contains no traps of $\cM^{*}$.
On the other hand, if~\eqref{eq:missing-hole-trap-ineq} does not hold 
then (due also to Lemma~\ref{lem:new-hop})
Lemma~\ref{lem:emerg-wall} is applicable to the interval $I'$ and
the interval $I$ that is the $X$ projection of $Q$, and we can conclude
that $I$ contains a vertical emerging wall.
But this was also excluded.
 \end{proof}

 \begin{lemma}\label{lem:if-correl}
Let rectangle $Q$ with $X$ projection $I$ contain no traps or vertical walls of
$\cM^{*}$, and no vertical walls of $\cM$ closer than $\f/3$ to its sides.
For $j\in\{1,2\}$, let $\l_{j}$ be as introduced in the definition of
correlated traps and emerging walls in Steps~\ref{defstep:correl}
and~\ref{defstep:emerg} of the scale-up construction.
Let $I' = \clint{a}{a + L_{j}}$,
$J = \clint{b}{b + 5\bub}$ with $I' \times J \sbsq Q$ be such that $I'$ is 
at a distance $\ge L_{j}+7\bub$ from the ends of $I$.
Then there is an interval $I'' \sbsq I'$ of size $\l_{j}$,
such that the rectangle $I'' \times J$ contains no trap of $\cM$.
The same holds if we interchange horizontal and vertical.
 \end{lemma}
 \begin{proof}
The proof of this lemma is completely analogous to the proof of
Lemma~\ref{lem:if-not-missing-hole}.
 \end{proof}

 \begin{lemma}\label{lem:clean}
The new mazery $\cM^{*}$ defined by the above construction satisfies
Conditions~\ref{cond:distr}.\ref{i:distr.clean.1}
and~\ref{cond:distr}.\ref{i:distr.clean.2}.
 \end{lemma}
\begin{Proof}
  \begin{step+}{step:clean.1dim}
Let us prove Condition~\ref{cond:distr}.\ref{i:distr.clean.1}.
  \end{step+}
   \begin{prooof}
Consider an interval $I$ of size $3 \bub^{*}$ containing no walls of
$\cM^{*}$.
Condition~\ref{cond:distr}.\ref{i:distr.cover} says that
the real line is spanned by an extended sequence
$\tup{W_{1},W_{2},\dots}$ of neighbor 
walls of $\cM$ separated from each other by hops of $\cM$.
As shown in the construction of part~\ref{step:cover.initial-pure} of the proof of
Lemma~\ref{lem:cover}, we can also assume that every dominant light wall is an
element of this sequence.
If any of these walls $W$ is contained in $I$ then it is light.
Indeed, $I$ contains no wall of $\cM^{*}$, so $W$ can only be heavy if it is
contained in a dominant light wall $W'$, but then it is equal to $W'$, 
as we assumed.

Since $I$ contains no wall of $\cM^{*}$, if two of these walls fall into
$I$ then they are separated by a hop of size $>\f$.

Let $I'$ be the middle third of $I$.
Then $|I'| \ge 2\f + \bub$, and removing the $W_{i}$ from $I'$ leaves a
subinterval $\rint{a}{b} \sbsq I'$ of size at least $\f$.
(If at least two $W_{i}$ intersect $I'$ take the interval between
consecutive ones, otherwise $I'$ is divided into at most two pieces of
total length at least $2\f$.)
Now $K = \rint{a+\bub+\f/3}{b-\bub-\f/3}$ is an interval of length at least
$\f/3-2\bub > 3\bub$ which has distance at least $\f/3$ from any wall of $\cM$. 
There will be a clean point in the middle of $K$ which will then be clean
in $\cM^{*}$.
   \end{prooof}

 \begin{step+}{step:clean.2dim}
Let us prove Condition~\ref{cond:distr}.\ref{i:distr.clean.2} now for $\cM^{*}$.
 \end{step+} 
 \begin{prooof}
We will confine ourselves to the statement in which the point $a$ is
assumed clean and we find a $b$ such that the point $\pair{a}{b}$ is clean.
The half clean cases are proved similarly.
Let $I,J$ be right-closed intervals of size $3\bub^{*}$,
suppose that the rectangle $I \times J$ contains no traps or horizontal walls
of $\cM^{*}$, and $a$ is a point in the middle third of $I$ that is
clean in $\cM^{*}$ for $X$.
We need to prove that 
there is an integer $b$ in the middle third of $J$ such that the point
$\pair{a}{b}$ is clean in $\cM^{*}$.

Just as in Part~\ref{step:clean.1dim} above, we find $K$ 
with $\f/3-2\bub \le |K|\le \f$ in the middle of $J$ 
which is at distance at least $\f/3$ from any horizontal wall of $\cM$.
Let $I' = \rint{a-\g-\bub}{a+\g+\bub}$, then $I'\sbsq I$.
We will find an interval $K''\sbsq K$ with $|K''|\ge\g'$
such that $I'\times K''$ contains no trap of $\cM$.
If there are no traps in $I' \times K$ let $K''=K$.
Assume now that $I' \times K$ contains a trap $T = \Rect(u,v)$ of $\cM$,
where $u=\pair{u_{0}}{u_{1}}$, $v=\pair{v_{0}}{v_{1}}$.
Since we assumed there are no traps of $\cM^{*}$ and thus no
uncorrelated traps, any trap must meet either
$\clint{u_{0}}{v_{0}} \times K$ or $I' \times \clint{u_{1}}{v_{1}}$ or be
at a distance at least $\f$ from $T$ (and hence outside $I'\times K$).
Let $K'$ be a subinterval of $K \setminus \clint{u_{1}}{v_{1}}$
of size $4\g'$
(which exists since $|K| \ge 2\cdot(4\g')+\bub$).
By Lemma~\ref{lem:if-correl}, there must exist a subinterval 
$K''$ of $K'$ of length $g' \ge 2\g + 3\bub$ such that
$\clint{u_{0}-2\bub}{u_{0}+3\bub} \times K''$ contains no trap.
Then also $I' \times K''$ contains no trap.

Now restrict $K''$ to an interval of size $3\bub$ in its middle
and then find a clean point $\pair{a}{b}$ in its middle third
applying Condition~\ref{cond:distr}.\ref{i:distr.clean.2} for $\cM$.
Then $\pair{a}{b}$ has distance at least $\g+\bub$ from the boundary of
$I'\times K''$ and so has distance at least $\g$ from any trap.
Since $b$ is at distance at least $\f/3$ from any wall, it is clean in
$\cM^{*}$.
Hence $\pair{a}{b}$ is clean in $\cM^{*}$.
 \end{prooof}
 \end{Proof}

\section{Probability bounds}\label{sec:bounds}

In this section, we derive all those bounds on probabilities in $\cM^{k}$ that are
possible to give without indicating the dependence on $k$.

\subsection{General bounds}\label{subsec:bounds}

Recall the definitions needed for the hole lower bound
condition, Condition~\ref{cond:distr}.\ref{i:distr.hole-lb}, in particular the
definition of the event $E$.
Since $(c-1)$ will be used often, we denote it by $\hat c$.
Let $u \le v < w$, and $a$ be given with $v-u \le 12\bub$, and define
$b = a + \cei{\frac{v-u}{2}}$, $c = a + (v-u) + 1$.
We need to extend the lower bound condition in several ways.
Since we will hold the sequence $y$ of values of the sequence $Y$ of random
variables fixed in this subsection, we take the liberty
and omit the condition $Y = y$ from the probabilities: it is always assumed
to be there.
For the following lemma, remember
Condition~\ref{cond:distr}.\ref{i:distr.ncln}. 

 \begin{lemma}\label{lem:hole-lb-clean}
Let $F_{t}$ be the event that the point $\pair{t}{w}$ is upper right rightward H-clean.
Let $\hat E$ be the event that $E$ is realized with a hole $\rint{d}{t}$,
and $F_{t}$ holds (that is the hole is good as seen from $\rint{a}{u}$, in the
sense of Definition~\ref{def:holes}).
We have
 \begin{equation}\label{eq:hole-lb-clean}
   \Prob(\hat E) \ge (1-\ncln) \Prob(E).
 \end{equation}
 \end{lemma}
 \begin{proof}
For $b \le t \le c+\bub$, let $E_{t}$ be the event
that $E$ is realized by a hole ending at $t$ but is not realized
by any hole ending at any $t' < t$.
Then $E = \bigcup_{t}E_{t}$, 
$\hat E \spsq \bigcup_{t} (E_{t} \cap F_{t})$.
Due to the Markov chain property of $X$ and the form of $E_{t}$, the fact that
$E_{t}$ depends only on $X(0),\dots,X(t)$ and using
inequalities~\eqref{eq:ncln.1dim} and~\eqref{eq:ncln.lower} of
Condition~\ref{cond:distr}.\ref{i:distr.ncln}, we have
 \[
 \Prob(E_{t} \cap F_{t}) = \Prob(E_{t})\Prob(F_{t} \mid E_{t}) 
  \ge \Prob(E_{t})(1-\ncln).
 \]
The events $E_{t}$ are mutually disjoint. 
Hence
 \[
   \Prob(\hat E) \ge \sum_{t} \Prob(E_{t} \cap F_{t})
            \ge (1-\ncln) \sum_{t} \Prob(E_{t}) 
             = (1-\ncln) \Prob(E).
 \]
 \end{proof}

Recall Remark~\ref{rem:distr}.\ref{i:rem.hole-lb}, referring to
the most important special case of the hole lower bound:
for any horizontal wall $B$ of rank $\r$, at any point $b$,
the probability that there is a vertical hole passing through $B$ at point
$b$ is at least $\h(\r)$.
We strengthen this observation in a way similar to
Lemma~\ref{lem:hole-lb-clean}.

 \begin{lemma}\label{lem:hole-lb-clean-2}
Let $v < w$, and let us fix the value $y$ of the sequence of random variables
$Y$ in such a way that
there is a horizontal wall $B$ with body $\rint{v}{w}$.
Let point $b$ be given.
Let $E$ be the event that a good hole $\rint{b}{b'}$ passes through $B$
(this event still depends on 
the sequence $X=\tup{X(1),X(2),\dots}$ of random variables).
Let $\s\in\{1,\dots,\m\}$ then
 \[
   \Pbof{E \mid X(b-\bub)=\s} \ge (1-\ncln)^{2}\h(\r).
 \]
 \end{lemma}
 \begin{proof}
Condition~\ref{cond:distr}.\ref{i:distr.ncln} implies that
the probability that point $b$ is lower left H-clean is lower-bounded by $(1-\ncln)$.
Conditioning on times during and before this event, Lemma~\ref{lem:hole-lb-clean} 
lower-bounds the probability that $\rint{b}{b'}$ is an upper right rightward
H-clean hole.
The lower bound is then the product of these two lower bounds.
 \end{proof}

Now, we prove a version of the hole lower bound condition
that will help proving the same bound for $\cM^{*}$.
This is probably the only part of the paper in which the probability
estimates are somewhat tricky.

 \begin{definition}\label{def:E-star}
Recall the definition of event $E$ in
Condition~\ref{cond:distr}.\ref{i:distr.hole-lb}, and that it refers to a
horizontal wall with body $\rint{v}{w}$ seen from a point $\pair{a}{u}$.
Take the situation described above, possibly without the bound on $(v-u)$.
Let 
 \begin{equation}\label{eq:E-star}
  E^{*} = E^{*}(u, v, w;\, a)
 \end{equation}
be the event (a function of the sequence $X$)
that there is a $d \in \clint{b}{c-1}$ with the following properties for
$Q=\Rect^{\rightarrow}(\pair{a}{u}, \pair{d}{v})$:
 \begin{enumerate}[\upshape (i*)]
  \item A vertical hole (of $\cM$) fitting $B$ starts at $d$.
  \item $Q$ contains no traps or vertical barriers of $\cM$ or $\cM^{*}$ and
is inner H-clean for $\cM^{*}$.
 \end{enumerate}
The difference between $E^{*}(\cdot)$ and $E(\cdot)$ is only
that $E^{*}$ requires the H-cleanness for $\cM^{*}$ and also absence of barriers
and traps for $\cM^{*}$ whenever possible.
 \end{definition}

 \begin{definition}[Barrier and trap probability upper bounds]\label{def:pub-intr}
Let 
 \begin{equation}\label{eq:pub-intr}
  \pub
 \end{equation}
be an upper bound of the probabilities over all possible points
$a$ of the line, and over all possible values of $X(a)$,
that a barrier of $\cM$ starts at $a$.
Let it also bound similarly the probability 
that a barrier of $\cM^{*}$ starts at $a$.
Let
 \begin{equation}\label{eq:tubh-intr}
   \tubh
 \end{equation}
be an upper bound of the conditional probabilities 
over $X$ (with $Y$ and $X(a-1)$ fixed in any possible way) 
over all possible points $\pair{a}{b}$ of the plane,
that a trap of $\cM$ starts at $\pair{a}{b}$.
Let it also bound similarly the probability that a
trap of $\cM^{*}$  starts there.
 \end{definition}

 \begin{lemma}\label{lem:hole-lb-2}
Suppose that the requirement $v-u \le 12\bub$ 
in the definition of the event $E^{*}$ is replaced with $v-u\le 12\bub^{*}$,
while the rest of the requirements are the same.
Let us fix $X(a)=\s$ arbitrarily.
Then the inequality
 \begin{equation}\label{eq:hole-lb-3}
   \Pbof{E^{*}\mid X(a)=\s,\;Y=y} \ge 0.5 \land ((c-b)^{\hxp} \h(\r)) - U
 \end{equation}
holds, with $U = 26 \pub \bub^{*} + 338 \tubh(\bub^{*})^{2}$.
If $v-u>12\bub$ then we also have the somewhat stronger inequality
 \begin{equation}\label{eq:hole-lb-2}
   \Pbof{E^{*}\mid X(a)=\s,\;Y=y} \ge 0.5 \land (1.1(c-b)^{\hxp} \h(\r)) - U.
 \end{equation}
 \end{lemma}
  \begin{figure}
\centering
   \includegraphics[scale=0.95]{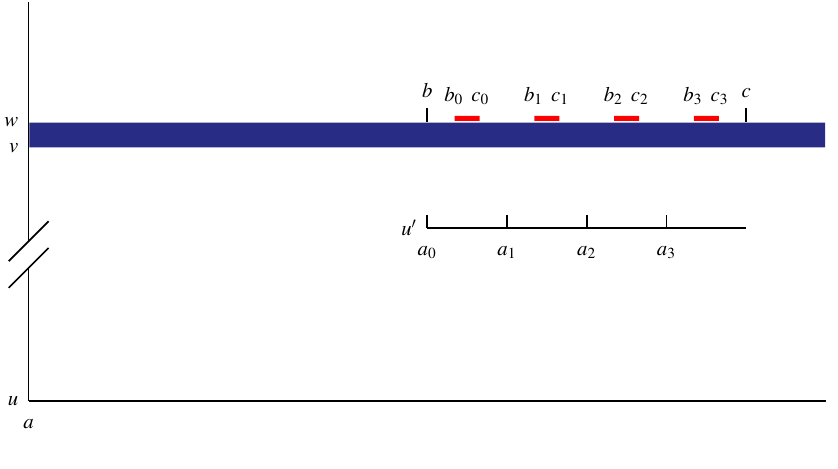}
 \caption{To the proof of Lemma~\protect\ref{lem:hole-lb-2}.}
 \end{figure}

 \begin{Proof}
For ease of reading, we will omit the conditions $X(a)=\s$, $Y=y$ from the
probabilities. 

For the case $v-u\le 12\bub$, even the original stronger inequality holds,
namely $\Prob(E^{*}) \ge (c-b)^{\hxp} \h(\r)$.
Condition~\ref{cond:distr}.\ref{i:distr.hole-lb} implies this already for
$\Prob(E)$, so it is sufficient to show $E\sbsq E^{*}$ in this case.

As remarked after its definition, the event $E^{*}$ differs from $E$ only in
requiring that rectangle $Q$ contains no traps or vertical barriers
of $\cM^{*}$, not only for $\cM$, and that 
points $\pair{a}{u}$ and $\pair{d}{v}$ are H-clean in $Q$ for $\cM^{*}$ also, not
only for $\cM$.

Consider a trap of $\cM^{*}$ in $Q$: this cannot be an uncorrelated or
correlated trap, since its components traps, which are traps of $\cM$, are
already excluded.
It cannot be a trap of the missing-hole kind either, since that trap is too big
for $Q$ when $v-u\le 12\bub$.
The same argument applies to vertical barriers of $\cM^{*}$.
The components of the compound barriers that belong to $\cM$
are excluded, and the emerging barriers are too big.

These considerations take care also of the issue of H-cleanness
for $\cM^{*}$, since the latter also boils down to the absence of traps and
barriers.

Let us move to the case of $v-u>12\bub$, which implies $\hat c - b \ge 6 \bub$.
We will use the following inequality, which can be checked by direct calculation.
Let $\ag = 1 - 1/e = 0.632\dots$, then for $x > 0$ we have
 \begin{equation}\label{eq:exp-lb}
   1 - e^{-x} \ge \ag \land \ag x.
 \end{equation}
Let $n = \flo{(c - b)/(3 \bub)}$, then $n\ge 2$ and hence
$(c-b)/(3\bub)\le n+1 \le 1.5 n$, implying
 \begin{equation}\label{eq:n-lb}
  n \bub \ge (c - b ) /4.5.
 \end{equation}
Let
 \begin{align*}
       u'    &= v - 2 \bub,
\quad   a_{i} = b + 3 i \bub,
\\    E'_{i} &= E(u', v, w;\, a_{i})
\quad \txt{ for } i = 0, \dots, n-1,
\\     E'  &= \bigcup_{i} E'_{i}.
 \end{align*}
Let $C$ be the event that point $\pair{a}{u}$ is upper right rightward H-clean in
$\cM$.
Then by Conditions~\ref{cond:distr}.\ref{i:distr.ncln}
 \begin{align}\label{eq:P(neg C)}
 \Prob(\neg C) \le 2(\ncln/2)\le 0.1.
 \end{align}
Let $D$ be the event that the rectangle $\rint{a}{c} \times \clint{u}{v}$
contains no trap or vertical barrier of $\cM$ or $\cM^{*}$.
(Then $C\cap D$ implies that $\pair{a}{u}$ is also upper right rightward H-clean
in the rectangle $\rint{a}{c} \times \clint{u}{v}$
in $\cM^{*}$.)
By definition,
 \begin{align*}
  \Prob(\neg D) &\le 2\pub(c-a) + 2\tubh(c-a)(v-u+1)
\\  &\le 2\cdot 13\pub\bub^{*} + 2\cdot 13 \cdot 13\tubh(\bub^{*})^{2}
      = 26 \pub \bub^{*} + 338 \tubh(\bub^{*})^{2}.
 \end{align*}

 \begin{step+}{step:hole-lb-2.CD}
Let us show $C \cap D \cap E' \sbsq E^{*}(u, v, w;\, a)$.
 \end{step+}
 \begin{prooof}
Indeed, suppose that $C \cap D \cap E'_{i}$ holds with some hole starting
at $d$. 
Then there is a 
rectangle $Q'_{i} = \Rect^{\rightarrow}(\pair{a_{i}}{u'}, \pair{d}{v})$
containing no traps or vertical barriers of $\cM$,
such that $\pair{d}{v}$ is H-clean in $Q'_{i}$.
It follows from $D$ that the rectangle 
 \[
  Q^{*}_{i} = \Rect^{\rightarrow}(\pair{a}{u}, \pair{d}{v}) \spsq Q'_{i}
 \]
contains no traps or vertical barriers of $\cM$ or $\cM^{*}$.
Since event $C$ occurs, the point $\pair{a}{u}$ is H-clean for $\cM$ in
$Q^{*}_{i}$.
The event $E'_{i}$ and the inequalities $d-a,v-u\ge\bub$ 
imply that $\pair{d}{v}$ is H-clean in $Q^{*}_{i}$,
and a hole passing through the potential wall starts at $d$ in $X$.
The event $D$ implies that
there is no trap or vertical barrier of $\cM$ or in $Q^{*}_{i}$.
Hence $Q^{*}_{i}$ is also inner H-clean in $\cM^{*}$, and so $E^{*}$ holds.
 \end{prooof}

We have $\Prob(E^{*}) \ge \Prob(C)\Prob(E' \mid C) - \Prob(\neg D)$.

 \begin{step+}{step:hole-lb-2.union-prob}
It remains to estimate $\Prob(E' \mid C)$.
 \end{step+}
 \begin{prooof}
Let us denote $\s = \bub^{\hxp}\h(\r)$.
Condition~\ref{cond:distr}.\ref{i:distr.hole-lb} is applicable to
$E'_{i}$, so we have for each $k \in \{1,\dots,\m\}$:
 \begin{align*}
\Pbof{E'_{i} \mid C\cap\{X(a_{i})=k\}}=\Pbof{E'_{i} \mid X(a_{i})=k} \ge \s,   
 \end{align*}
where we could delete the condition $C$ due to the Markov property.
Hence
 \[
 \Pbof{\neg E'_{i} \mid C \cap \{X(a_{i})=k\}} \le 1 - \s \le e^{-\s}.
 \]
Due to the Markov property of the sequence $X$, this implies
$\Prob\bigparen{\neg E'_{i} \mid C\cap\bigcap_{j<i} \neg E'_{j}} 
\le e^{-\s}$, and hence
 \begin{equation}\label{eq:union-E-R}
  \Prob(E' \mid C) = 1 - \Prob\bigparen{\bigcap_{i}\neg E'_{i} \mid C}
                 \ge 1 - e^{- n \s} \ge \ag \land (\ag n \s),
 \end{equation}
where in the last step we used~\eqref{eq:exp-lb}.
By~\eqref{eq:n-lb}:
 \begin{equation*}
 \ag n \bub^{\hxp} = \ag n^{1-\hxp}(\bub n)^{\hxp}
      \ge \ag 2^{1-\hxp}(\bub n)^{\hxp} 
      \ge \ag (2^{1-\hxp}/4.5^{\hxp})(c-b)^{\hxp}
      \ge 1.223 (c-b)^{\hxp},
 \end{equation*}
where we used the value of $\hxp$ from~\eqref{eq:hxp}.
Substituting into~\eqref{eq:union-E-R}:
 \begin{align*}
   \Prob(E' \mid C) &\ge \ag \land (1.223(c-b)^{\hxp}\h(\r)),
\\ \Prob(C)\Prob(E' \mid C) &\ge 0.9 \cdot (\ag \land (1.223 (c-b)^{\hxp} \h(\r)))
                      > 0.5 \land (1.1(c-b)^{\hxp} \h(\r)),
 \end{align*}
where we used~\eqref{eq:P(neg C)}.
 \end{prooof}
 \end{Proof}

\subsection{New traps}

Recall the definition of uncorrelated compound traps in
Step~\ref{defstep:uncorrel} of the scale-up construction in
Section~\ref{sec:plan}.

 \begin{lemma}[Uncorrelated Traps]\label{lem:uncorrel-trap}
Given a string $x = \tup{x(0), x(1),\dots}$, a point $\pair{a_{1}}{b_{1}}$,
% and $v > a_{1}$, 
let $\cF$ be the event that an uncorrelated compound
trap of $\cM^{*}$ starts at $\pair{a_{1}}{b_{1}}$.
% with projection lying in the interval $I=\lint{a_{1}}{v}$.
Let $\s\in\{1,\dots,\m\}$, then
 \begin{equation}\label{eq:uncorrel-trap}
 \Pbof{\cF \mid X = x,\; Y(b_{1}-1)=\s} \le 2 \f^{2} \tub^{2}.
 \end{equation}

 \end{lemma}
 \begin{proof}
Let $\cG(a, b)$ be the event that a trap of $\cM$ starts at $\pair{a}{b}$.
Let $\cG(a, b;\; a', b')$ 
be the event that a trap of $\cM$ starts at $\pair{a}{b}$,
and is contained in $\lint{a}{a'}\times\lint{b}{b'}$.
Since the new trap is the smallest rectangle containing two old traps,
it must contain these in two of its opposite corners:
let $\cE$ be the event that one of these corners is $\pair{a_{1}}{b_{1}}$.

Let $N = \pair{a_{1}}{b_{1}} + \rint{0}{\f}^{2}$.
Then
 \[
 \cE \sbsq \bigcup_{\pair{a_{2}}{b_{2}} \in N} 
  \cG(a_{1}, b_{1}; a_{2}, b_{2}) \cap \cG(a_{2}, b_{2}).
 \]
% For an arbitrary interval $I$, let $\cX(I)$ be the event $X(I) = x(I)$.
% Thus, for any sufficiently large interval $J$ we can write
%  \[
%    \cX(J) \cap \cE \sbsq \bigcup_{\pair{a_{2}}{b_{2}} \in N} 
%  (\cX\lint{a_{1}}{a_{2}} \cap \cG(a_{1}, b_{1};\, a_{2}, b_{2})) 
%   \cap (\cX\clint{a_{2}}{v} \cap \cG(a_{2}, b_{2}; v, \infty)).
%  \]
The events $\cG(a_{1}, b_{1};\, a_{2}, b_{2})$
and $\cG(a_{2}, b_{2})$ belong to rectangles whose projections are disjoint.
Denoting by $\cC$ the event $Y(b_{1}-1)=\s$,
by Condition~\ref{cond:distr}.\ref{i:distr.trap-ub} and the Markov property: 
 \begin{align*}
    \Pbof{\cG(a_{1}, b_{1};\, a_{2}, b_{2}) \mid \{X=x\}\cap \cC}&\le \tub,
\\ \Pbof{\cG(a_{2}, b_{2})\mid \cG(a_{1}, b_{1};\, a_{2}, b_{2}) \cap\{X=x\}\cap
\cC} &\le \tub.
 \end{align*}
Hence by the union bound $\Pbof{\cE \mid \{X=x\}\cap\cC} \le  \f^{2} \tub^{2}$.
If $\cF \setminus \cE$ holds then there is a pair $\pair{A}{B} \in N$
such that $\cG(a_{1}, B;\, A, \infty)$ and $\cG(A, b_{1};\, \infty, B)$ holds.
A computation similar to the above one gives the upper bound 
$\f^{2} \tub^{2}$ for $\Pbof{\cF \setminus \cE \mid \{X=x\}\cap\cC}$.
\end{proof}

Recall the definition of correlated traps in
Step~\ref{defstep:correl} of the scale-up construction in Section~\ref{sec:plan}.

 \begin{lemma}[Correlated Traps]\label{lem:correl}
Let a site $\pair{a}{b}$ be given.
For $j=1,2$, let $\cF_{j}$ be the event that a horizontal correlated trap
of type $j$ starts at $\pair{a}{b}$.
 \begin{alphenum}

  \item\label{i:correl-trap-vert} 
Let us fix a string $x = \tup{x(0), x(1), \dots}$, and 
also $\s\in\{1,\dots,m\}$ arbitrarily.
We have
 \begin{equation}\label{eq:correl-trap-vert}
 \Pbof{\cF_{j} \mid X=x,\; Y(b-1)=\s} \le w^{2}.
 \end{equation}

  \item\label{i:correl-trap-horiz} 
Let us fix a string $y = \tup{y(0), y(1), \dots}$, and 
also $\s\in\{1,\dots,m\}$ arbitrarily.
We have
 \begin{equation}\label{eq:correl-trap-horiz}
 \Pbof{\cF_{j} \mid Y = y,\; X(a)=\s} \le (5\bub\l_{j}\tub)^{4}.
 \end{equation}
 \end{alphenum}
 \end{lemma}
 \begin{proof}
Part~\eqref{i:correl-trap-vert} is an immediate consequence of 
requirement~\eqref{eq:correl-trap-ineq} of the definition of correlated traps.
It remains to prove part~\eqref{i:correl-trap-horiz}.
Note that this result implies the same bounds also if we fix $X(a-1)$ arbitrarily.
If there is a correlated trap with $X$-projection starting at some $a$ then
there must be traps with $X$-projections
in $\rint{a+r\l_{j}}{a+(r+1)\l_{j}}$ for $r=0,1,2,3$. 
Due to Condition~\ref{cond:distr}.\ref{i:distr.trap-ub}
(the trap upper bound) and the Markov property,
the probability of a trap in any one of these is at most $5\bub\l_{j}\tub$,
even conditioned on the values of $X$ before.
Hence the probability of such a compound trap happening is at most
$(5\bub\l_{j}\tub)^{4}$.
 \end{proof}

Recall the definition of traps of the missing-hole kind in
Step~\ref{defstep:missing-hole} of the scale-up algorithm in Section~\ref{sec:plan}.

 \begin{lemma}[Missing-hole traps]\label{lem:missing-hole}
For $a,b \in \bbZ_{+}$, let $\cF$ be the event 
that a horizontal trap of the missing-hole kind starts at $\pair{a}{b}$.
 \begin{alphenum}
  \item\label{i:missing-hole-trap-vert}
Let us fix a string $x = \tup{x(0), x(1), \dots}$, and 
also $\s\in\{1,\dots,m\}$ arbitrarily.
We have
 \begin{equation}\label{eq:missing-hole-trap-vert}
 \Pbof{\cF \mid X=x,\; Y(b-1)=\s} \le w^{2}.
 \end{equation}

  \item\label{i:missing-hole-trap-horiz}
Let us fix a string $y = \tup{y(0), y(1), \dots}$, and 
also $s\in\{1,\dots,m\}$ arbitrarily.
Let $n = \Flo{\frac{\g}{3 \bub}}$.
We have 
 \begin{equation}\label{eq:missing-hole-trap-horiz}
 \Pbof{\cF \mid Y=y,\; X(a)=\s} \le e^{- (1-\ncln)^{2} n \h(\R^{*})}.
 \end{equation}
 \end{alphenum}

 \end{lemma}
 \begin{proof}
Part~\eqref{i:missing-hole-trap-vert} is an immediate consequence of 
requirement~\eqref{eq:missing-hole-trap-ineq} 
of the definition of missing-hole traps.
It remains to prove part~\eqref{i:missing-hole-trap-horiz}.
Note that this result implies the same bounds also if we fix 
$X(a-1)$ arbitrarily.
Let $J = \clint{b}{b + 3\bub}$.
According to the definition of missing-hole traps above, we can assume
without loss of generality that, with $b_{1}=b+\bub$, there
is a $b_{2}>b_{1}$ such that
$\rint{b_{1}}{b_{2}}$ is a potential light horizontal wall $W$.
For $i = 0,\dots,n-1$, let $\cA(d, i)$ be the event that no good hole 
$\rint{a_{1}}{a_{2}}$ with $a_{1}=a+3i\bub+\bub$ passes through $W$.
All these events must hold if
a horizontal trap of the missing-hole kind starts at $\pair{a}{b}$.
Using the Markov property and Lemma~\ref{lem:hole-lb-clean-2}:
 \begin{equation*}
 \Prob\bigparen{\cA(d, i)\mid \bigcap_{j<i}\cA(d, j)} 
\le 1 - (1-\ncln)^{2} \h(\R^{*})\le e^{-(1-\ncln)^{2} \h(\R^{*})}.
  \end{equation*}
Therefore 
$\Prob\bigparen{\bigcap_{i} \cA(d, i)} \le e^{-n (1-\ncln)^{2} \h(\R^{*})}$.
\end{proof}

\subsection{Emerging walls}

Recall the definition of emerging walls in
Step~\ref{defstep:emerg} of the scale-up algorithm in Section~\ref{sec:plan}.

 \begin{lemma}\label{lem:emerg}
For any point $u$, let $\cF(t)$ be the event that a barrier $\rint{u}{v}$ 
of $X$ of the emerging kind, of length $t$ starts at $u$.
Let $k\in\{1,\dots,\m\}$.
We have, with $n = \Flo{\frac{\g}{3 \bub}}$:
 \begin{equation}\label{eq:emerg}
 \sum_{t}\Pbof{\cF(t) \mid X(u)=k} \le 4\m\bub^{2} \tub^{2}
 \bigparen{2\cdot (5\bub\g')^{4} + \tub^{-4} e^{-(1-\ncln)^{2} n \h(\R^{*})}}.
 \end{equation}
 \end{lemma}
 \begin{proof}
For interval $I'=\clint{u'}{v'}$ let event $\cL_{j}(x, Y, I',1)$ be defined as in
Steps~\ref{defstep:correl} and~\ref{defstep:missing-hole}
of the scale-up algorithm in Section~\ref{sec:plan}.
Let us fix an arbitrary $k'\in\{1,\dots,\m\}$.
By the proof of Lemma~\ref{lem:correl}, for $j=1,2$:
 \begin{equation*}
   \Pbof{\cL_{j}(X,Y,I',1) \mid X(u)=k,\; Y(0) = k'}
 \le (5\bub\l_{j}\tub)^{4} =: U_{j},
 \end{equation*}
where we took inequality~\eqref{eq:correl-trap-horiz} and unconditioned on all
$Y(b)$ for $b\ge 1$.
Similarly by the proof of Lemma~\ref{lem:missing-hole}:
 \begin{equation*}
 \Pbof{\cL_{3}(X,Y,I',1) \mid X(u)=k,\; Y(0) = k'}
   \le e^{- (1-\ncln)^{2} n \h(\R^{*})}=: U_{3}.
 \end{equation*}
Let us define the function $\pi(x(I'))$ depending on the sequence $x(I')$
as $\pi(x(I'))=\Pbof{X(I')=x(I')\mid X(u)=k}$.
Then
\begin{align*}
\sum_{x(I')}\pi(x(I'))
 \Pbof{\cL_{j}(X,Y,I',1) \mid X(u)=k,\; X(I')=x(I'),\; Y(0) = k'}\le U_{j}.
  \end{align*}
The Markov inequality implies that the probability (conditioned on $X(u)=k$) of those
$x$ for which the inequality
 \begin{equation*}
   \Pbof{\cL_{j}(X, Y, I',1) \mid X(u)=k,\; X(I')=x(I'),\; Y(0)=k'} > \tub^{2}
 \end{equation*}
holds will be upper-bounded by $\tub^{-2}U_{j}$.

The length of interval $I'=\rint{u'}{v'}$ is defined by the type, but not its
starting point $u'\in u+\lint{0}{2\bub}$, neither is
the endpoint $v\in v'+\lint{0}{2\bub}$ of the barrier. 
For every one of the $(2\bub)^{2}$ possible choices of these,
we obtain a particular length $t$ of the barrier $\rint{u}{v}$.
Multiplying by the number $\m$ of possible choices of $k'$ we obtain an
upper bound on desired sum for an emerging wall of type $j$.
Adding up the three values and recalling $\max(\l_{1},\l_{2})=\g'$ gives
 \[
  4\m\bub^{2}\tub^{-2}(U_{1} + U_{2} + U_{3}) <
  4\m\bub^{2}\tub^{2}\bigparen{2\cdot (5\bub\g')^{4} + 
      \tub^{-4} e^{-(1-\ncln)^{2} n \h(\R^{*})}}.
 \]
 \end{proof}

\subsection{Compound walls}\label{subsec:compound}

Let us use the definition of compound walls given in
Step~\ref{defstep:compound} of the scale-up algorithm of
Section~\ref{sec:plan}.

 \begin{lemma}\label{lem:compound-wall-ub} 
Consider ranks $\r_{1},\r_{2}$ at any stage of the scale-up construction.
Assume that Condition~\ref{cond:distr}.\ref{i:distr.wall-ub} already holds for
rank values $\r_{1},\r_{2}$.
For a given point $x_{1}$, let us fix $X(x_{1})=k$ for some
$k\in\{1,\dots,\m\}$ arbitrarily.
Then the sum, over all $l$, of the probabilities for the occurrence of a 
compound barrier of type $\ang{\r_{1},\r_{2},i}$
and width $l$ at $x_{1}$ is bounded above by
 \begin{equation}\label{eq:compound-wall-ub}
   \lg^{i}\p(\r_{1})\p(\r_{2}).
 \end{equation}
 \end{lemma}
 \begin{proof}
Noting $d_{i+1}-d_{i} \le \lg^{i}$ for all $i$,
we will prove an upper bound $(d_{i+1}-d_{i})\p(\r_{1})\p(\r_{2})$.
For fixed $\r_{1},\r_{2},x_{1},d$, let $B(d,l)$ be the event that a 
compound barrier of any type $\ang{\r_{1},\r_{2},i}$ with
distance $d$ between the component barriers, and size $l$ appears
at $x_{1}$. 
For any $l$, let $A(x,\r,l)$ be the event that a barrier of rank
$\r$ and size $l$ starts at $x$.
We can write 
 \[
   B(d,l) = \bigcup_{l_{1}+d+l_{2}=l} 
    A(x_{1}, \r_{1}, l_{1}) \cap A(x_{1}+l_{1}+d,\, \r_{2},\, l_{2}).
 \]
where events $A(x_{1},\r_{1},l_{1})$, $A(x_{1}+l_{1}+d,\,\r_{2},\,l_{2})$
belong to disjoint intervals.
Recall the definition of  $\p(\r,l)$ in~\eqref{eq:p(,)}.
By the Markov property,
 \[
   \Pbof{B(d,l)\mid X(x_{1})=k} \le \sum_{l_{1}+d+l_{2}=l} 
          \p(r_{1}, l_{1}) \p(\r_{2}, l_{2}).
 \]
Hence Condition~\ref{cond:distr}.\ref{i:distr.wall-ub} implies
$\sum_{l} \Pbof{B(d,l)\mid X(x_{1})=k} \le \sum_{l_{1}} \p(\r_{1}, l_{1})
          \sum_{l_{2}} \p(\r_{2}, l_{2}) \le \p(\r_{1})\p(\r_{2})$,
which completes the proof.
 \end{proof}

In the lemma below, we use $w_{1},w_{2}$: please note that these are integer
coordinates,
and have nothing to do with the trap probability upper bound $\tub$: we will
never have these two uses of $w$ in a place where they can be confused.

\begin{lemma}\label{lem:compound-hole-lb}
Let $u \le v_{1} < w_{2}$, and $a$ be given with $v_{1}-u\le 12\bub^{*}$, and let
 \begin{equation*}
                 b = a + \cei{(v_{1}-u)/2},
  \qquad         c = b + (v_{1}-u) + 1.
 \end{equation*}
Assume that $Y=y$ is fixed in such a way that $W$ is a compound
horizontal wall with body $\rint{v_{1}}{w_{2}}$, and type $\ang{\r_{1}, \r_{2}, i}$,
with rank $\r$ as given in~\eqref{eq:compound-rank}.
Assume also that the component walls $W_{1}, W_{2}$ already satisfy the hole
lower bound, Condition~\ref{cond:distr}.\ref{i:distr.hole-lb}.
Let 
 \[
  E_{2} = E_{2}(u, v_{1}, w_{2};\, a) = E^{*}(u, v_{1}, w_{2};\, a)
 \]
where $E^{*}$ was defined in~\eqref{eq:E-star}.
Assume
 \begin{equation}\label{eq:s-ub}
   (\bub^{*})^{\hxp} \h(\r_{j}) \le 0.07,\txt{ for } j=1,2.
 \end{equation}
Let $k\in\{1,\dots,\m\}$.
Then
 \begin{equation}\label{eq:compound-hole-lb}
  \Pbof{E_{2}\mid X(a)=k,\; Y=y} \ge 
   (c-b)^{\hxp} (\lg^{i}/2)^{\hxp} \h(\r_{1}) \h(\r_{2}) \cdot(1 - V)
 \end{equation}
with
 $V = 2\cdot (26 \pub \bub^{*} + 338 \tubh(\bub^{*})^{2}) / \h(\r_{1} \lor \r_{2})$. 
 \end{lemma}
 \begin{Proof}
Figure~\ref{fig:compound-hole-lb} shows the role of the various coordinates in
the proof.
\begin{figure}
  \centering
  \includegraphics[scale=0.5]{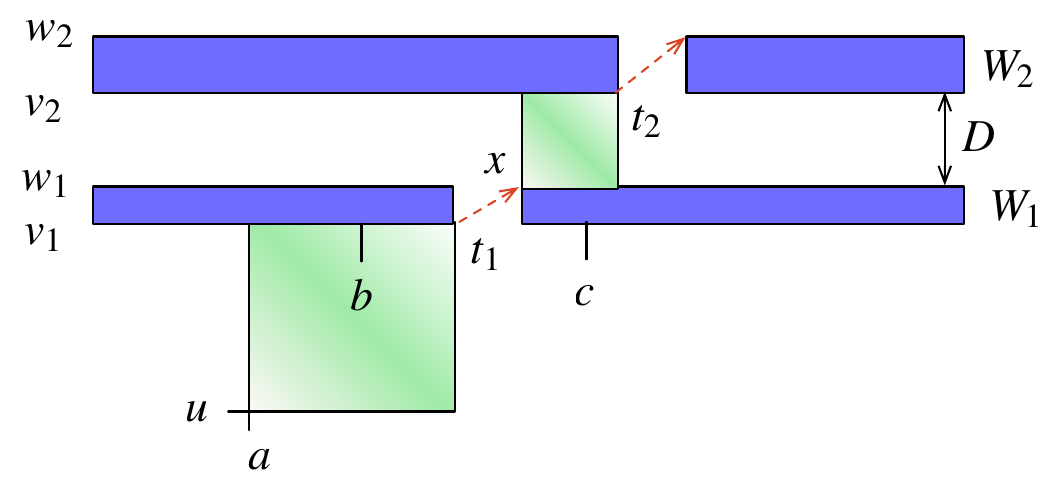}
  \caption{Hole through a compound wall.  The light rectangles are inner H-clean.}
\label{fig:compound-hole-lb}
\end{figure}
Let $\D$ be the distance between the component walls
$W_{1}, W_{2}$ of the wall $W$, where the body of $W_{i}$ is $\rint{v_{i}}{w_{i}}$.
Consider first passing through $W_{1}$.
For each $x \in [b, c + \bub - 1]$, 
let $A_{x}$ be the event that $E^{*}(u, v_{1}, w_{1};\, a)$ holds
with the vertical projection of the hole ending at $x$,
and that $x$ is the smallest possible number with this property.
Let $B_{x} = E^{*}(w_{1}, v_{2}, w_{2};\, x)$.

 \begin{step+}{step:compound-hole-lb.Ax-Bx}
We have $E_{2} \spsq \bigcup_{x} (A_{x} \cap B_{x})$.
 \end{step+}
 \begin{pproof}
If for some $x$ we have $A_{x}$, then there is a rectangle
$\Rect(\pair{a}{u},\pair{t_{1}}{v_{1}})$ satisfying the requirements 
of $E^{*}(u, v_{1}, w_{1};\, a)$ and also a 
hole $\Rect(\pair{t_{1}}{v_{1}},\pair{x}{w_{1}})$ through the first wall.
If also $B_{x}$ holds, then there is a rectangle
$\Rect(\pair{x}{w_{1}},\pair{t_{2}}{v_{2}})$ satisfying the requirements 
of $E^{*}(w_{1}, v_{2}, w_{2};\, x)$,
and also a hole $\Rect(\pair{t_{2}}{v_{2}},\pair{x'}{w_{2}})$ through the second
wall.

Let us show $\pair{t_{1}}{v_{1}}\leadsto\pair{x'}{w_{2}}$.
Since $|x'-t_{1}|\le|w_{2}-v_{1}|$, this will imply that
the interval $\rint{t_{1}}{x'}$ is a hole that
passes through the compound wall $W$.

We already know $\pair{t_{1}}{v_{1}}\leadsto\pair{x}{w_{1}}$
and $\pair{t_{2}}{v_{2}}\leadsto(x', w_{2})$; we still need to prove
$\pair{x}{w_{1}}\leadsto\pair{t_{2}}{v_{2}}$.

The requirements imply that $\Rect(\pair{x}{w_{1}},\pair{t_{2}}{v_{2}})$ is a
hop of $\cM$.
Indeed, the inner H-cleanness of $\rint{x}{t_{2}}$ in the process $X$
follows from $B_{x}$.
The inner cleanness of $\rint{w_{1}}{v_{2}}$ in the process $Y$ is implied by the
fact that $\rint{v_{1}}{w_{2}}$ is a compound wall.
The fact that $W$ is a compound wall also implies that the interval
$\rint{w_{1}}{v_{2}}$ contains no horizontal walls.

According to $B_{x}$, this rectangle has the necessary slope constraints, 
hence by 
the reachability condition of $\cM$, its endpoint is reachable from
its starting point.
 \end{pproof}

It remains to lower-bound $\Prob\bigparen{\bigcup_{x} (A_{x} \cap B_{x})}$.
For each $x$, the events $A_{x},B_{x}$ belong to disjoint intervals, 
and the events $A_{x}$ are disjoint of each other.
 
 \begin{step+}{step:compound-hole-lb.Ax}
Let us lower-bound $\sum_{x} \Prob(A_{x})$.
 \end{step+}
 \begin{prooof}
We have, using the notation of Lemma~\ref{lem:hole-lb-2}:
 $\sum_{x} \Prob(A_{x}) = \Prob(E^{*}(u, v_{1}, w_{1};\, a))$. 
Lemma~\ref{lem:hole-lb-2} is applicable and we get 
$\Prob(E^{*}(u, v_{1}, w_{1};\, a)) \ge F_{1} - U$ with
 \begin{equation}\label{eq:compound-hole-lb.U}
 \begin{split}
   F_{1} &= 0.5 \land ((c-b)^{\hxp} \h(\r_{1})),
\\     U &= 26 \pub \bub^{*} + 338 \tubh(\bub^{*})^{2}.
 \end{split}
 \end{equation}
By the assumption~\eqref{eq:s-ub}:
$(c-b)^{\hxp} \h(\r_{1}) \le (7\bub^{*})^{\hxp} \h(\r_{1}) \le 0.5$,
hence the operation $0.5 \land$ can be deleted from $F_{1}$:
 \begin{equation}\label{eq:compound-hole-lb.F1}
   F_{1} = G_{1} := (c-b)^{\hxp} \h(\r_{1}).
 \end{equation}
 \end{prooof}

 \begin{step+}{step:compound-hole-lb.Bx}
Let us now lower-bound $\Prob(B_{x})$, for an arbitrary 
condition $X(x) = k$ for $k\in\{1,\dots,\m\}$.
 \end{step+}
 \begin{prooof}
We have $B_{x} = E^{*}(w_{1}, v_{2}, w_{2};\, x)$.
The conditions of Lemma~\ref{lem:hole-lb-2} are satisfied for
$u = w_{1}$, $v = v_{2}$, $w = w_{2}$, $a = x$.
It follows that $\Prob(B_{x}) \ge F_{2} - U$ with
 \begin{equation*}
   F_{2} = 0.5 \land ((\flo{\D/2}+1)^{\hxp} \h(\r_{2})),
 \end{equation*}
which can again be simplified using assumption~\eqref{eq:s-ub} and
$\D\le\f\le\bub^{*}$:
 \[
  F_{2} = G_{2} := (\flo{\D/2}+1)^{\hxp} \h(\r_{2}).
 \]
 \end{prooof}

 \begin{step+}{step:compound-hole-lb.combine}
Let us combine these estimates, using 
$G = G_{1} \land G_{2} > \h(\r_{1}\lor \r_{2})$.
 \end{step+}
 \begin{prooof}
By the Markov property, we find that the lower bound on 
$\Prob(B_{x})$ (for arbitrary $X(x)=k$) is also a lower bound on
$\Prob(B_{x} \mid A_{x})$:
 \begin{align*}
   \Prob(E_{2}) &\ge \sum_{x} \Prob(A_{x})\Prob(B_{x} \mid A_{x})
      \ge (G_{1} - U)(G_{2} - U)
\\    &\ge G_{1} G_{2}(1 - U(1/G_{1} + 1/G_{2}))
       \ge G_{1} G_{2}(1 - 2 U/G)
\\    &= (c-b)^{\hxp} (\flo{\D/2}+1)^{\hxp} \h(\r_{1}) \h(\r_{2}) 
       (1 - 2 U / G)
\\    &\ge (c-b)^{\hxp} (\flo{\D/2}+1)^{\hxp} \h(\r_{1}) \h(\r_{2}) 
       (1 - 2 U/h(\r_{1}\lor \r_{2})).
 \end{align*}
 \end{prooof}
 \begin{step+}{step:compound-hole-lb.distance}
We conclude by showing $\flo{\D/2}+1 \ge \lg^{i}/2$.
 \end{step+}
 \begin{prooof}
Recall $d_{i} \le D < d_{i+1}$ where $d_{i}$ was defined
in~\eqref{eq:jump-types}.
For $i=0,1$, we have $\flo{D/2}+1 = 1 > \lg^{1}/2$.
For $i\ge 2$, we have $\flo{D/2}+1 \ge D/2 \ge \lg^{i}/2$.
 \end{prooof}
 \end{Proof}

\section{The scale-up functions}\label{sec:params}

Lemma~\ref{lem:main} says
that there is an $\tilde\m$ such that if $\m>\tilde\m$ then the
sequence $\cM^{k}$ can be constructed in such a way that the 
claim~\eqref{eq:main} of the main lemma holds.
If we computed some $\tilde\m$ explicitly then all parameters of the 
construction could be
turned into constants: but this is unrewarding work and it would only make
the relationships between the parameters less intelligible.
We prefer to name all these parameters, to point out the necessary
inequalities among them, and finally to show that if $\m$ is sufficiently
large then all these inequalities can be satisfied simultaneously.

Mazery $\cM^{1}$ is defined in Example~\ref{xmp:base}, with
a parameter $\tub=\tub_{1}$ that can be anyting not less than $\frac{1}{\m-1}$;
let $\tub_{1}=\frac{1}{\m-1}$.
Instead of $\m$, it will be more convenient to use the parameter
$\R_{0}=\R_{0}(\m)$ introduced below which defines $\tub_{1}=\lg^{-\R_{0}\tubxp\txp}$
via Definitions~\ref{def:ranks} and~\ref{def:exponential} via constants $\lg,\tubxp,\txp$:
 \begin{align}\label{eq:m-expr}
 \R_{0}=\frac{2\log_{\lg}(\m-1)}{\tubxp\txp},\quad 
 \m=1+1/\tub_{1}\le 2/\tub_{1}=2\lg^{\R_{0}\tubxp\txp}.
 \end{align}
In what follows, rather than asking $\m$ to be sufficiently large, we will ask,
equivalently, $\R_{0}$ to be sufficiently large.

The following definition introduces some of the parameters needed for scale-up.
Recall that the slope lower bound $\slb$ must satisfy $\slb < 1/2$.

 \begin{definition}\label{def:ranks}
We set
 \begin{equation}\label{eq:slb1}
  \slb_{1} = 0.
 \end{equation}
To obtain the new rank lower bound, we multiply $\R$ by a constant:
 \begin{equation}\label{eq:txp}
           \R = \R_{k}  = \R_{0} \txp^{k},
\quad       \R_{k+1} = \R^{*} = \R \txp,
\quad                1 < \txp < 2,\quad 1<\R_{0}.
 \end{equation}
The rank of emerging walls, introduced in~\eqref{eq:emerg-rank}, is defined
using a new parameter $\txp'$:
 \[
   \hat\R = \txp'\R.
 \]
 \end{definition}

We require
 \begin{equation}\label{eq:txp'}
 \txp < \txp' < \txp^{2}.
 \end{equation}

We need some bounds on the possible rank values.

 \begin{definition} Let $\txpub = 2\txp / (\txp - 1)$.
 \end{definition}

  \begin{lemma}[Rank upper bound]\label{lem:rank-bds}
In a mazery, all ranks are upper-bounded by $\txpub \R$.
  \end{lemma}
 \begin{proof}
The proof is by induction on $k$.
The statement is true for $k=1$, where not being any barriers, certainly all
their ranks are bounded by $\txpub\R_{1}$.
Assume the statement for $k$, we will prove it for $k+1$.
Since $\txp'<\txp^{2}<\txpub$, the rank upper bound in $\cM$ is larger than the
rank of emerging walls.
New walls in $\cM^{*}$ are either emerging walls, or are obtained by
applying the compounding operation, possibly twice, to a wall of $\cM$ or to an
emerging wall, that is to a wall of rank $\le\txpub R$.
Each compounding operation adds the rank of a light wall, less than $\R^{*}$.
The ranks in $\cM^{*}$  are thus less than $\txpub R+2\R^{*}=(\txpub/\tau+2)\R^{*}=\txpub\R^{*}$.
 \end{proof}

 \begin{corollary}\label{c.rank-lifetime}
Every rank exists in $\cM^{k}$ for at most
$\cei{\log_{\txp}\frac{2 \txp}{\txp - 1}}$ values of $k$.
 \end{corollary}
 \begin{proof} Immediate.    \end{proof}

Recall $\lg = 2^{1/2}$, as defined in~\eqref{eq:lg-def}.
It can be seen from the definition of compound ranks
in~\eqref{eq:compound-rank} and from
Lemma~\ref{lem:compound-wall-ub} that the probability bound $\p(\r)$ 
of a wall should be approximately $\lg^{-\r}$.
The actual definition makes the bound a little smaller:

It is convenient to express 
several other parameters of $\cM$ and the scale-up in terms of a single
one, $\T$:

 \begin{definition}[Exponential relations]\label{def:exponential}
 Let $\T = \lg^{\R}$,
 \begin{align*}
      \bub     &= \T^{\bubxp},
\quad     \f        = \T^{\fxp},
\quad     \g        = \T^{\gxp},
\quad     \tub      = \T^{-\tubxp}.
 \end{align*}
As we will see, $\tubxp$ just needs to be sufficiently large with respect to the
other constants.
On the other hand, we require
 \begin{equation}\label{eq:bubxp-etc}
           0 < \bubxp < \gxp < \fxp < 1.
 \end{equation}
The values $\bubxp,\fxp,\gxp,\txp,\txp'$
will be chosen independent of the mazery level.  
A bound on $\fxp$ has been indicated in the requirement~\eqref{eq:R-cond}
which will be satisfied by
 \begin{equation}\label{eq:txp-ub}
   \txp \le 2 - \fxp.
 \end{equation}
We turn this into equality: $\txp = 2 - \fxp$.
 \end{definition}

Let us estimate $\bub^{*}$.
Emerging walls can have size as large as $4\g'+ 4\bub$, % was $3\g'$, for RS referee
and at the time of their creation, they
are the largest existing ones.
We get the largest new walls when the compound operation combines these
with light walls on both sides, leaving the largest gap possible,
so the largest new wall size is
 \[
  4 \g' + 2\f + 6\bub < 3\f,
 \]
where we used~\eqref{eq:bub-g-f}.
Hence any value larger than $3 \f$ can be chosen as 
$\bub^{*} = \bub^{\txp}$. 
With $\R_{0}$ large enough, we always get this if
 \begin{equation}\label{eq:bubxp-lb}
  \fxp < \txp \bubxp.
 \end{equation}
As a reformulation of one of the inequalities of~\eqref{eq:bub-g-f},
we require
 \begin{equation}\label{eq:g-smaller}
 \gxp \ge \frac{\bubxp + \fxp}{2}.
 \end{equation}
We also need
 \begin{align}\label{eq:correl-trap-xp}
                               4(\gxp + \bubxp) &< \tubxp(4 - \txp),
\\\label{eq:emerg-xp.2}         4\gxp + 6\bubxp + \txp'  &< \tubxp,
\\\label{eq:emerg-xp.3}         \txp(\bubxp+1)   &< \txp'.
 \end{align}
Using the exponent $\hxp$ introduced in~\eqref{eq:hxp}, we require
 \begin{align}
\label{eq:hxp-ub.1}                  \txp \hxp &< \gxp - \bubxp,
\\\label{eq:hxp-ub.2}              \txpub \hxp &< 1 - \txp \bubxp,
\\\label{eq:hxp-ub.3}              \txpub \hxp &< \tubxp - 2 \txp \bubxp.
 \end{align}
Note that all these inequalities require $\hxp$ to be just sufficiently small.

 \begin{lemma}\label{lem:exponents-choice}
 The exponents $\bubxp, \fxp, \gxp, \txp, \txp', \hxp$ can be chosen to
satisfy the inequalities~\eqref{eq:txp},\eqref{eq:txp'},
\eqref{eq:bubxp-etc}-\eqref{eq:hxp-ub.3}.
 \end{lemma}
 \begin{proof}
It can be checked that the choices
$\bubxp = 0.15$, $\gxp = 0.2$, $\fxp = 0.25$, $\txp=1.75$, $\txp' = 2.5$,
$\tubxp=4.5$, $\txpub=4.66\dots$ satisfy all the inequalities in question.
 \end{proof}

 \begin{definition}
Let us fix now 
the exponents $\bubxp, \fxp, \gxp, \txp, \txp', \hxp$  as chosen in the lemma.
In order to satisfy all our requirements also for small $k$,
we will fix $\aux_{2}$ sufficiently small,
then $\aux_{1}$ sufficiently large, then $\aux_{3}$ sufficiently
large, and finally $\R_{0}$ sufficiently large.
 \end{definition}

We need to specify some additional parameters.

 \begin{definition}\label{def:tub-def} We define
 \begin{alignat}{4}\label{eq:pub-def}
            &\pub &&= \T^{-1},&\quad\tubh &&= \tub,
\\\label{eq:new-ncln}
           &\ncln^{*} &&= \ncln + \bub^{*}\pub.
 \end{alignat}
 \end{definition}

\section{Probability bounds after scale-up}\label{sec:after-scale-up}

\subsection{Bounds on traps}

The structures $\cM^{k}$ are now defined but we have not proved that
they are mazeries, since not all inequalities required in the definition of
mazeries have been verified yet.

 \begin{lemma}\label{lem:trap-scale-up}
For any value of the constant $\aux_{3}$, if $\R_{0}$ is sufficiently large then
the following holds: 
if $\cM=\cM^{k}$ is a mazery then $\cM^{*}$ satisfies the trap upper 
bound~\ref{cond:distr}.\ref{i:distr.trap-ub}.
 \end{lemma}
 \begin{proof}
For some string $x = \tup{x(0), x(1), \dots}$, for a point $\pair{a}{b}$, 
let $\cE$ be the event that a trap starts at $\pair{a}{b}$.
We assume $Y(b-1)=\s$ fixed arbitrarily.
We need to bound $\Pbof{\cE \mid X= x,\;Y(b-1)=\s}$.
There are three kinds of trap in $\cM^{*}$: uncorrelated and correlated
compound traps, and traps of the missing-hole kind.
Let $\cE_{1}$ be the event that an uncorrelated trap occurs.
According to~\eqref{eq:uncorrel-trap}, using $\txp = 2 - \fxp$ (and
recalling that $\tub^{*}$ plays the role of the
trap probability upper bound $\tub$  for $\cM^{*}$):
 \begin{equation*}
  \begin{split}
   \Pbof{\cE_{1} \mid X= x,\;Y(b-1)=\s} &\le 2 \f^{2} \tub^{2} 
    = 2 \T^{2\fxp -2\tubxp} 
\\ &= 2 \T^{-\txp\tubxp - (2 - \txp)\tubxp + 2\fxp}
    = \tub^{*}\cdot 2/\f^{\tubxp - 2}.
  \end{split}
 \end{equation*}
This can be made smaller than $\tub^{*}$ by an
arbitrarily large factor if $\R_{0}$ is large.

Let $\cE_{2}$ be the event that a vertical correlated trap appears.
By Lemma~\ref{lem:correl}, using~\eqref{eq:g-prime}:
  \begin{multline*}
   \Pbof{\cE_{2} \mid X = x,\;Y(b-1)=\s} \le \sum_{j=1}^{2} (5\bub\l_{j}\tub)^{4}
    \le 2\cdot (5\bub\g'\tub)^{4}
\\ = 2\cdot 11^{4}\T^{4\gxp + 4\bubxp - 4\tubxp -\txp\tubxp + \txp\tubxp}
  = 2\tub^{*}\cdot 11^{4}\T^{4(\gxp + \bubxp) - \tubxp(4 - \txp)}.
  \end{multline*}
Due to~\eqref{eq:correl-trap-xp}, this can be made smaller than $\tub^{*}$
by an arbitrarily large factor if $\R_{0}$ is large.

Let $\cE_{3}$ be the event that a vertical trap of the missing-hole kind appears
at $\pair{a}{b}$.
Lemma~\ref{lem:missing-hole} implies for $n = \Flo{\frac{\g}{3 \bub}}$:
 \begin{equation*}
 \Pbof{\cE_{3} \mid X=x,\;Y(b-1)=\s} \le e^{-(1-\ncln)^{2} n \h(\R^{*})}.
 \end{equation*}
Further, using inequality~\eqref{eq:bub-g-f} and the largeness of $\R_{0}$:
 \begin{align*}
  n > \g / (3 \bub) - 1 > \g/(4\bub) = \T^{\gxp-\bubxp}/4. 
% changed -2 to -1, for RS referee.
 \end{align*}
Now,
 \begin{align*}
   \h(\R^{*})  &= \aux_{3}\T^{-\txp\hxp},
\\  (1-\ncln)^{2} n \h(\R^{*}) &>  0.8 n \h(\R^{*})
         > 0.2 \aux_{3} \T^{\gxp - \bubxp - \txp\hxp},
\\ \Pbof{\cE_{3} \mid X=x,\;Y(b-1)=\s} &\le  e^{-0.2 \aux_{3}\T^{\gxp - \bubxp - \txp\hxp}}.
 \end{align*}
Due to~\eqref{eq:hxp-ub.1}, this can be made smaller than $\tub^{*}$ by an
arbitrarily large factor if $\R_{0}$ is large.

For $j=1,2$, let $\cE_{4,j}$ be the event that a
horizontal trap of the correlated kind of type $j$ starts at $\pair{a}{b}$.
Let  $\cE_{4,3}$ be the event that a horizontal % horizontal added for RS referee
trap of missing-hole kind starts at $\pair{a}{b}$.
Lemmas~\ref{lem:correl} and~\ref{lem:missing-hole} imply
 \[
  \Pbof{\cE_{4,j} \mid X=x,\;Y(b-1)=\s} \le \tub^{2} = \tub^{*}\T^{-\tubxp(2-\txp)}.
 \]
Due to~\eqref{eq:txp-ub}, this can be made smaller than $\tub^{*}$ by an
arbitrarily large factor if $\R_{0}$ is large.

Thus, if $\R_{0}$ is sufficiently large then the sum of these six
probabilities is still less than $\tub^{*}$.
 \end{proof}

\subsection{Bounds on walls}

Recall the definition of $\p(\r)$ in~\eqref{eq:wall-prob}.

 \begin{lemma}\label{lem:emerg-contrib}
For every possible value of $\aux_{1},\aux_{2},\aux_{3}$, 
if  $\R_{0}$ is sufficiently large then the following holds.
Assume that $\cM = \cM^{k}$ is a mazery.
Fixing any point $a$ and fixing $X(a)$ in any way, the
sum of the probabilities over $l$ that a barrier of the emerging kind of size
$l$ starts at $a$ is at most $\p(\hat\R) / 2 = \p(\txp'\R) / 2$.
 \end{lemma}
 \begin{proof}
We use the result and notation of Lemma~\ref{lem:emerg}, and also
the estimate of $\Prob(\cE_{3})$ in the proof of Lemma~\ref{lem:trap-scale-up},
replacing $\m$ with the upper bound from~\eqref{eq:m-expr}
 \begin{align*}
 \sum_{l}\Prob(\cF(l)) &\le 8\tub_{1}^{-1}\bub^{2}\tub^{2}\bigparen{2\cdot(5\bub\g')^{4} 
     + \tub^{-4} e^{- (1-\ncln)^{2} n \h(\R^{*})}},
\\ 8\tub_{1}^{-1}\bub^{2}\tub^{-2} e^{- (1-\ncln)^{2} n \h(\R^{*})} &\le 
   8\tub_{1}^{-1}\T^{2\tubxp+2\bubxp}e^{- 0.2\aux_{3}\T^{\gxp - \bubxp - \txp\hxp}}.
 \end{align*}
Due to~\eqref{eq:hxp-ub.1}, the last expression decreases exponentially in
$\T$, so for sufficiently large $\R_{0}$ it is less than $\p(\txp'\R)$ by an
arbitrarily large factor.
On the other hand, using~\eqref{eq:g-prime}:
 \begin{align*}
      8\tub_{1}^{-1}\bub^{2}\tub^{2}\cdot 2\cdot (5\bub\g')^{4} 
   &\le 16\cdot 11^{4}\T^{-\tubxp + 4\gxp + 6\bubxp}
 = 16\T^{-\txp'}\cdot 11^{4}\T^{4\gxp + 6\bubxp + \txp' - \tubxp}.
 \end{align*}
If $\R_{0}$ is sufficiently large then, due to~\eqref{eq:emerg-xp.2}, this
is less than $\p(\txp'\R)$ by an arbitrarily large factor.
\end{proof}

 \begin{lemma}\label{lem:compound-contrib}
For a given value of $\aux_{2}$, if 
we choose the constants $\aux_{1}, \R_{0}$ sufficiently large
in this order then the following holds.
Assume that $\cM=\cM^{k}$ is a mazery.
After one operation of forming compound barriers,
fixing any point $a$ and fixing $X(a)$ in any way,
for any rank $\r$, the sum, over all widths $l$, of the probability that
a compound barrier of rank $\r$ and width $l$ starts at $a$ is at most 
$\p(\r)\R^{-\aux_{1}/2}$.
 \end{lemma}
 \begin{proof}
Let $\r_{1} \le \r_{2}$ be two ranks,
and assume that $\r_{1}$ is light: $\r_{1} < \R^{*}=\txp\R$.
With these, we can form compound barriers of type
$\ang{\r_{1},\r_{2},i}$.
The bound~\eqref{eq:compound-wall-ub} and the definition of $\p(r)$ 
in~\eqref{eq:wall-prob} shows
that the contribution by this term to the sum of probabilities, over all
widths $l$, that a barrier of rank $\r = \r_{1}+\r_{2}-i$ and size $l$
starts at $x$ is at most
 \[
 \lg^{i}\p(\r_{1})\p(\r_{2}) 
  = \aux_{2}^{2}\lg^{-\r}(\r_{1}\r_{2})^{-\aux_{1}} = \aux_{2}(\r/\r_{1}\r_{2})^{\aux_{1}}\p(\r).
 \]
Now
$\r_{1}\r_{2} \ge \R\r_{2} \ge (\R/2)(\r_{1}+\r_{2}) \ge \r\R/2$,
hence the above bound reduces to $\aux_{2}(\R/2)^{-\aux_{1}}\p(\r)$.
The same rank $\r$ can be obtained by the compound operation at most
the following number of times:
 \[
 |\setof{(i,\r_{1}): i \le \R\fxp,\; \r_{1} < \txp\R}| \le (\fxp\R+1)\txp\R.
 \]
The total probability contributed to rank $\r$ is therefore at most
 \[
  \begin{split}
    & \aux_{2}(\R/2)^{-\aux_{1}}\p(\r)(\fxp\R+1)\txp\R
    <  \p(\r)\R^{-\aux_{1}/2}
  \end{split}
 \]
if $\R_{0}$ and $\aux_{1}$ are sufficiently large.
 \end{proof}

 \begin{lemma}\label{lem:all-wall-ub}
For every choice of $\aux_{1},\aux_{2},\aux_{3}$ if we choose $\R_{0}$
sufficiently large then the following holds.
Suppose that each structure $\cM^{i}$ for $i \le k$ is a mazery.
Then Condition~\ref{cond:distr}.\ref{i:distr.wall-ub} holds for $\cM^{k+1}$.
\end{lemma}
 \begin{proof}
By Corollary~\ref{c.rank-lifetime}, 
each rank $\r$ occurs for at most a constant
number $n = \cei{\log_{\txp}\frac{2 \txp}{\txp - 1}}$ values of $i \le k$.
For any rank, a barrier can be formed only as an emerging
barrier or compound barrier.
The first can happen for one $i$ only, and Lemma~\ref{lem:emerg-contrib}
bounds the probability contribution by $\p(\r)/2$.
A compound barrier can be formed for at most $n$ values of $i$, and for
each value in at most two steps.
Lemma~\ref{lem:compound-contrib} bounds each contribution
by $\p(\r)\R^{-\aux_{1}/2}$.
After these increases, the probability becomes
at most $\p(\r)(1/2 + 2n \R^{-\aux_{1}/2}) < \p(\r)$
if $\R_{0}$ is sufficiently large.
 \end{proof}

\subsection{Auxiliary bounds}

The next two lemmas show that the choices made in
Definition~\ref{def:tub-def} satisfy the requirements
imposed in Definition~\eqref{eq:pub-intr}.
Recall the introduction of the wall probability upper bound $\pub$
in Definition~\ref{def:pub-intr} and its value assignment $\T^{-1}$
in~\eqref{eq:pub-def}.

 \begin{lemma}\label{lem:pub}
For small enough $\aux_{2}$, the probability of a barrier
of $\cM$ starting at a given point $b$ is bounded by $\pub$.
 \end{lemma}
  \begin{proof}
We have 
 $\sum_{\r \ge \R} \p(\r) < \aux_{2} \sum_{\r \ge \R} \lg^{-\r}
  = \lg^{-\R}\aux_{2}(1 - 1/\lg)^{-1} < \lg^{-\R}$
if $\aux_{2} < 1 - 1/\lg$.
  \end{proof}

 \begin{lemma}\label{lem:0.6}
If $\R_{0}$ is sufficiently large then
 $\sum_{k} (2\bub_{k+1}\pub_{k} + \bub_{k+1}^{2}\tub_{k}) < 0.5$.
 \end{lemma}

 \begin{proof} Substituting the definitions of $\pub$ in~\eqref{eq:pub-def},
further the values of all other parameters given in Definitions~\ref{def:ranks} 
and~\ref{def:exponential}:
 \begin{align*}
 \sum_{k} (2\bub_{k+1}\pub_{k} + \bub_{k+1}^{2}\tub_{k})
 \le 2 \sum_{k} \lg^{-\R_{0}\txp^{k}(1 - \bubxp\txp)}
  + \sum_{k} \lg^{-\R_{0}\txp^{k}(\tubxp - 2\bubxp\txp)}   
 \end{align*}
 which because of~\eqref{eq:hxp-ub.2} and~\eqref{eq:hxp-ub.3}, 
is less than 0.5 if $\R_{0}$ is large.
 \end{proof}

Note that for $\R_{0}$ large enough, the relations
 \begin{align}
  \label{eq:strong-ncln-ub}
                     \bub^{*}\pub &< 0.5(0.1 - \ncln),   % 0.05
\\\label{eq:strong-f-intro}
     \slopeincr \g / \f &< 0.5(0.5 - \slb)
 \end{align}
hold for $\cM=\cM^{1}$ and $\slb = \slb_{1}$.
This is clear for~\eqref{eq:strong-ncln-ub}. 
For~\eqref{eq:strong-f-intro}, since $\slb_{1}=0$ according
to~\eqref{eq:slb1}, we only need
$0.25 > \slopeincr \g / \f = \slopeincr \T^{-(\fxp - \gxp)}$,
which is satisfied if $\R_{0}$ is large enough.

\begin{lemma}\label{lem:ncln-ub}
Suppose that the structure $\cM=\cM^{k}$ is a mazery and it
satisfies~\eqref{eq:strong-ncln-ub} and~\eqref{eq:strong-f-intro}.
Then $\cM^{*}=\cM^{k+1}$ also satisfies these inequalities if $\R_{0}$ is
chosen sufficiently large (independently of $k$), % RS referee
and also satisfies Condition~\ref{cond:distr}.\ref{i:distr.ncln}.
 \end{lemma}

 \begin{proof} % expanded and simplified for RS referee.
Let us show first that $\cM^{*}$ also satisfies the inequalities if $\R_{0}$ is
chosen sufficiently large.

For sufficiently large $\R_{0}$, we have $\bub^{**} \pub^{*} < 0.5 \bub^{*} \pub$.
Indeed, this says $\T^{(\txp\bubxp - 1)(\txp - 1)} < 0.5$.
Hence using~\eqref{eq:strong-ncln-ub} 
and the definition of $\ncln^{*}$ in~\eqref{eq:new-ncln}:
\begin{align*}
   \bub^{**} \pub^{*} &\le 0.5\bub^{*} \pub 
    \le 0.5(0.1 - \ncln) - 0.5\bub^{*}\pub
\\           &\le 0.5(0.1 - \ncln) - 0.5(\ncln^{*}-\ncln) = 0.5(0.1 - \ncln^{*}).
  \end{align*}
This is inequality~\eqref{eq:strong-ncln-ub} for $\cM^{*}$.

For inequality~\eqref{eq:strong-f-intro}, the scale-up definition Definition~\ref{def:new-slb} says 
$\slb^{*} - \slb = \slopeincr \g / \f$.
The inequality $\g^{*} / \f^{*} < 0.5\g/\f$ is guaranteed if $\R_{0}$ is large.
From here, we can conclude the proof as for $\ncln$.

To verify Condition~\ref{cond:distr}.\ref{i:distr.ncln} for
$\cM^{*}$, recall the definition in~\eqref{eq:new-ncln} of
 \begin{align*}
  \ncln^{*}= \ncln + \bub^{*}\pub.
 \end{align*}
For the first inequality of Condition~\ref{cond:distr}.\ref{i:distr.ncln},
we upper-bound the conditional probability that a point $a$ of the line is 
strongly clean in $\cM$ but not in $\cM^{*}$ by 
 \begin{equation}\label{eq:ncln-traps}
  (2\f/3+\bub)\pub,
 \end{equation}
which upper-bounds the probability that a 
horizontal barrier of $\cM$ starts in $\rint{a - \f/3-\bub}{a + \f/3}$.
This can be upper-bounded by $\f\pub < \bub^{*}\pub/3$ by~\eqref{eq:bub-star}.
Hence an upper bound on the conditional probability of not strong cleanness in
$\cM^{*}$ is $\ncln/2+\bub^{*} \pub/3<\ncln^{*}/2$ as required.

For the other inequalities in Condition~\ref{cond:distr}.\ref{i:distr.ncln},
consider a rectangle $Q=\Rect^{\to}(u,v)$ and fix $X(u_{0})=k$, and $Y=y$.
The conditional probability that a point $u$ is not trap-clean in $Q$ for $\cM$ but
not for $\cM^{*}$ is upper-bounded by the probability of the appearance of a
trap of $\cM$ within a distance $\g$ of point $u$ in $Q$.
There are at most $\g^{2}$ positions for the trap, so a bound is
\[
    \g^{2}\tub =  \T^{2\gxp - \tubxp} < 0.5\T^{\txp \bubxp - 1}.
\]
For the latter inequality, for large $\R_{0}$, we need to check
$2\gxp-\txp\bubxp+1<\tubxp$.
But~\eqref{eq:emerg-xp.2} and~$\txp'>1$ imply even
the stronger $4\gxp+1<\tubxp$.
We conclude the same way for the first inequality.
The argument for the other inequalities in
Condition~\ref{cond:distr}.\ref{i:distr.ncln} is identical.
 \end{proof}

\subsection{Lower bounds on holes}

We will make use of the following estimate.

 \begin{lemma}\label{lem:hop-lb}
Let $\rint{a_{0}}{b_{0}}$, $\clint{a_{1}}{b_{1}}$ be intervals with length 
$\le 12\bub^{*}$. 
Suppose that the sequence $Y$ and the value $\s\in\{1,\dots,\m\}$
are fixed arbitrarily. 
If $\aux_{3}$ and then $\R_{0}$ are chosen sufficiently large then
the following event holds with probability at least $0.75$ even if
conditioned on $X(a_{0})=\s$:
The rectangle $Q = \rint{a_{0}}{b_{0}}\times\clint{a_{1}}{b_{1}}$ is inner
H-clean for $\cM^{*}$, 
and contains no traps or vertical barriers of $\cM$ or $\cM^{*}$.
 \end{lemma}
 \begin{proof}
We will just write ``probability'' but will understand conditional probability,
when $Y=y$ and $X(a_{0})$ is fixed.
According to Lemma~\ref{lem:ncln-ub},
the probability that one of the two cleanness conditions is not
satisfied is at most $0.2$.
Using Lemmas~\ref{lem:trap-scale-up}, \ref{lem:all-wall-ub} and~\ref{lem:pub},
the probability that a vertical barrier of $\cM$ or $\cM^{*}$ is contained in
$Q$ is at most
 \begin{align*}
 12\bub^{*}(\pub+\pub^{*})\le 24\bub^{*}\pub=24\T^{\txp\bubxp - 1}.
 \end{align*}
The probability that a trap of $\cM$ or $\cM^{*}$ 
is contained in $Q$ is at most 
 \begin{align*}
 12\bub^{*}(12\bub^{*}+1)(\tub+\tub^{*})<2\cdot 156(\bub^{*})^{2}\tub
=312\T^{2 \txp\bubxp - \tubxp}.
 \end{align*}
If $\R_{0}$ is sufficiently large, then
the sum of the last two terms is at most 0.05.
 \end{proof}

\begin{lemma}\label{lem:emerg-hole-lb}
For emerging walls, the fitting holes
satisfy Condition~\ref{cond:distr}.\ref{i:distr.hole-lb}
if $\R_{0}$ is sufficiently large.
\end{lemma}
 \begin{Proof}
   \begin{figure}
     \centering
     \includegraphics[scale=0.5]{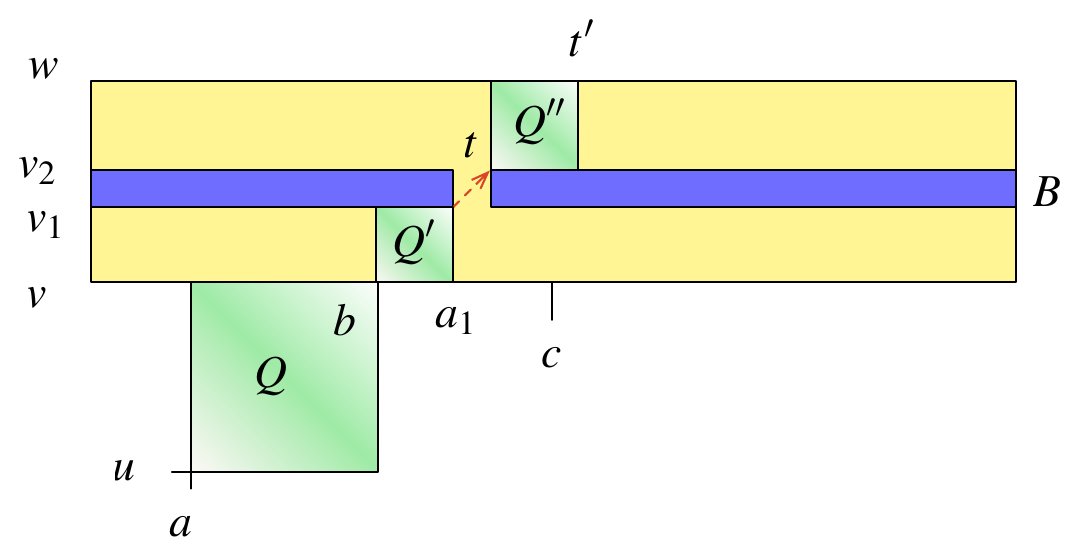}
     \caption{Hole through an emerging wall: 
to the proof of Lemma~\ref{lem:emerg-hole-lb}}
     \label{fig:emerg-hole-lb}
   \end{figure}
Recall Condition~\ref{cond:distr}.\ref{i:distr.hole-lb}
applied to the present case.
Let $u \le v < w$, $a$ be given with $v-u \le 12\bub^{*}$, and define
$b = a + \cei{(v-u)/2}$, $c = a + (v-u) + 1$.
Assume that $Y=y$ is fixed in such a way that
$B$ is a horizontal wall of the emerging kind
with body $\rint{v}{w}$.
Let $E^{*} = E^{*}(u, v, w;\, a)$ be as given in Definition~\ref{def:E-star}.
Let $\s\in\{1,\dots,\m\}$.
We will prove
 \begin{equation*}
   \Pbof{E^{*} \mid X(a)=\s,\;Y = y} \ge (c-b)^{\hxp} \h(\hat\R).
 \end{equation*}
Recall the definition of emerging walls in Step~\ref{defstep:emerg}
of the scale-up construction.
The condition at the end says, in our case, that
either $\rint{v}{w}$ is a hop of $Y$ or it can be partitioned
into a light (horizontal) wall $\rint{v_{1}}{v_{2}}$ of some rank $\r$, 
and two (possibly empty) hops surrounding it: so, $v \le v_{1} < v_{2} \le w$.
Without loss of generality, assume this latter possibility.
Let
 \[
      a_{1} = b + (v_{1} - v).
 \]
Let $\cF$ be the event that 
 \begin{enumerate}[\upshape (a)]

  \item\label{i:emerg.hops} 
Rectangle $Q = \Rect^{\rightarrow}(\pair{a}{u}, \pair{b}{v})$ contains no vertical
barriers or traps of $\cM$ or $\cM^{*}$, further is inner H-clean 
for $\cM^{*}$.
Rectangle $Q'=\Rect^{\rightarrow}(\pair{b}{v}, \pair{a_{1}}{v_{1}})$ contains
no vertical barriers or traps of $\cM$, and is inner H-clean for $\cM$.

  \item\label{i:middle-hole}
For an arbitrary $t\in \rint{a_{1}}{a_{1}+\bub}$, let
 \[
   t' = t+(w-v_{2}).
 \]
We require that event $E(v_{1},v_{1},v_{2};\,a_{1})$
is realized with some hole $\rint{a_{1}}{t}$, and the rectangle
$Q''=\Rect^{\rightarrow}(\pair{t}{v_{2}}, \pair{t'}{w})$ contains no
vertical barriers or traps, and is inner H-clean for $\cM$.

 \end{enumerate}
 \begin{step+}{step:emerg-hole-lb.implies}
Event $\cF$ implies the event $E^{*}$ of Definition~\ref{def:E-star}, with $d$
taken to be equal to $b$.
 \end{step+}
 \begin{pproof}
Assume that $\cF$ holds.
Rectangle $\Rect^{\rightarrow}(\pair{a}{u},\pair{b}{v})$ has the necessary inner
cleanness properties, the absence of barriers and traps, and the slope lower
bound to have $\pair{a}{u}\leadsto\pair{b}{v}$:
it remains to show $\pair{b}{v}\leadsto\pair{t'}{w}$.
We have $\pair{b}{v}\leadsto \pair{a_{1}}{v_{1}}$
since $\Rect^{\rightarrow}(\pair{b}{v}, \pair{a_{1}}{v_{1}})$ is a hop that is
also a square, hence satisfies the slope lower bound condition.
For similar reasons, $\pair{t}{v_{2}}\leadsto \pair{t'}{w}$.
Also, since $\rint{a_{1}}{t}$ is a hole through $\rint{v_{1}}{v_{2}}$, we
have $\pair{a_{1}}{v_{1}}\leadsto \pair{t}{v_{2}}$.
 \end{pproof}

 \begin{step+}{step:emerg-hole-lb.estim}
We have $\Pbof{\cF\mid X(a)=\s,\;Y = y} \ge 0.75^{3}\aux_{3}\T^{-\hxp\txp}$.
 \end{step+}
 \begin{pproof}
Condition~\eqref{i:emerg.hops} in the definition of $\cF$ 
is coming from two rectangles with disjoint
projections, therefore by the method used throughout the paper, we can
multiply their probability lower bounds, which are given as 0.75
by Lemma~\ref{lem:hop-lb}.

Condition~\eqref{i:middle-hole}
also refers to an event with a projection to the $x$ axis disjoint from the previous ones.
The probability of the existence of a hole is lower-bounded via
Condition~\ref{cond:distr}.\ref{i:distr.hole-lb}, by
$\h(\r) \ge \h(\R^{*})=\aux_{3}\T^{-\hxp\txp}$.
A reasoning similar to the proof of Lemma~\ref{lem:hole-lb-clean-2} shows
that the whole condition~\eqref{i:middle-hole} is satisfied at the expense
of another factor 0.75 via Lemma~\ref{lem:hop-lb}.
 \end{pproof}
Recall the definition of $\h(\r)$ in~\eqref{eq:h-def}.
The required lower bound of Condition~\ref{cond:distr}.\ref{i:distr.hole-lb} is 
 \begin{align*}
   (c-b)^{\hxp} \h(\hat\R) &\le (6 \bub^{*}+1)^{\hxp} \h(\hat\R)
   \le (7\T^{\txp \bubxp})^{\hxp} \h(\txp'\R)
   = \aux_{3}7^{\hxp}\T^{\hxp(\txp\bubxp-\txp')} 
\\ &< 0.75^{3}\aux_{3}\T^{-\hxp\txp}
 \end{align*}
if $\R_{0}$ is sufficiently large, due to~\eqref{eq:emerg-xp.3}.
  \end{Proof}

 \begin{lemma}\label{lem:all-compound-hole-lb}
After choosing $\aux_{1},\aux_{3},\R_{0}$ sufficiently large in this order, 
the following holds.
Assume that $\cM=\cM^{k}$ is a mazery: then every
compound wall satisfies the hole lower bound,
Condition~\ref{cond:distr}.\ref{i:distr.hole-lb}, provided its components
satisfy it.
 \end{lemma}

 \begin{Proof}
The proof starts from the setup in
Lemma~\ref{lem:compound-hole-lb}, making the appropriate substitutions in
the estimates.
 \begin{step+}{step:all-compound-hole-lb.recall}
Recall what is required.
 \end{step+}
 \begin{prooof}
Let $u \le v_{1} < w_{2}$, $a$ be given with $v_{1}-u \le 12\bub^{*}$, and define
 \begin{equation*}
                 b = a + \cei{(v_{1}-u)/2},
  \qquad         c = a + (v_{1}-u) + 1.
 \end{equation*}
Assume that $Y = y$ is fixed in such a way that there is a
compound horizontal wall $W$ with body $\rint{v_{1}}{w_{2}}$, and type
$\ang{\r_{1}, \r_{2}, i}$, with rank
\begin{align}
  \label{eq:r-def}
  \r = \r_{1} + \r_{2} - i
\end{align}
 as in~\eqref{eq:compound-rank}.
Also, let $X(a)=\s$ be fixed in an arbitrary way.
Let $E_{2} = E^{*}(u, v_{1}, w_{2};\, a)$ as defined
in~\eqref{eq:E-star}.
We need to prove $\Pbof{E_{2}\mid X(a)=\s} \ge (c-b)^{\hxp}\h(\r)$.
 \end{prooof}
 
\begin{step+}{step:all-compound-hole-lb.apply}
 Let us apply Lemma~\ref{lem:compound-hole-lb}.
\end{step+}
 \begin{prooof}
Recall the definition of $\h(\r)$ in~\eqref{eq:h-def}.
The assumption $(\bub^{*})^{\hxp} \h(\r_{i}) \le 0.07$ of the lemma holds since
\begin{align}
  \label{eq:h-ub}
      \h(\r_{i}) = \aux_{3} \lg^{-\hxp\r_{i}} \le \aux_{3} \T^{-\hxp},
\quad (\bub^{*})^{\hxp} \h(\r_{i}) \le \aux_{3}\T^{-\hxp(1 - \bubxp\txp)}
\end{align}
which, due to~\eqref{eq:hxp-ub.2}, is always smaller than $0.07$ if $\R_{0}$
is sufficiently large.
We conclude
 \begin{equation*}
  \Pbof{E_{2}\mid X(a)=\s} \ge 
   (c-b)^{\hxp} (\lg^{i}/2)^{\hxp} \h(\r_{1}) \h(\r_{2}) \cdot (1-V)
 \end{equation*}
with  $V = 2\cdot (26 \pub \bub^{*} + 338 \tubh(\bub^{*})^{2})
      / \h(\r_{1} \lor \r_{2})$. 
 \end{prooof}

 \begin{step+}{step:all-compound-hole-lb.mainpart}
Let us estimate the part of this expression before $1-V$.
 \end{step+}
 \begin{prooof}
Using and the formula for $\h(\r_{i})$ in~\eqref{eq:h-ub} and the definition
of $\r$ in~\eqref{eq:r-def}:
 \begin{align*}
   (\lg^{i}/2)^{\hxp}\h(\r_{1})\h(\r_{2}) &=2^{-\hxp}\aux_{3}^{2} \lg^{-\hxp(\r_{1}+\r_{2}-i)} 
     = 2^{-\hxp}\aux_{3}^{2}\lg^{-\r\hxp},
\\ (c-b)^{\hxp} (\lg^{i}/2)^{\hxp} \h(\r_{1}) \h(\r_{2})
     &\ge  2^{-\hxp}\aux_{3} (c-b)^{\hxp}\h(\r)
       > 2 (c-b)^{\hxp}\h(\r)
 \end{align*}
if $\aux_{3}$ is sufficiently large.
 \end{prooof}

 \begin{step+}{step:all-compound-hole-lb.estim}
To complete the proof, we show $1 - V \ge 0.5$ for large enough $\R_{0}$.
 \end{step+}
 \begin{prooof}
Lemma~\ref{lem:rank-bds} gives $\r_{1}\lor\r_{2} \le \txpub \R$, hence
 \[
  \h(\r_{1}\lor\r_{2}) = \aux_{3}\lg^{-(\r_{1}\lor\r_{2})\hxp} 
              \ge \aux_{3}\lg^{- \hxp\txpub\R} = \aux_{3} \T^{-\hxp\txpub}.
 \]
 Let us estimate both parts of $V$:
 \begin{equation*}
  \pub\bub^{*}/\h(\r_{1}\lor\r_{2}) 
       \le \aux_{3}^{-1}\T^{\hxp \txpub + \txp \bubxp - 1}, 
\quad \tubh(\bub^{*})^{2}/\h(\r_{1}\lor \r_{2}) 
       \le \aux_{3}^{-1}\T^{\hxp \txpub + 2 \txp \bubxp - \tubxp}.
  \end{equation*}
Conditions~\eqref{eq:hxp-ub.2}-\eqref{eq:hxp-ub.3} imply that $V$
can be made arbitrarily small if $\R_{0}$ is sufficiently large.
 \end{prooof}
 \end{Proof} % l.all-compound-hole-lb.

For the hole lower bound condition for $\cM^{*}$, there is one more case to consider.

 \begin{lemma}\label{lem:heavy-hole-lb}
After choosing $\aux_{1},\aux_{3},\R_{0}$ sufficiently large in this order, 
the following holds.
Assume that $\cM=\cM^{k}$ is a mazery: then every
wall of $\cM^{k+1}$ that is also a heavy wall of $\cM^{k}$
satisfies the hole lower bound,
Condition~\ref{cond:distr}.\ref{i:distr.hole-lb}.
 \end{lemma}
 \begin{proof}
Recall Condition~\ref{cond:distr}.\ref{i:distr.hole-lb}
applied to the present case.
Let $u \le v < w$, $a$ be given with $v-u \le 12\bub^{*}$, and define
 $b = a + \cei{(v-u)/2}$, $c = a + (v-u) + 1$.
Assume that $Y=y$ is fixed in such a way that
$B$ is a horizontal wall of $\cM$ with body $\rint{v}{w}$, with rank 
$\r\ge\R^{*}$ (since it is a heavy wall).  % RS referee
Assume also that $X(a)=\s$ is fixed arbitrarily.
Let $E^{*} = E^{*}(u, v, w;\, a)$ be defined
as after~\eqref{eq:E-star}.
We will prove
 \begin{equation*}
   \Pbof{E^{*} \mid X(a)=\s,\;Y = y} \ge (c-b)^{\hxp} \h(\r).
 \end{equation*}
Suppose first $v-u\le 12\bub$.
Then the fact that $\cM^{k}$ is a mazery implies the same inequality with $E$ in
place of $E^{*}$.
In our case, however, the event $E$ implies $E^{*}$, as shown already at the
beginning of the proof of Lemma~\ref{lem:hole-lb-2}.

It remains to check the case $v-u > 12\bub$:
for this, Lemma~\ref{lem:hole-lb-2} says 
 \begin{equation*}
   \Prob(E^{*}) \ge 0.5 \land (1.1(c-b)^{\hxp} \h(\r)) - U
 \end{equation*}
with  $U = 26 \pub \bub^{*} + 338 \tubh(\bub^{*})^{2}$.
The operation $0.5\land$ can be 
omitted since $1.1(c-b)^{\hxp} \h(\r)\le 0.5$.
Indeed, $c-b\le 7\bub^{*}$ implies
 \begin{equation*}
  1.1(c-b)^{\hxp} \h(\r) \le 7.7 \aux_{3}\lg^{\R\txp\bubxp\hxp}\lg^{-\r\hxp}
   = 7.7\aux_{3}\lg^{\hxp(\R\txp\bubxp-\r)}.
 \end{equation*}
It follows from~\eqref{eq:hxp-ub.2} that $\txp\bubxp<1$.
Since $\r\ge\R$, the right-hand side can be made $<0.5$ for
large enough $\R_{0}$.
Now
 \begin{equation*}
   1.1(c-b)^{\hxp} \h(\r) - U \ge (c-b)^{\hxp}\h(\r)(1.1-U/\h(\r)).
 \end{equation*}
The part subtracted from 1.1 is less than 0.1 if $\R_{0}$
is sufficiently large, by the same argument as the estimate of $V$ at the
end of the proof of Lemma~\ref{lem:all-compound-hole-lb}.
 \end{proof}

 \section{The approximation lemma}
\label{sec:approx-proof}

The crucial combinatorial step in proving the main lemma is the following.

 \begin{lemma}[Approximation]\label{lem:approx}
The reachability condition, Condition~\ref{cond:reachable}, holds for $\cM^{*}$ if
$\R_{0}$ is sufficiently large.
 \end{lemma}

The name suggest to view our renormalization method as
successive approximations: the lemma shows reachability in the absence of
some less likely events (traps, walls, and uncleanness in the corners of the
rectangle).
The present section is taken up by the proof of this lemma.
Recall that we are considering a bottom-open or left-open or closed rectangle $Q$
with starting point $u$ and endpoint $v$ with 
 \[
   \minslope(u, v) \ge \slb^{*} = \slb + \slopeincr \g / \f.
 \]
Denote $u = \pair{u_{0}}{u_{1}}$, $v = \pair{v_{0}}{v_{1}}$.
We require $Q$ to be a hop of $\cM^{*}$.
Thus, the points $u, v$ are clean for $\cM^{*}$ in $Q$, and $Q$
contains no traps or walls of $\cM^{*}$.
We have to show $u\leadsto v$.
Without loss of generality, assume 
 \[
  Q = I_{0} \times I_{1} = \Rect^{\eps}(u, v)  
 \]
with $|I_{1}| \le |I_{0}|$, where $\eps=\rightarrow,\uparrow$ or nothing.

 \subsection{Walls and trap covers}

Let us determine the properties of the set of walls in $Q$.

 \begin{lemma}\label{lem:the-grate}
Under conditions of Lemma~\ref{lem:approx}, with the notation given in the
discussion after the lemma, the following holds.
 \begin{alphenum}

  \item\label{i:grate.pre-hop}
For $\d=0,1$, for some $n_{\d} \ge 0$, there
is a sequence $W_{d,1},\dots,W_{d,n_{d}}$ of dominant 
light neighbor walls of $\cM$ separated from each other 
by external hops of $\cM$ of size $>\f$, and from the ends of $I_{\d}$
(if $n_{\d}>0$) by hops of $\cM$ of size  $\ge \f/3$.

  \item\label{i:grate.hole}
For every (horizontal) wall $W$ of $\cM$ occurring in $I_{1}$,
for every subinterval $J$ of $I_{0}$ of size $\g$ such that $J$ is at a
distance $\ge\g+7\bub$ from the ends of $I_{0}$, there is an outer 
rightward-clean hole fitting $W$, its endpoints at a distance of at
least $\bub$ from the endpoints of $J$.
The same holds if we interchange vertical and horizontal.
 \end{alphenum}
 \end{lemma}

 \begin{proof}
This is a direct consequence of Lemmas~\ref{lem:new-hop}
and~\ref{lem:if-not-missing-hole}.
 \end{proof}

From now on, in this proof, whenever we mention a \df{wall} we mean one of
the walls $W_{\d, i}$, and whenever we mention a trap then, unless said
otherwise, we mean only traps of $\cM$ not intersecting any of these walls.
Let us limit the places where traps can appear in $Q$.

 \begin{definition}[Trap cover]
A set of the form $I_{0} \times J$
with $|J| \le 4 \bub$ containing the starting point of a
trap of $\cM$ will be called a \df{horizontal trap cover}.
Vertical trap covers are defined similarly.
 \end{definition}

In the following lemma, when
we talk about the distance between two traps, we mean the distance between
their starting points.

 \begin{lemma}[Trap cover]\label{lem:trap-cover}
Let $T_{1}$ be a trap of $\cM$ contained in $Q$.
Then there is a horizontal or vertical trap cover $U \spsq T_{1}$ such that
the starting point of
every other trap in $Q$ is either contained in $U$ or is
at least at a distance $\f - \bub$ from $T_{1}$.
If the trap cover is vertical, it intersects none of the vertical walls
$W_{0, i}$; if it is horizontal, it intersects none of the horizontal
walls $W_{1,j}$.
 \end{lemma}
 \begin{proof}
Let $\pair{a_{1}}{b_{1}}$ be the starting point of $T_{1}$.
If there is no trap $T_{2} \sbsq Q$, with starting point 
$\pair{a_{2}}{b_{2}}$, closer than $\f - \bub$ to $T_{1}$,
such that $|a_{2} - a_{1}| \ge 2 \bub$, then 
 $U = \clint{a_{1} - 2 \bub}{a_{1} + 2 \bub} \times I_{1}$ will do.
Otherwise, let $T_{2}$ be such a trap and let
$U = I_{0} \times \clint{b_{1} - 2\bub}{b_{1} + 2\bub}$.
We have $|b_{2} - b_{1}| < \bub$, since otherwise $T_{1}$ and $T_{2}$
would form together an uncorrelated compound trap, which was excluded.

Consider now a trap $T_{3} \sbsq Q$, with starting point $\pair{a_{3}}{b_{3}}$,
at a distance $< \f - \bub$ from $\pair{a_{1}}{b_{1}}$.
We will show $\pair{a_{3}}{b_{3}} \in U$.
Suppose it is not so:
then $|a_{3} - a_{1}| < \bub$, otherwise $T_{1}$ and $T_{3}$ would
form an uncorrelated compound trap.
Also, the distance of $\pair{a_{2}}{b_{2}}$ and $\pair{a_{3}}{b_{3}}$ must be at least
$\f$, since otherwise they would form an uncorrelated compound trap.
Since $|a_{2} - a_{1}| < \f - \bub$ and $|a_{3} - a_{1}| < \bub$, we have
$|a_{2} - a_{3}| < \f$.
Therefore necessarily $|b_{2} - b_{3}| \ge \f$.
Since $|b_{2} - b_{1}| < \bub$, it follows $|b_{3} - b_{1}| > \f - \bub$,
so $T_{3}$ is at a distance at least $\f - \bub$ from $T_{1}$, contrary to
our assumption.

If the trap cover thus constructed is vertical and intersects some vertical
wall, just decrease it so that it does not intersect any such walls.
Similarly with horizontal trap covers.
 \end{proof}

Let us measure distances from the diagonal.

\begin{definition}[Relations to the diagonal]\label{def:diagonal}
Define, for a point $a = \pair{a_{0}}{a_{1}}$:
 \[
   d(a) = (a_{1}-u_{1}) - \slope(u,v)(a_{0} - u_{0})
 \]
to be the distance of $a$ above the diagonal of $Q$, then for
$w=\pair{x}{y}$, $w'=\pair{x'}{y'}$:
 \begin{equation}\label{eq:d-diff}
   \begin{aligned}
            d(w')-d(w) &= y'-y - \slope(u,v)(x'-x),
\\            |d(w')-d(w)| &\le |y'-y|+|x'-x|,
 \end{aligned}
 \end{equation}
as $\slope(u,v)\le 1$.
We define the strip
 \[
 C^{\eps}(u, v, h_{1}, h_{2})
= \setof{w \in \Rect^{\eps}(u,v) :
 h_{1} < d(w) \le h_{2}},
 \]
a channel of vertical width $h_{2}-h_{1}$ parallel to the diagonal of
$\Rect^{\eps}(u,v)$. 
\end{definition}
 \begin{figure}
   \centering
   \includegraphics[scale=0.8]{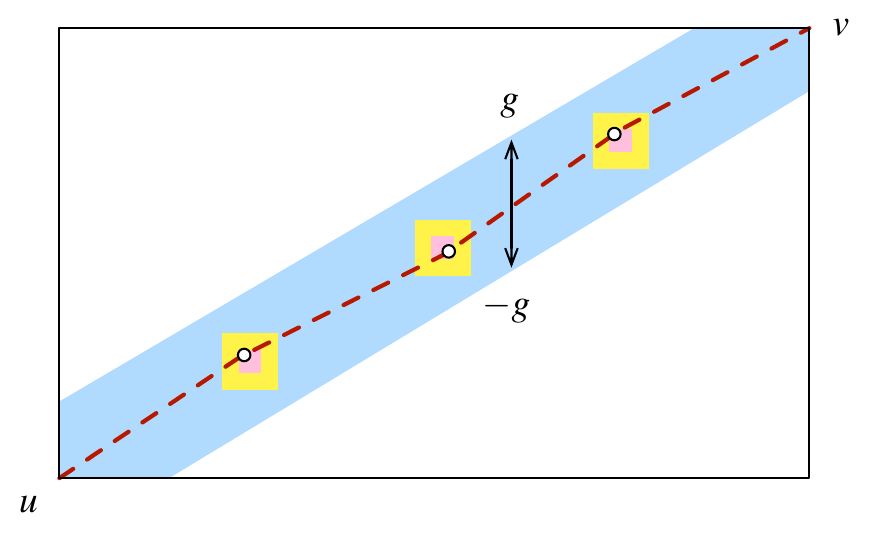}
   \caption{To Lemma~\ref{lem:approx.no-walls}.  
The shaded area is $ C^{\eps}(u, v, -\g, \g)$.}
   \label{fig:channel}
 \end{figure}

 \begin{lemma}\label{lem:approx.no-walls}
Assume that points $u, v$ are clean for $\cM$ in $Q=\Rect^{\eps}(u,v)$,
with 
 \[
   1\ge \slope(u, v) \ge \slb + 4 \g / \f.
 \]
If $C = C^{\eps}(u, v, -\g, \g)$ contains 
no traps or walls of $\cM$ then $u \leadsto v$.
 \end{lemma}
 \begin{proof}
If $|I_{0}| < \g$ then there is no trap or wall in $Q$, therefore $Q$ is a
hop and we are done via Condition~\ref{cond:reachable} for $\cM$.
Suppose $|I_{0}| \ge \g$.
Let 
 \[
      n = \Cei{\frac{|I_{0}|}{0.9 \g}},
\quad h = \frac{|I_{0}|}{n}.
 \]
Then $\g / 2 \le h \le 0.9 \g$.
Indeed, the second inequality is immediate.
For the first one, if $n\le 2$, we have $\g \le |I_{0}| = n h \le 2 h$, and
for $n\ge 3$:
 \begin{align*}
   \frac{|I_{0}|}{0.9\g}&\ge n-1,
\\ |I_{0}|/n &\ge (1-1/n) 0.9\g \ge 0.6\g.   
 \end{align*}
For $i = 1, 2, \dots, n-1$, let
 \[
      a_{i} = u_{0} + i h,
\quad b_{i} = u_{1} + i h\cdot \slope(u, v),
\quad w_{i} = \pair{a_{i}}{b_{i}},
\quad S_{i} = w_{i} + \clint{-\bub}{2\bub}^{2}.
 \]
Let us show $S_{i}\sbsq C$.
For all elements $w$ of $S_{i}$, we have $|d(w)| \le 3\bub$, and we know
$3\bub < \g$ from~\eqref{eq:bub-g-f}.
To see $S_{i} \sbsq \Rect^{\eps}(u,v)$,
we need (from the worst case $i=n-1$) 
$\slope(u,v) h > 2\bub$.
Using~\eqref{eq:bub-g-f} and the assumptions of the lemma:
 \[
   \frac{2\bub}{h} \le \frac{2\bub}{\g/2}
   = 4\bub/\g \le 4\g/\f \le \slope(u,v).
 \]
By Remark~\ref{rem:distr}.\ref{i:middle-third},
there is a clean point $w'_{i} = \pair{a'_{i}}{b'_{i}}$ in the middle third 
$w_{i}+\clint{0}{\bub}^{2}$ of $S_{i}$.
Let $w'_{0} = u$, $w'_{n} = v$.
By their definition, each rectangle $\Rect^{\eps}(w'_{i}, w'_{i+1})$
has size $< 0.9 \g + \bub < \g-\bub$, hence
falls into the channel $C$ and is consequently trap-free.

Let us show $\minslope(w'_{i}, w'_{i+1}) \ge \slb$: this will imply
$w'_{i}\leadsto w'_{i+1}$.
It is sufficient to show $\slope(w'_{i}, w'_{i+1}) \ge \slb$.
Note that,
using $s=\slope(u,v)$, the distance from $w'_{i}$ to $w'_{i+1}$ is at most $h+\bub$ in the $x$
coordinate and at least $s h -\bub$ in the $y$ coordinate. 
% Let $x \ge h-\bub > \g/3$ 
% be the horizontal projection of $w'_{i+1} - w'_{i}$.
We have
  \begin{align*}
 \slope(w'_{i},w'_{i+1}) \ge \frac{s h - \bub}{h+\bub} 
   = s - \frac{\bub + s\bub}{h+\bub} \ge s - \frac{2\bub}{\g/2} 
   \ge s - 4\g/\f\ge \slb,
  \end{align*}
where the last inequality used again~\eqref{eq:bub-g-f}.
 \end{proof}

We introduce particular strips around the diagonal.

 \begin{definition}
Let
\begin{align}
  \label{eq:H-def}  H  &= 12,
\\\label{eq:C-def}  C  &=  C^{\eps}(u, v, -3H\g, 3H\g).
 \end{align}
 \end{definition}

Inequalities~\eqref{eq:bub-g-f} imply
 \begin{equation}\label{eq:H-slopeincr}
  \slopeincr \ge 33H + 7.
 \end{equation}

Let us introduce the system of walls and trap covers we will have to
overcome.

 \begin{definition}
Let us define a sequence of trap covers $U_{1}, U_{2},\dots$ as follows.
If some trap $T_{1}$ is in $C$, then let
$U_{1}$ be a (horizontal or vertical) trap cover covering it according to
Lemma~\ref{lem:trap-cover}.
If $U_{i}$ has been defined already and there is a trap $T_{i+1}$ in
$C$ not covered by $\bigcup_{j \le i} U_{j}$ then let $U_{i+1}$
be a trap cover covering this new trap.
To each trap cover $U_{i}$ we assign a real number $a_{i}$ as follows.
Let $\pair{a_{i}}{a'_{i}}$ be the intersection of the diagonal of $Q$ and the
left or bottom edge of $U_{i}$ (if $U_{i}$ is vertical or horizontal
respectively).
Let $\pair{b_{i}}{b'_{i}}$ be the intersection of the diagonal and the left edge
of the vertical wall $W_{0,i}$ introduced in Lemma~\ref{lem:the-grate},
and let $\pair{c'_{i}}{c_{i}}$ be the intersection of the diagonal and the
bottom edge of the horizontal wall $W_{1,i}$.
Let us define the finite set
 \begin{equation*}%\label{eq:union-of-three}
   \set{s_{1}, s_{2},\dots} 
= \set{a_{1}, a_{2}, \dots} \cup \set{b_{1}, b_{2}, \dots} \cup 
  \set{c'_{1}, c'_{2}, \dots}
  \end{equation*}
where $s_{i} \le s_{i+1}$.

We will call the objects (trap covers or walls) belonging to the points
$s_{i}$ our \df{obstacles}.
 \end{definition}

 \begin{lemma}\label{lem:one-of-three}
If $s_{i},s_{j}$ belong to the same obstacle category among the three
(horizontal wall, vertical wall, trap cover), then
$|s_{i}-s_{j}| \ge 3\f/4$.
 \end{lemma}
It follows that for every $i$ at least one of the three numbers
$(s_{i+1} - s_{i})$, $(s_{i+2} - s_{i+1})$, $(s_{i+3} - s_{i+2})$ is larger
than $\f / 4$.
 \begin{proof}
If both $s_{i}$ and $s_{j}$ belong to walls of the same orientation
then they are farther than $\f$ from each other, 
since the walls from which they come are at least $\f$ apart.
(For the numbers $c'_{i}$, this uses $\slope(u,v) \le 1$.)

Suppose that both belong to the set $\set{a_{1}, a_{2}, \dots}$, say 
they are $a_{1} \le a_{2}$, coming from $U_{1}$ and $U_{2}$.
Let $\pair{x_{j}}{y_{j}}$ be the starting point of some trap $T_{j}$ in
$U_{j} \cap C$ (with $C$ defined in~\eqref{eq:C-def}).
If $U_{j}$ is vertical then $|x_{j} - a_{j}| \le 4 \bub$,
and $|y_{j} - a'_{j}| \le 3H \g + 4\bub$.
If $U_{j}$ is horizontal then 
$|x_{j} - a_{j}| \le (3H \g + 4\bub) / \slope(u,v)$,
and $|y_{j} - a'_{j}| \le 4 \bub$.

Suppose that $a_{2} - a_{1} \le 0.75 \f$, then also 
$a'_{2} - a'_{1} \le 0.75 \f$.
From the above estimates it follows that
 \begin{align*}
  |x_{2}-x_{1}| \lor |y_{2} - y_{1}| 
    &\le 0.75\f + (2 \cdot 3 H \g+8\bub)/ \slope(u,v) 
     \le 0.75\f + 2.1\cdot 3 H \f / \slopeincr
\\  &= \f - 0.05 \f - (0.2 - 2.1\cdot 3 H / \slopeincr)\f
     \le \f - 0.05\f < \f - \bub,
 \end{align*}
where we used $\slope(u, v) \ge \slb^{*} \ge \slopeincr \g / \f$,
\eqref{eq:H-slopeincr}
and $\bub < 0.05 \f$ which follows from~\eqref{eq:bub-g-f}.
But this would mean that the starting points of the traps $T_{j}$
are closer than $\f - \bub$, in contradiction to Lemma~\ref{lem:trap-cover}. 
 \end{proof}

 \subsection{Passing through the obstacles}

The remark after Lemma~\ref{lem:one-of-three} allows us to break up the sequence
of obstacles into groups of size at most three, which can be dealt with separately.
This leads to 
the following, weaker, form of the Approximation Lemma, which still
takes a lot of sweat (case distinctions) to prove:

 \begin{lemma}\label{lem:three}
Assume $\slope(u,v) \le 1$, $(\slb^{*}=)\slb + \slopeincr \g / \f < 1/2$, and let
$u, v$ be points with
 \begin{equation}\label{eq:triangle-cnd.three}
   \slb + (\slopeincr-1) \g / \f \le \slope(u, v).
 \end{equation}
Assume that the set $\set{s_{1}, s_{2}, \dots}$ defined above 
consists of at most three elements, with
the consecutive elements less than $\f/4$ apart.
Assume also 
 \begin{equation}\label{eq:far-from-ends}
   v_{0} - s_{i},\; s_{i} - u_{0} \ge 0.1 \f.
 \end{equation}
Then if $\Rect^{\rightarrow}(u, v)$ or $\Rect^{\uparrow}(u, v)$ is a 
hop of $\cM^{*}$ then $u\leadsto v$.
 \end{lemma}
\begin{Proof}
%% Let us note that 
%% inequality~\eqref{eq:triangle-cnd.three} implies, using~\eqref{eq:bub-g-f}:
%% $\slope(u,v) \ge (\slopeincr-1)\g/\f \ge (\slopeincr-1)\bub/\g$, and hence
%%  \begin{equation*}%\label{eq:slope-u-b-lb}
%%   \bub/\slope(u,v) \le \g/(\slopeincr-1).
%%  \end{equation*}
We can assume without loss of generality that there are indeed three
points $s_{1}, s_{2}, s_{3}$.
By Lemma~\ref{lem:one-of-three}, they must then come from three obstacles of
different categories: $\set{s_{1}, s_{2}, s_{3}} = \set{a, b, c'}$ where
$b$ comes from a vertical wall, $c'$ from a horizontal wall, 
and $a$ from a trap cover.
There is a number of cases: we illustrated the most complex one in
Figure~\ref{fig:close-trap}.

If the index $i\in\{1,2,3\}$ 
of a trap cover is adjacent to the index of a wall of the same
orientation, then this pair will be called a \df{parallel pair}.
A parallel pair is either horizontal or vertical.
It will be called a \df{trap-wall pair} if
the trap cover comes first, and the \df{wall-trap pair} if the wall comes
first.
If $s_{i}-s_{i-1} < 1.1\g$ 
for a vertical pair or $(s_{i}-s_{i-1})\slope(u,v) < 1.1\g$  for a horizontal
pair then we say that the pair is \df{bound}.
Thus, a pair is bound if the distance between the starting edges of its
obstacles is less than $1.1\g$.
A bound pair (if exists)
is more difficult to pass, therefore its crossing points
will be chosen in a coordinated way, starting from the trap cover side.

We will call an obstacle $i$ \df{free}, if it is not part of a bound pair.
Consider the three disjoint channels 
 \[
  C(u, v, K, K + 2H\g), \txt{ for } K=-3H\g,\, -H\g,\, H\g.
 \]
The three lines (bottom or left edges) of the trap covers or walls
corresponding to $s_{1},s_{2},s_{3}$ can intersect in at most two places,
so at least one of the above channels does not contain such an
intersection.
Let $K$ belong to such a channel.
For $i \in\{1,2,3\}$, we shall choose points 
 \[
         w_{i} = \pair{x_{i}}{y_{i}},
\quad   w'_{i} = \pair{x'_{i}}{y'_{i}},
\quad  w''_{i} = \pair{x''_{i}}{y''_{i}}
 \]
in the channel $C(u, v, K+2\g, K+(2H-2)\g)$ in such a way that $w_{i}$ is on
the (horizontal or vertical) line corresponding to $s_{i}$.
Not all these points will be defined.
The points $w_{i},w'_{i},w''_{i}$ will always be defined if $i$ is free.
Their role in this case is the following:
$w'_{i}$ and $w''_{i}$ are points on the two sides of the trap cover or
wall with $w'_{i}\leadsto w''_{i}$.
Point $w_{i}$ will be on the starting edge of the obstacle, and it
will direct us in locating $w'_{i},w''_{i}$.
However, $w_{i}$ by itself will not determine $w'_{i},w''_{i}$:
other factors are involved.
Correspondingly, for each free obstacle, we will distinguish
a \df{forward} way of crossing (when $d(w_{i})$ will be made equal to
$d(w''_{j})$ for some $j<i$)
and a \df{backward} way of crossing (when $d(w_{i})$ will be made equal to
$d(w'_{j})$ for some $j>i$).

Part~\ref{step:three.x-free} of the proof
collects some estimates on crossing free obstacles.
Part~\ref{step:three.free} proves the lemma in the
cases when all obstacles are free.
In the case when there is a horizontal trap-wall pair,
finding the wall crossing point requires more
freedom than finding the trap cover crossing point.
To be able do it first, all obstacles are crossed in the
backward direction.

Part~\ref{step:three.x-bound} collects estimates for
different cases of crossing of a bound pair.
Finally, part~\ref{step:three.bound} proves the lemma also for the cases
where a bound pair exists.

 \begin{step+}{step:three.x-free}
Consider crossing a free obstacle $s_{i}$,
assuming that $w_{i}$ has been defined already.
 \end{step+}
 \begin{prooof} % x-free
There are cases corresponding to whether the obstacle is a trap
cover or a wall, and whether it is vertical or horizontal.
Backward crossings are quite similar to forward ones.

  \begin{step+}{step:three.x-free.trapc} 
 Consider crossing a trap cover $s_{i}$.
  \end{step+}
  \begin{prooof} % trapc.vert
  \begin{step+}{step:three.x-free.trapc.vert}
 Assume that the trap cover is vertical.
  \end{step+}
 \begin{prooof}
 Consider crossing a vertical trap cover forwards.
Recall
 \begin{align*}
   \l_{1}=7\bub,\; L_{1}=4\l_{1}.
 \end{align*}
Let us apply Lemma~\ref{lem:if-correl} to vertical correlated traps
$J\times I'$, with 
$J=\clint{x_{i}}{x_{i}+5\bub}$, $I'=\clint{y_{i}}{y_{i}+L_{1}}$. 
The lemma is applicable since $w_{i}\in C(u, v, K+2\g, K+(2H-2)\g)$ implies
$u_{1}< y_{i}-L_{1}-\l_{1} < y_{i}+2L_{1}+\l_{1} < v_{1}$.
Indeed, formula~\eqref{eq:far-from-ends} implies, using~\eqref{eq:triangle-cnd.three}:
 \begin{align*}
 y_{i} > u_{1} + 0.1\f\cdot\slope(u,v)-3H\g
\ge u_{1}+\Paren{0.1(\slopeincr-1)-3H}\g \ge u_{1} +13\g.
 \end{align*}
This implies $u_{1}< y_{i}-L_{1}-\l_{1}$, using $\bub\ll\g$
from~\eqref{eq:bub-g-f}.
The inequality about $v_{1}$ is similar, using the other inequality
of~\eqref{eq:far-from-ends}.
It implies that there is a region 
$\clint{x_{i}}{x_{i}+5\bub} \times \clint{y}{y+\l_{1}}$ containing no traps,
with $\lint{y}{y+\l_{1}} \sbsq \lint{y_{i}}{y_{i} + L_{1}}$.
Thus, there is a $y$ in $\lint{y_{i}}{y_{i}+L_{1}-\l_{1}}$
such that $\clint{x_{i}}{x_{i}+5\bub} \times \clint{y}{y+\l_{1}}$ contains no
traps.
(In the present proof, all other arguments finding a region with no traps
in trap covers are analogous, so we will not mention
Lemma~\ref{lem:if-correl} explicitly again.)
However, all traps must start in a trap cover, so the region
$\clint{x_{i}-2\bub}{x_{i}+6\bub} \times \clint{y}{y+\l_{1}}$ contains no trap
either. 
Thus there are clean points 
$w'_{i}\in (x_{i}-\bub, y+\bub) + \clint{0}{\bub}^{2}$ and
$w''_{i}\in (x_{i}+4\bub, y+5\bub) + \clint{0}{\bub}^{2}$.
Note that $\minslope(w'_{i}, w''_{i}) \ge 1/2$, giving
$w'_{i} \leadsto w''_{i}$.
We have, using~\eqref{eq:d-diff} and $\slope(u,v) \le 1$:
 \begin{align}
\nonumber%\label{eq:x-free.trapc.vert.fw.x}
 -\bub &\le x'_{i}-x_{i} \le 0,
    & 4\bub &\le x''_{i} - x_{i} \le 5\bub,
\\\nonumber 
  \bub &\le y'_{i}-y_{i} \le 23\bub,
    & 5\bub &\le y''_{i}-y_{i} \le 27\bub,
\\%\label{eq:x-free.trapc.vert.fw.d}
  \bub &\le d(w'_{i}) -  d(w_{i})  \le 24\bub,
    &     0 &\le d(w''_{i}) - d(w_{i})  \le 27\bub.
 \end{align}
 Consider crossing a vertical trap cover backwards.
There is a $y$ in $\lint{y_{i}-L_{1}}{y_{i}-\l_{1}}$
such that the region
$\clint{x_{i}-2\bub}{x_{i}+6\bub} \times \clint{y}{y+\l_{1}}$ contains no trap.
There are clean points 
$w'_{i}\in (x_{i}-\bub, y+\bub) + \clint{0}{\bub}^{2}$ and
$w''_{i}\in (x_{i}+4\bub, y+5\bub) + \clint{0}{\bub}^{2}$ with
$\minslope(w'_{i}, w''_{i}) \ge 1/2$, giving $w'_{i}\leadsto w''_{i}$.
We have
\begin{align*}
 %\label{eq:x-free.trapc.vert.bw.x}
     -\bub &\le x'_{i}-x_{i} \le 0, 
            &   4\bub &\le x''_{i} - x_{i} \le 5\bub,
\\ -27\bub &\le y'_{i}-y_{i}  \le -5\bub,
            & -23\bub &\le y''_{i}-y_{i} \le -\bub,
\\%\label{eq:x-free.trapc.vert.bw.d}
   -27\bub &\le d(w'_{i}) -  d(w_{i})  \le -4\bub,
            & -28\bub &\le d(w''_{i}) - d(w_{i}) \le -\bub.
  \end{align*}
 \end{prooof} % trapc.vert

 \begin{step+}{step:three.x-free.trapc.horiz}
 Assume that the trap cover is horizontal.
 \end{step+}
 \begin{prooof} % trapc.horiz
 Consider crossing a horizontal trap cover forwards.
There is an $x$ in $\lint{x_{i}-L_{1}}{x_{i}-\l_{1}}$ such that
$\clint{x}{x+\l_{1}} \times \clint{y_{i}-2\bub}{y_{i}+6\bub}$ contains no
trap.
Thus there are clean points 
$w'_{i}\in (x+\bub, y_{i}-\bub) + \clint{0}{\bub}^{2}$ and
$w''_{i}\in (x+5\bub, y_{i}+4\bub) + \clint{0}{\bub}^{2}$ with
$w'_{i} \leadsto w''_{i}$.
We have similarly to the above, the inequalities
 \begin{align*}
%\label{eq:x-free.trapc.horiz.fw.x}
       -27\bub  &\le x'_{i} - x_{i}  \le -5\bub,
             &  -23\bub  &\le x''_{i} - x_{i} \le -\bub,
\\        -\bub &\le y'_{i}-y_{i}\le 0,
             &     4\bub &\le y''_{i} - y_{i} \le 5\bub,
\\%\label{eq:x-free.trapc.horiz.fw.d}
         -\bub  &\le d(w'_{i}) -  d(w_{i}) \le 27\bub, 
             &    4\bub  &\le d(w''_{i}) - d(w_{i}) \le 28\bub.
 \end{align*}
 Consider crossing a horizontal trap cover backwards.
There is an $x$ in $\lint{x_{i}}{x_{i}+L_{1}-\l_{1}}$ such that
$\clint{x}{x+\l_{1}} \times \clint{y_{i}-2\bub}{y_{i}+6\bub}$ contains no
trap.
Thus there are clean points 
$w'_{i}\in (x+\bub, y_{i}-\bub) + \clint{0}{\bub}^{2}$ and
$w''_{i}\in (x+5\bub, y_{i}+4\bub) + \clint{0}{\bub}^{2}$ with
$w'_{i} \leadsto w''_{i}$.
We again have
 \begin{align*}
%\label{eq:x-free.trapc.horiz.bw.x}
        \bub  &\le x'_{i} - x_{i}  \le 23\bub, 
         &  5\bub  &\le x''_{i} - x_{i} \le 27\bub,
\\      -\bub &\le y'_{i} - y_{i} \le 0,
         &   4\bub &\le y''_{i}-y_{i} \le 5\bub,
\\%\label{eq:x-free.trapc.horiz.bw.d}
     -24\bub  &\le d(w'_{i}) -  d(w_{i}) \le 0,   
         & -23\bub &\le d(w''_{i}) - d(w_{i}) \le 5\bub.
 \end{align*}
 \end{prooof} % trapc.horiz
 \end{prooof} % trapc

 \begin{step+}{step:three.x-free.wall}
 Consider crossing a wall.
 \end{step+}
 \begin{prooof} %  % wall
 \begin{step+}{step:three.x-free.wall.vert}
 Assume that the wall is vertical.
 \end{step+}
 \begin{prooof} % vert
 Consider crossing a vertical wall forwards.
Let us apply Lemma~\ref{lem:the-grate}\eqref{i:grate.hole}, with
$I'=\clint{y_{i}}{y_{i}+\g}$. 
The lemma is applicable since $w_{i}\in C(u, v, K+2\g, K+(2H-2)\g)$ implies
$u_{1} \le y_{i}-\g-7\bub < y_{i}+2\g+7\bub<v_{1}$.
It implies that our wall contains an outer upward-clean
hole $\rint{y'_{i}}{y''_{i}} \sbsq y_{i} + \rint{\bub}{\g-\bub}$ 
passing through it.
(In the present proof, all other arguments finding a hole through walls
are analogous, so we will not mention
Lemma~\ref{lem:the-grate}\eqref{i:grate.hole} explicitly again.)
Let $w'_{i} = \pair{x_{i}}{y'_{i}}$, and let $w''_{i}=\pair{x''_{i}}{y''_{i}}$
be the point on the other side of the wall reachable from $w'_{i}$.
We have
  \begin{align}
\nonumber%\label{eq:x-free.wall.vert.fw.x}
  x'_{i} &= x_{i}, 
      &      0 &\le x''_{i} - x_{i} \le \bub, 
\\  \bub &\le y'_{i}-y_{i} \le y''_{i}-y_{i} \le \g-\bub,
\\\nonumber%\label{eq:x-free.wall.vert.fw.d}
    \bub &\le d(w'_{i}) -  d(w_{i})  \le \g-\bub,
      &      0 &\le d(w''_{i}) - d(w_{i})  \le \g-\bub.
  \end{align}
 Consider crossing a vertical wall backwards.
This wall contains an outer upward-clean hole
$\rint{y'_{i}}{y''_{i}} \sbsq y_{i} + \rint{-\g+\bub}{-\bub}$
passing through it.
Let $w'_{i} = \pair{x_{i}}{y'_{i}}$, and let $w''_{i}=\pair{x''_{i}}{y''_{i}}$ 
be the point on the other side of the wall reachable from $w'_{i}$.
We have
  \begin{align*}
%\label{eq:x-free.wall.vert.bw.x}
     x'_{i} &= x_{i}, 
       &      0 &\le x''_{i} - x_{i} \le \bub, 
\\ -\g+\bub &\le y'_{i} -  y_{i} \le y''_{i}-y_{i} \le -\bub,
\\%\label{eq:x-free.wall.vert.bw.d}
   -\g+\bub &\le d(w'_{i}) -  d(w_{i})  \le -\bub,
       &    -\g &\le d(w''_{i}) - d(w_{i})  \le -\bub.
 \end{align*} 
%   -\g+\bub \le y'_{i} -  y_{i} \le y''_{i}-y_{i} \le -\bub.
 \end{prooof} % vert

 \begin{step+}{step:three.x-free.wall.horiz}
 Assume that the wall is horizontal.
 \end{step+}
 \begin{prooof} % horiz
 Consider crossing a horizontal wall forwards.
Similarly to above, this wall contains an outer rightwards clean hole 
$\rint{x'_{i}}{x''_{i}} \sbsq x_{i} + \rint{-\g+\bub}{-\bub}$ 
passing through it.
Let $w'_{i} = \pair{x'_{i}}{y_{i}}$ and let $w''_{i}=\pair{x''_{i}}{y''_{i}}$ 
be the point on the other side of the wall reachable from $w'_{i}$.
We have
  \begin{align*}
 -\g+\bub &\le x'_{i}-x_{i} \le x''_{i}-x_{i} \le -\bub,
\\ y'_{i} &= y_{i},
      &   0 &\le y''_{i}-y_{i} \le \bub,
\\%\label{eq:x-free.wall.horiz.fw.d}
        0 &\le d(w'_{i}) -  d(w_{i}) \le \g-\bub,
      &   0 &\le d(w''_{i}) - d(w_{i})  \le \g.
 \end{align*} 
 Consider crossing a horizontal wall backwards.
This wall contains an outer rightward-clean hole 
$\rint{x'_{i}}{x''_{i}} \sbsq x_{i} + \rint{\bub}{\g-\bub}$ 
passing through it.
Let $w'_{i} = \pair{x'_{i}}{y_{i}}$ and let $w''_{i}=\pair{x''_{i}}{y''_{i}}$ 
be the point on the other side of the wall reachable from $w'_{i}$.
We have
 \begin{align*}
%\nonumber
      \bub &\le x'_{i} -  x_{i} \le x''_{i}-x_{i} \le \g-\bub,
\\  y'_{i} &= y_{i},
      &        0 &\le y''_{i}-y_{i} \le \bub,
\\%\label{eq:x-free.wall.horiz.bw.d}
  -\g+\bub &\le d(w'_{i}) -  d(w_{i}) \le 0,
      & -\g+\bub &\le d(w''_{i}) - d(w_{i}) \le \bub.
 \end{align*} 
 \end{prooof} % horiz
 \end{prooof} % wall

 \begin{step+}{step:three.x-free.summary}
Let us summarize some of the inequalities proved above,
with 
 \[
  D=d(w)-d(w_{i}),
 \]
where $w$ is equal to any one of the defined $w'_{i}, w''_{i}$.
 \end{step+}
 \begin{prooof}% summary
 \begin{equation}\label{eq:free-bounds}
 \begin{aligned}
   \txt{trap covers going forwards:}  &&    -\bub &\le D \le 28\bub,
\\ \txt{trap covers going backwards:} &&  -28\bub &\le D \le 5\bub,
\\ \txt{walls going forwards:}        &&        0 &\le D \le \g,
\\ \txt{walls going backwards:}       &&      -\g &\le D \le \bub.
 \end{aligned}
 \end{equation}
Further
 \begin{equation}\label{eq:free-bounds-vert-horiz}
 \begin{aligned}
   \txt{vertical obstacles:}     &&
   -\bub &\le x'_{i}-x_{i} \le x''_{i}-x_{i} \le 5\bub,
\\ \txt{horizontal obstacles:}   &&
   -\bub &\le y'_{i}-y_{i} \le y''_{i}-y_{i} \le 5\bub,
\\ \txt{horizontal trap covers:} &&
  -27\bub &\le x'_{i}-x_{i} \le x''_{i}-x_{i} \le 27\bub,
\\ \txt{horizontal walls:}       &&
 -\g+\bub &\le x'_{i}-x_{i} \le x''_{i}-x_{i} \le \g-\bub.
 \end{aligned}
 \end{equation}
Let $\pi_{x} w,\pi_{y} w\in\bbR$ be the $X$ and $Y$ projections of point $w$,
and let $\pi_{i} w\in\bbR^{2}$ be the projection of point $w$
onto the (horizontal or vertical) line corresponding to $s_{i}$.
Then the above inequalities and~\eqref{eq:d-diff} imply, 
with $\hat w=\pi_{i}w - w$ where $w=w'_{i},w''_{i}$:
 \begin{equation}\label{eq:projection}
  -5\bub \le d(\pi_{i}w)-d(w), \pi_{x}\hat w, \pi_{y}\hat w \le 5\bub.
 \end{equation}
Indeed, for example for $\pi_{x}\hat w$, for a horizontal wall or trap cover,
the difference between the projection of $w$ and the projection of the
projection onto the wall is 0.
For a vertical wall or trap cover, one of $w'_{i},w''_{i}$ 
is at a distance at most $\bub$ from the line corresponding to $s_{i}$, the
other one is inside the trap cover or within $\bub$ of the other side of the
wall, therefore is at most at a distance $5\bub$.
The inequality for $d(\pi_{i}w)-d(w)$ follows from this, \eqref{eq:d-diff} and that
only one of the coordinates changes.

For crossing a wall we have
 \begin{equation}\label{eq:wall-across}
   -\bub \le d(w''_{i}) - d(w'_{i}) \le \bub.
 \end{equation}
Indeed, $w'_{i}$ and $w''_{i}$ are the two opposite corners of a hole through a
wall of width $\le\bub$.
 \end{prooof}% summary
 \end{prooof} % x-free

 \begin{step+}{step:three.free}
 Assume that there is no bound pair: then $u\leadsto v$.
 \end{step+} 
 \begin{pproof} % free
 \begin{step+}{step:three.free.fw}
 Assume that there is no horizontal trap-wall pair.
 \end{step+}
 \begin{prooof} % free.fw
We choose $w_{1}$ with $d(w_{1})=K + 3\g$.
For each $i>1$ we choose $w_{i}$ with $d(w_{i})=d(w''_{i-1})$, and we cross
each obstacle in the forward direction.

  \begin{step+}{step:three.free.fw.d-bds}
For all $i$, the points we created are inside a certain channel:
 \begin{equation}\label{eq:d-global-bds}
   d(w), d(\pi_{i} w) \in K+\clint{2\g}{(2H-2)\g},
 \end{equation}
where $w$ is any one of $w_{i},w'_{i},w''_{i}$.
  \end{step+}
  \begin{pproof} % d-bds
It follows from~\eqref{eq:free-bounds} that,
for $w \in \{w'_{i}, w''_{i}\}$:
 \begin{align}
\label{eq:d-trap-fw-incr}
      -\bub &\le d(w) - d(w_{i}) \le 28\bub &&\text{for trap covers},
\\\label{eq:d-wall-fw-incr}
           0 &\le d(w) - d(w_{i}) \le \g    &&\text{for walls}.
 \end{align}
To estimate for example $d(w_{3})$, we write
 \begin{multline*}
      d(w''_{3})=d(w_{1})+d(w_{2})-d(w_{1})+d(w_{3})-d(w_{2})+d(w''_{3})-d(w_{3})
\\              =K+3\g + (d(w''_{1})-d(w_{1}))+(d(w''_{2})-d(w_{2}))+(d(w''_{3})-d(w_{3})),
 \end{multline*}
were we used $d(w_{1})=K+3\g$.
Since we have two walls and a trap cover, using~\eqref{eq:d-trap-fw-incr} once
and~\eqref{eq:d-wall-fw-incr} twice gives
  \[
   K+3\g-\bub \le d(w''_{3}) \le K + 5\g + 28\bub.
 \]
The same argument works for all $w_{i},w'_{i},w''_{i}$.
Then~\eqref{eq:projection} implies
 $K+3\g-6\bub \le d(w),d(\pi_{i} w) \le K + 5\g + 33\bub$,
where $w$ is any one of $w_{i},w'_{i},w''_{i}$.
  \end{pproof} % d-bds

 \begin{step+}{step:three.free.fw.dist-lb}
For $i=2,3$ the inequality $x'_{i}-x''_{i-1} \ge \g$ holds.
 \end{step+}
 \begin{pproof}
If $s_{i-1},s_{i}$ come from trap covers or walls in
different orientations, then the intersection of their lines lies outside
$C(u,v,K,K+2H\g)$.
Part~\ref{step:three.free.fw.d-bds} above says 
 \[
  \pi_{i}w'_{i}, \pi_{i-1}w''_{i-1}\in C(u,v,K+2\g,K+(2H-2)\g).
 \]
Now if two segments $A$, $B$ of different orientation (horizontal and
vertical) are contained in $C(u,v,K+2\g,K+(2H-2)\g)$ and
are such that $A$ is to the left of $B$ and their lines
intersect outside $C(u,v,K,K+2H\g)$, then for any points $a\in A$, $b\in B$
we have $\pi_{x}(b-a)\ge 2\g$.
In particular, $\pi_{x}\pi_{i}w'_{i} - \pi_{x}\pi_{i-1} w''_{i-1} \ge 2\g$.
Using~\eqref{eq:projection} we get:
\begin{align*}
  x'_{i}-x''_{i-1} &= \pi_{x}\pi_{i}w'_{i} - \pi_{x}\pi_{i-1} w''_{i-1}
 + (x'_{i}-\pi_{x}\pi_{i}w'_{i}) - (x''_{i-1}-\pi_{x}\pi_{i-1} w''_{i-1}) 
\\ &\ge 2\g - 10\bub.
 \end{align*}
If $s_{i-1},s_{i}$ 
come from a vertical trap-wall or wall-trap pair, then freeness
implies that elements of this pair are farther than $1.1\g$ from each
other.
Then we get $x'_{i}-x''_{i-1} \ge 1.1\g-5\bub > \g$, as can be seen by
considering the two possible orders: trap-wall, wall-trap separately.

If $s_{i-1},s_{i}$ come from a horizontal wall-trap pair then, using
$\slope(u,v) \le 1$ and~\eqref{eq:free-bounds-vert-horiz}:
\begin{align*}
  x_{i}-x''_{i-1} &= (y_{i} - y''_{i-1})/\slope(u,v) \ge y_{i} - y''_{i-1}
  = y_{i}-y_{i-1} - (y''_{i-1}-y_{i-1}) 
\\ &\ge 1.1\g-5\bub.
\end{align*}
By~\eqref{eq:free-bounds-vert-horiz} we have $x'_{i}-x_{i} \ge -27\bub$.
Combination with the above estimate and~\eqref{eq:bub-g-f} gives
$x'_{i}-x''_{i-1} \ge 1.1\g-32\bub > \g$ due to~\eqref{eq:bub-g-f}.
 \end{pproof} % dist-lb

 \begin{step+}{step:three.free.fw.reach}
Let us show $u\leadsto v$.
 \end{step+}
 \begin{pproof} % reach
We defined all $w'_{i},w''_{i}$ as clean points with
$w'_{i} \leadsto w''_{i}$ and the sets 
$C^{\eps_{i}}(w''_{i},w'_{i+1}, -\g,\g)$ are trap-free, where
$\eps_{i}=\rightarrow$ for horizontal walls, $\uparrow$ for vertical walls
and nothing for trap covers.
If we are able to show that the minslopes between the endpoints of these sets
are lower-bounded by $\slb + 4\g / \f$
then Lemma~\ref{lem:approx.no-walls} will imply $u\leadsto v$.
For this, it will be sufficient to show that the 
slopes are lower-bounded by $\slb + 4\g / \f$ and upper-bounded by 2,
since~\eqref{eq:bub-g-f} implies $1/(\slb + 4\g/f)>2$.
We will make use of the following relation 
for arbitrary $a=\pair{a_{0}}{a_{1}}$, $b=\pair{b_{0}}{b_{1}}$:
 \begin{equation}\label{eq:slope-by-d}
   \slope(a, b) = \slope(u,v) + \frac{d(b)-d(a)}{b_{0}-a_{0}}.
 \end{equation}
Let us bound the end slopes first.
We have
 \begin{equation}\label{eq:end-slopes}
         \slope(u, w'_{1}) = \slope(u,v) + \frac{d(w'_{1})}{x'_{1}-u_{0}},
\; \slope(w''_{3}, v) = \slope(u,v) -\frac{d(w''_{3})}{v_{0}-x''_{3}}.
 \end{equation}
Inequalities~\eqref{eq:d-global-bds} yield the bounds
$K+2\g\le d(w'_{1}),d(w''_{3})\le K + (2H-2)\g$, hence the restriction of $K$
to $\clint{-3H\g}{H\g}$ implies 
$|d(w'_{1})|, |d(w''_{3})|\le (3H-2)\g$.
By~\eqref{eq:far-from-ends} and~\eqref{eq:bub-g-f}, using 
also~\eqref{eq:free-bounds-vert-horiz}:
 \begin{equation*}
 (x'_{1}-u_{0}),\; (v_{0}-x''_{3}) \ge 0.1\f-\g+\bub \ge \f/11.
 \end{equation*}
This shows
 \begin{equation}
    |\slope(u, w'_{1}) - \slope(u,v)|,\;
    |\slope(w''_{3}, v) - \slope(u,v)| \le 33H\g/\f.
 \end{equation}
Using $1\ge\slope(u,v) \ge \slb+(\slopeincr-1)\g/\f$:
 \begin{align*}
  2 > 1+33H\g/\f &\ge \slope(u, w'_{1}),\; \slope(w''_{3}, v) 
\\           &\ge \slb + (\slopeincr - 1 - 33H)\g/f \ge \slb + 4\g/\f  % was 6
 \end{align*}
by~\eqref{eq:bub-g-f} and~\eqref{eq:H-slopeincr}.

Let us proceed to lower-bounding $\minslope(w''_{i-1}, w'_{i})$ for
$i=2,3$.
Using~\eqref{eq:slope-by-d}:
 \begin{equation}\label{eq:mid-slope-expr}
   \slope(w''_{i-1}, w'_{i}) = \slope(u,v) +
  \frac{d(w'_{i}) - d(w''_{i-1})}{x'_{i}-x''_{i-1}}.
 \end{equation}
Using~\eqref{eq:free-bounds} and Part~\ref{step:three.free.fw.dist-lb} above:
 \[
  -5\bub/\g \le \slope(w''_{i-1}, w'_{i}) - \slope(u,v) \le 1.
 \]
(For just the case of crossing forwards, we could have used the lower bound
$-\bub$, but a similar application below to backward crossing requires $-5\bub$.)
By the conditions of the lemma and~\eqref{eq:bub-g-f}:
 \begin{align*}
   \slope(w''_{i-1}, w'_{i}) &\ge \slope(u,v)-5\bub/\g 
   \ge \slope(u,v)-5\g/\f 
\\             &\ge \slb + (\slopeincr-6) \g / \f 
             \ge \slb + 4\g/\f, 
\\ \slope(w''_{i-1}, w'_{i}) &\le \slope(u,v)+ 1 <  2.
 \end{align*}
 \end{pproof} % reach
 \end{prooof} % free.fw

 \begin{step+}{step:three.free.bw}
 Assume now that there is a horizontal trap-wall pair.
 \end{step+}
 \begin{prooof} % free.bw
What has been done in part~\ref{step:three.free.fw} will be
repeated, going backwards through $i=3,2,1$ rather than forwards.
Thus, we choose $w_{3}$ with $d(w_{3}) = d(v)+K+(2H-3)\g$.
Assuming that $w'_{i+1}$ has been chosen already, we choose
$w_{i}$ with $d(w_{i})=d(w'_{i+1})$, and we cross each obstacle in the
backward direction.

It follows from~\eqref{eq:free-bounds}, that
for all $i$ we have~\eqref{eq:d-global-bds} again.

 \begin{step+}{step:three.free.bw.dist-lb}
The inequality $x'_{i}-x''_{i-1} \ge \g$ holds.
% \label{eq:free.bw.x-lb}
 \end{step+}
 \begin{pproof}
If $s_{i-1}$ and $s_{i}$ come from trap covers or walls in
different orientations, then
we can reason as in Part~\ref{step:three.free.fw.dist-lb}.
If $s_{i-1},s_{i}$ come from a horizontal trap-wall pair then
$d(w'_{i}) = d(w_{i-1})$ and $y'_{i}=y_{i}$ imply
 \[
    x'_{i}-x_{i-1} = (y_{i} - y_{i-1})/\slope(u,v) \ge 1.1\g.
 \]
By~\eqref{eq:free-bounds-vert-horiz}, 
$x_{i-1}-x''_{i-1} \ge -27\bub$, hence
$x'_{i}-x''_{i-1} \ge 1.1\g-27\bub \ge\g$.
 \end{pproof} % dist-lb

 \begin{step+}{step:three.free.bw.reach}
We have $u\leadsto v$.
 \end{step+}
 \begin{pproof} % reach
Let us estimate the minslopes as in part~\ref{step:three.free.fw.reach},
using~\eqref{eq:end-slopes} again.
The estimates for $\slope(u, w'_{1})$ and $\slope(w''_{3},v)$
are as before.
We conclude for  $\minslope(w''_{i-1}, w'_{i})$ 
using~\eqref{eq:mid-slope-expr} and 
Part~\ref{step:three.free.bw.dist-lb} just as we did in
Part~\ref{step:three.free.fw.reach} above.
 \end{pproof} % reach
 \end{prooof} % free.bw
 \end{pproof} % free

 \begin{step+}{step:three.x-bound}
Consider crossing a bound pair.
 \end{step+}
 \begin{prooof}
A bound trap-wall or wall-trap pair will be crossed with an approximate
slope 1 rather than $\slope(u,v)$.
We first find a big enough (size $\l_{2}=\g'$) hole in the trap cover, and
then locate a hole in the wall that allows to pass, with slope 1,
through the big hole of the trap cover.
There are cases according to whether we have a trap-wall pair or a
wall-trap pair, and whether it is vertical or horizontal.

We will prove
 \begin{equation}\label{eq:bd-bounds}
 -1.2\g \le d(w) - d(w_{i}) \le 10\g,
 \end{equation}
for $w$ equal to any one of the defined $w'_{j},w''_{j}$ where
$j\in\{i-1,i\}$ (wall-trap) and $j\in\{i,i+1\}$ (trap-wall).
The inequalities~\eqref{eq:free-bounds-vert-horiz} 
and~\eqref{eq:projection} will hold also if the obstacle is within a bound
pair.

  \begin{step+}{step:three.x-bound.trapc} 
 Consider crossing a trap-wall pair $(i,\,i+1)$, assuming that
$w_{i}$ has been defined already.
  \end{step+}
  \begin{prooof}
  \begin{step+}{step:three.x-bound.trapc.vert}
 Assume that the trap-wall pair is vertical.
  \end{step+}
 \begin{prooof}
Let us apply Lemma~\ref{lem:if-correl} with $j=2$, so taking $\l_{2}=\g'=2.2\g$,
$L_{2}=4\l_{2}$
similarly to the forward crossing in Part~\ref{step:three.x-free.trapc.vert}.
We find 
a $y^{(1)}$ in $\lint{y_{i}}{y_{i}+L_{2}-\l_{2}}$ such that the region
$\clint{x_{i}-3\bub}{x_{i+1}} \times \clint{y^{(1)}}{y^{(1)}+\l_{2}}$
contains no trap.

Let $w_{i+1}$ be defined by  $y_{i+1}= y^{(1)} + (s_{i+1}-s_{i}) + 2\bub$.
Thus, it is the point on the left edge of the wall if we intersect
it with a slope 1 line from $\pair{s_{i}}{y^{(1)}}$ and then move up $2\bub$.
Similarly to the forward crossing in Part~\ref{step:three.x-free.wall.vert}, the wall  
starting at $s_{i+1}$ is passed through by an outer upward-clean hole 
$\rint{y'_{i+1}}{y''_{i+1}} \sbsq y_{i+1} + \rint{\bub}{\g-\bub}$.
Let $w'_{i+1} = \pair{x_{i+1}}{y'_{i+1}}$, and let $w''_{i+1}=\pair{x''_{i+1}}{y''_{i+1}}$ 
be the point on the other side of the wall reachable from $w'_{i+1}$.

Let $w = \pair{x_{i}}{y^{(2)}}$ be defined by $y^{(2)}= y'_{i+1} - (s_{i+1}-s_{i})$.
Thus, it is the point on the left edge of the trap cover if we intersect
it with a slope 1 line from $w'_{i+1}$.
Then $3\bub \le y^{(2)}-y^{(1)}\le\g+\bub$, therefore $w +
\clint{-3\bub}{0}^{2}$ contains no trap, and
there is a clean point $w'_{i}\in w+\clint{-2\bub}{-\bub}^{2}$.
(Point $w''_{i}$ is not needed.)

Let us estimate $d(w'_{i})-d(w_{i})$ and $d(w''_{i+1})-d(w_{i})$.
Recalling $\l_{2}=\g'$, $L_{2}=4\g'$ gives
 \begin{align}
\label{eq:three.x-bound.trapc.vert.dw}
             3\bub &\le y^{(2)}-y_{i} = d(w) - d(w_{i}) \le 3\g'+\g+\bub, 
\\\nonumber
            -2\bub &\le d(w'_{i}) - d(w) \le 2\bub,
\\\nonumber      0 &\le d(w'_{i+1}) - d(w)  \le 1.1\g.
 \end{align}
Combining the last inequalities with~\eqref{eq:three.x-bound.trapc.vert.dw}
gives 
 \begin{align*}
     \bub  &\le d(w'_{i}) - d(w_{i}) \le 3\g'+ \g + 3\bub,
\\   3\bub &\le d(w'_{i+1}) - d(w_{i}) \le 3\g'+2.1\g + \bub.
 \end{align*}
These prove~\eqref{eq:bd-bounds} for our case
if we also invoke~\eqref{eq:wall-across} to infer about $w''_{i+1}$.
Let us show $w'_{i} \leadsto w'_{i+1}$.
We apply Condition~\ref{cond:reachable}.
It is easy to see that
the rectangle $\Rect^{\eps_{i}}(w'_{i},w'_{i+1})$ is trap-free.
Consider the slope condition.
We have $1/2 \le \slope(w'_{i}, w) \le 2$, and 
$\slope(w, w'_{i+1}) = 1$.
Hence, $1/2 \le \slope(w'_{i}, w'_{i+1}) \le 2$, which implies
$\minslope(w'_{i}, w'_{i+1}) \ge \slb + 4 \g / \f$
as needed.
 \end{prooof} % trapc.vert

 \begin{step+}{step:three.x-bound.trapc.horiz}
 Assume that the trap-wall pair is horizontal.
 \end{step+}
 \begin{prooof} % trapc.horiz
There is an $x^{(1)}$ in $\lint{x_{i}-L_{2}}{x_{i}-\l_{2}}$
such that the region 
$\clint{x^{(1)}}{x^{(1)}+\l_{2}} \times \clint{y_{i}-3\bub}{y_{i+1}}$
contains no trap.
Let $w_{i+1}$ be defined by $x_{i+1}= x^{(1)} + (s_{i+1}-s_{i})\slope(u,v) + 2\bub$.
The wall starting at $s_{i+1}$ contains an outer rightward-clean hole 
$\rint{x'_{i+1}}{x''_{i+1}} \sbsq x_{i+1} + \rint{\bub}{\g-\bub}$ 
passing through it.
Let $w'_{i+1} = \pair{x'_{i+1}}{y_{i+1}}$, and let 
$w''_{i+1}=\pair{x''_{i+1}}{y''_{i+1}}$ 
be the point on the other side of the wall reachable from $w'_{i+1}$.
Let $w = \pair{x^{(2)}}{y_{i}}$ be defined by 
$x^{(2)}= x'_{i+1} - (s_{i+1}-s_{i})\slope(u,v)$.
Then there is a clean point $w'_{i}\in w+\clint{-2\bub}{-\bub}^{2}$ as before.
We have 
 \begin{equation*}
        3\bub\le x^{(2)}-x^{(1)}\le \g+\bub,
\quad
  - 4\g' + 3\bub \le x^{(2)} - x_{i} \le - \g' + \g + \bub. 
 \end{equation*}
Using this, $d(w)-d(w_{i}) = -(x^{(2)}-x_{i})\slope(u,v)$ and
$\slope(u,v) \le 1$ gives
 \begin{equation*}
  0 \le d(w) - d(w_{i}) \le  4\g' - 3\bub.
 \end{equation*}
As in Part~\ref{step:three.x-bound.trapc.vert}, this gives
 \begin{align*}
   -2\bub &\le d(w'_{i}) - d(w_{i})   \le 4\g' - \bub,
\\    0   &\le d(w'_{i+1}) - d(w_{i}) \le 4\g' + 1.1\g - 3\bub.
 \end{align*}
These and~\eqref{eq:wall-across} prove~\eqref{eq:bd-bounds} for our case.
We show $w'_{i}\leadsto w'_{i+1}$ similarly to
Part~\ref{step:three.x-bound.trapc.vert}.
 \end{prooof} % trapc.horiz
 \end{prooof} % trapc

 \begin{step+}{step:three.x-bound.wall}
 Consider crossing a wall-trap pair $(i-1,\,i)$, assuming that
$w_{i}$ has been defined already.
 \end{step+}
 \begin{prooof} %  % wall
 \begin{step+}{step:three.x-bound.wall.vert}
 Assume that the wall-trap pair is vertical.
 \end{step+}
 \begin{prooof} % vert
This part is somewhat similar to Part~\ref{step:three.x-bound.trapc.vert}, and
is illustrated in Figure~\ref{fig:close-trap}.
There is a $y^{(1)}$ in $\lint{y_{i}+\l_{2}}{y_{i}+L_{2}}$
such that the region 
$\clint{x_{i}}{x_{i}+6\bub} \times \clint{y^{(1)}-\l_{2}}{y^{(1)}}$
contains no trap.
Let $w_{i-1}$ be defined by  $y_{i-1}= y^{(1)} - (s_{i}-s_{i-1}) - 5\bub$.
The wall starting at $s_{i-1}$ contains an outer upward-clean hole
$\rint{y'_{i-1}}{y''_{i-1}} \sbsq y_{i-1} + \rint{-\g + \bub}{-\bub}$
passing through it. 
We define $w'_{i-1}$, and $w''_{i-1}$ accordingly.
Let $w = \pair{x_{i}}{y^{(2)}}$ where $y^{(2)}= y'_{i-1} + (s_{i}-s_{i-1})$.
There is a clean point
$w''_{i}\in w+(4\bub, 4\bub) + \clint{0}{\bub}^{2}$.
We have
 \begin{align}
 \nonumber 
           -\g-4\bub &\le  y^{(2)} - y^{(1)}\le - 6\bub,
\\\label{eq:three.x-bound.wall.vert.dw}
   \g' - \g - 4\bub &\le y^{(2)} - y_{i} = d(w) - d(w_{i}) \le 4\g'- 6\bub.
\\\label{eq:three.x-bound.wall.vert.w''(i)}
             -\bub  &\le d(w''_{i}) - d(w) \le 5\bub,
\\\label{eq:three.x-bound.wall.vert.w'(i-1)}
             -1.1\g &\le d(w'_{i-1}) - d(w) \le 0.
 \end{align}
Combining the last two inequalities with~\eqref{eq:three.x-bound.wall.vert.dw} 
gives 
 \begin{align*}
    \g'-2.1\g-4\bub &\le d(w'_{i-1}) - d(w_{i}) \le 4\g'- 6\bub, 
\\    \g'-\g -5\bub &\le d(w''_{i}) - d(w_{i}) \le 4\g'-\bub.    
 \end{align*}
These and~\eqref{eq:wall-across} prove~\eqref{eq:bd-bounds} for our case.
The reachability $w''_{i-1}\leadsto w''_{i}$
is shown similarly to Part~\ref{step:three.x-bound.trapc.vert}.
For this note
 \begin{align*}
      y_{i-1}   &\ge y^{(1)} - 1.1\g - 5\bub,
\\    y''_{i-1} &\ge y^{(1)} - 2.1\g - 4\bub \ge y^{(1)}-2.2\g = y^{(1)}-\g'.
 \end{align*}
This shows that the rectangle $\Rect^{\eps}(w''_{i-1}, w''_{i})$ is
trap-free.
The bound $1/2 \le \minslope(w''_{i-1}, w''_{i})$ is easy to check.
 \end{prooof} % vert

 \begin{step+}{step:three.x-bound.wall.horiz}
 Assume that the wall-trap pair is horizontal.
 \end{step+}
 \begin{prooof} % horiz
This part is somewhat similar to Parts~\ref{step:three.x-bound.trapc.horiz}
and~\ref{step:three.x-bound.wall.vert}.
There is an $x^{(1)}$ in $\lint{x_{i}-L_{2}+\l_{2}}{x_{i}}$
such that the region 
$\clint{x^{(1)}-\l_{2}}{x^{(1)}} \times \clint{y_{i}}{y_{i}+6\bub}$
contains no trap.
Let $w_{i-1}$ be defined by $x_{i-1}= x^{(1)} - (s_{i}-s_{i-1})\slope(u,v) - 5\bub$.
The wall starting at $s_{i-1}$ contains an outer rightward-clean hole 
$\rint{x'_{i-1}}{x''_{i-1}} \sbsq x_{i-1} + \rint{-\g+\bub}{-\bub}$ 
passing through it.
We define $w'_{i-1}$, $w''_{i-1}$ accordingly.
Let $w = \pair{x^{(2)}}{y_{i}}$ where $x^{(2)}= x'_{i-1} + (s_{i}-s_{i-1})\slope(u,v)$.
There is a clean point 
$w''_{i}\in w+(4\bub, 4\bub) + \clint{0}{\bub}^{2}$. 
We have
 \begin{equation*}
         x^{(2)} \in x^{(1)} - 5\bub + \rint{-\g+\bub}{-\bub},
\quad
  - 3\g' -\g - 4\bub \le x^{(2)} - x_{i} \le -6\bub. 
 \end{equation*}
This gives
 $0  \le d(w) - d(w_{i}) \le 3\g' + \g + 4\bub$. 
Combining with~\eqref{eq:three.x-bound.wall.vert.w''(i)}  
and~\eqref{eq:three.x-bound.wall.vert.w'(i-1)} 
which holds just as in Part~\ref{step:three.x-bound.wall.vert}, we get
 \begin{align*}
          -1.1\g &\le d(w'_{i-1}) - d(w_{i}) \le 3\g' + \g + 4\bub,
\\        -\bub  &\le d(w''_{i}) - d(w_{i}) \le 3\g' + \g + 9\bub.
 \end{align*}
These and~\eqref{eq:wall-across} prove~\eqref{eq:bd-bounds} for our case.
The reachability $w''_{i-1}\leadsto w''_{i}$
is shown similarly to Part~\ref{step:three.x-bound.wall.vert}.
 \end{prooof} % horiz
 \end{prooof} % wall
 \end{prooof} % x-bound

 \begin{step+}{step:three.bound}
 Assume that there is a bound pair: then $u\leadsto v$.
 \end{step+}
 \begin{pproof} % bound
We define $w_{i}$ with
 \begin{align}\label{eq:wi}
 d(w_{i}) = K + 5\g    
 \end{align}
if we have a trap-wall pair $(i,\,i+1)$ or a wall-trap
pair $(i-1,\,i)$.
As follows from the starting discussion of the proof of the lemma,
the third obstacle, outside the bound pair, is a wall.

\begin{figure}
\centering
\includegraphics[scale=1]{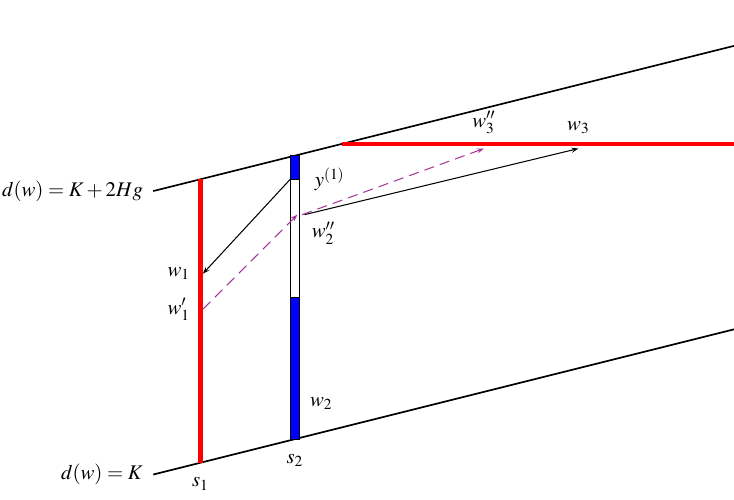}
 \caption{Approximation Lemma: the case of a bound wall-trap pair $(1,2)$.
The arrows show the order of selection.
First $w_{2}$ is defined.
Then the trap-free segment of size $\g'$ above $w_{2}$ is found.
Its starting point $y^{(1)}$ is projected back by a slope 1 line
onto the vertical wall to find $w_{1}$ after moving down by $5\bub$.
The hole starting with $w'_{1}$ is found within $\g$ below $w_{1}$.
Then $w''_{2}$ is found near the back-projection of $w''_{1}$.
Then $w''_{2}$ is projected forward, by a $\slope(u,v)$ line onto the
horizontal wall, to find $w_{3}$.
Finally, the hole ending in $w''_{3}$ is found within $\g$ backwards from
$w_{3}$.}\label{fig:close-trap}
\end{figure}

  \begin{step+}{step:three.bound.fw}
Assume that the bound pair is $(1,2)$.
  \end{step+}
 \begin{prooof} % fw
  \begin{step+}{step:three.bound.fw.trapc}
Assume that we have a trap-wall pair.
  \end{step+}
 \begin{prooof} % trapc
We defined $d(w_{1}) = K + 5\g$,
further define $w'_{1},w''_{2}$ as in Part~\ref{step:three.x-bound.trapc}, and
$w_{3},w'_{3},w''_{3}$ as in Part~\ref{step:three.free.fw}.
Let us show that these points do not leave $C(u,v,K+2\g,K+(2H-2)\g)$:
for all $i$, we have~\eqref{eq:d-global-bds}.
Inequalities~\eqref{eq:bd-bounds} imply
 \begin{equation*}
    -1.2\g \le d(w'_{1}) - d(w_{1}),\;
 d(w''_{2}) - d(w_{1}) = d(w_{3}) - d(w_{1}) \le 10\g,  % was 8
 \end{equation*}
while inequalities~\eqref{eq:free-bounds} imply
 $0 \le d(w'_{3}) - d(w_{3}),\; d(w''_{3}) - d(w_{3})  \le \g$.    
Combining with~\eqref{eq:wi} gives for
$w\in\set{w'_{1},w''_{2},w'_{3},w''_{3}}$:
 \[
   K+3.8\g \le d(w)\le K + 16\g < K + (2H - 2)\g
 \]
according to~\eqref{eq:H-def}.

We have shown 
$w'_{1}\leadsto w''_{2}$ and $w'_{3}\leadsto w''_{3}$, further such that
the sets $C^{\eps_{1}}(u,w'_{1}, -\g,\g)$,
$C^{\eps_{2}}(w''_{2}, w'_{3}, -\g,\g)$ and 
$C^{\eps_{3}}(w''_{3}, v, -\g,\g)$ for the chosen $\eps_{i}$
are trap-free.
It remains to show that the minslopes between the endpoints of these sets
are lower-bounded by $\slb + 4\g / \f$:
then a reference to Lemma~\ref{lem:approx.no-walls} will imply $u\leadsto v$.
This is done for all three pairs $\pair{u}{w'_{1}}$, $\pair{w''_{3}}{v}$ 
and $\pair{w''_{2}}{w'_{3}}$ just as in Part~\ref{step:three.free.fw.reach}.
 \end{prooof} % trapc

  \begin{step+}{step:three.bound.fw.wall}
Assume that we have a wall-trap pair.
  \end{step+}
 \begin{prooof} % wall
We defined $d(w_{2}) = K + 5\g$; we further 
define $w'_{1},w''_{2}$ as in Part~\ref{step:three.x-bound.wall}, and
$w_{3},w'_{3},w''_{3}$ as in Part~\ref{step:three.free.fw}.
The proof is finished similarly to Part~\ref{step:three.bound.fw.trapc}.
 \end{prooof} % wall
 \end{prooof} % fw

  \begin{step+}{step:three.bound.bw}
Assume now that the bound pair is $(2,3)$.
  \end{step+}
 \begin{prooof} % bw
  \begin{step+}{step:three.bound.bw.trapc}
Assume that we have a trap-wall pair.
  \end{step+}
 \begin{prooof} % trapc
We defined $d(w_{2}) = K + 5\g$; we further
define $w'_{2},w''_{3}$ as in Part~\ref{step:three.x-bound.trapc}, and
$w_{1},w'_{1},w''_{1}$ as in Part~\ref{step:three.free.bw}.
Let us show that these points do not leave $C(u,v,K+2\g,K+(2H-2)\g)$.
Inequalities~\eqref{eq:bd-bounds} imply
 \begin{equation*}
    -1.2\g \le d(w'_{2}) - d(w_{2}) = d(w_{1}) - d(w_{2}),\; 
               d(w''_{3}) - d(w_{2}) \le 10\g,
 \end{equation*}
while inequalities~\eqref{eq:free-bounds} imply
 $-\g \le d(w'_{1}) -  d(w_{1}),\; d(w''_{1})-d(w_{1}) \le \bub$.
Combining with~\eqref{eq:wi} gives for
$w\in\set{w'_{1},w''_{1},w'_{2},w''_{3}}$:
 \begin{equation*}
   K+2.8\g \le d(w)\le K + 15\g + \bub < K + (2H - 2)\g.
 \end{equation*}
Reachability is proved as in Part~\ref{step:three.bound.fw.trapc}.
 \end{prooof} % trapc

  \begin{step+}{step:three.bound.bw.wall}
Assume that we have a wall-trap pair.
  \end{step+}
 \begin{prooof} % wall
We defined $d(w_{3}) = K + 5\g$; we further
define $w'_{2},w''_{3}$ as in Part~\ref{step:three.x-bound.wall}, and
$w_{1},w'_{1},w''_{1}$ as in Part~\ref{step:three.free.bw}.
The proof is finished as in Part~\ref{step:three.bound.bw.trapc}.
 \end{prooof} % wall
 \end{prooof} % bw
  \end{pproof} % bound
\end{Proof}

\begin{Proof}[Proof of Lemma~\protect\ref{lem:approx}(Approximation)]
Recall that the lemma says that if a rectangle $Q=\Rect^{\eps}(u,v)$ contains no
walls or traps of $\cM^{*}$, is inner clean in $\cM^{*}$ and
satisfies the slope lower bound condition
$\minslope(u,v)\ge\slb^{*}$ for $\cM^{*}$, then $u\leadsto v$.

The proof started by recalling, in Lemma~\ref{lem:the-grate},
that walls of $\cM$ in $Q$ can be grouped to a
horizontal and a vertical sequence, whose members are well separated from each
other and from the sides of $Q$.
Then it showed, in Lemma~\ref{lem:trap-cover},
that all traps of $\cM$ are covered by certain
horizontal and vertical stripes called trap covers.
Walls of $\cM$ and trap covers were called obstacles.

Next it showed, in Lemma~\ref{lem:approx.no-walls}, that in case there are no
traps or walls of $\cM$ in $Q$ then there is a path through $Q$ that stays close
to the diagonal.

Next, a series of obstacles (walls or trap covers) was defined,
along with the points $s_{1},s_{2},\dots$ that are obtained by the intersection points of
the obstacle with the diagonal, and projected to the $x$ axis.
It was shown in Lemma~\ref{lem:one-of-three} that these obstacles are well separated
into groups of up to three.
The most laborious lemma of the paper,
Lemma~\ref{lem:three}, showed how to pass each triple of obstacles.
It remains to conclude the proof.

For each pair of numbers $s_{i}, s_{i+1}$ with 
$s_{i+1} - s_{i} \ge 0.22\f$, define its midpoint $(s_{i} + s_{i+1})/2$.
Let $t_{1} < t_{2} < \dots < t_{n}$ be the sequence of all these
midpoints.
Let us define the square
 \[
  S_{i} = \pair{t_{i}}{u_{1} + \slope(u, v)(t_{i} - u_{0})} + 
  \clint{0}{\bub} \times \clint{-\bub}{0}.
 \]
By Remark~\ref{rem:distr}.\ref{i:middle-third}, each of these squares
contains a clean point $p_{i}$.

 \begin{step+}{approx.middle}
For $1 \le i < n$, the 
rectangle $\Rect(p_{i}, p_{i+1})$ satisfies the conditions
of Lemma~\ref{lem:three}, and therefore $p_{i}\leadsto p_{i+1}$.
The same holds also for $\Rect^{\eps}(u,p_{1})$
if the first obstacle is a wall, and for $\Rect(p_{n},v)$
if the last obstacle is a wall.
Here $\eps=\uparrow,\to$ or nothing, depending on the nature of the original
rectangle $\Rect^{\eps}(u,v)$.
 \end{step+}
 \begin{pproof} % middle
By Lemma~\ref{lem:one-of-three}, there are at most
three points of $\set{s_{1}, s_{2}, \dots}$ between
$t_{i}$ and $t_{i+1}$.
Let these be $s_{j_{i}}, s_{j_{i}+1}, s_{j_{i}+2}$.
Let $t'_{i}$ be the $x$ coordinate of $p_{i}$, then 
$0 \le t'_{i} - t_{i} \le \bub$.
The distance of each $t'_{i}$ from the closest point $s_{j}$
is at most $\f/8 - \bub \ge 0.1\f$.
It is also easy to check that 
$p_{i}, p_{i+1}$ satisfy~\eqref{eq:triangle-cnd.three}, so 
Lemma~\ref{lem:three} is indeed applicable.
 \end{pproof} % middle

 \begin{step+}{approx.u}
We have $u \leadsto p_{1}$ and $p_{n}\leadsto v$.
 \end{step+}
 \begin{pproof}
If $s_{1} \ge 0.1\f$, then the statement is proved by an
application of Lemma~\ref{lem:three}, so suppose $s_{1} < 0.1\f$.
Then $s_{1}$ belongs to a trap cover.

If $s_{2}$ belongs to a wall then $s_{2}\ge\f/3$, so
$s_{2}-s_{1}>0.23\f$.
If $s_{2}$ also belongs to a trap cover then
the reasoning used in Lemma~\ref{lem:one-of-three} gives $s_{2}-s_{1} > \f/4$.
In both cases, a midpoint $t_{1}$ was chosen between $s_{1}$ and $s_{2}$ with $t_{1}-s_{1}>0.1\f$,
and there is only $s_{1}$ between $u$ and $t_{1}$.     

If the trap cover belonging to $s_{1}$ is closer than $\g-6\bub$ then the fact
that $u$ is clean in 
$\cM^{*}$ implies that it contains a large trap-free region where it is
easy to get through.

Assume now that it is at a distance $\ge \g-6\bub$ from $u$.
Then we will pass through it, going from $u$ to $p_{1}$ similarly to
Part~\ref{step:three.x-free.trapc.vert} of the proof of Lemma~\ref{lem:three},
though the computations are a little different.

Define $w_{1}=\pair{x_{1}}{y_{1}}$ by $d(w_{1})=-L_{1}$.
We will apply Lemma~\ref{lem:if-correl} to vertical correlated traps 
$J\times I'$, with 
$J=\clint{x_{1}}{x_{1}+5\bub}$, $I'=\clint{y_{1}}{y_{1}+L_{1}}$.
The lemma is applicable since 
$u_{1}< y_{1}-L_{1}-\l_{1} < y_{1}+2L_{1}+\l_{1} < v_{1}$.
Indeed, using $\bub\ll\g$ and $\g/\f>\bub/\g$ from~\eqref{eq:bub-g-f},
as well as~\eqref{eq:triangle-cnd.three}, the inequality $x_{1}\ge\g-6\bub\ge\g/2$ implies
 \begin{align*}
 y_{1} &\ge u_{1} + (\g/2)\slope(u,v)-L_{1} \ge u_{1}+(\g/2)(\slopeincr-1)(\g/\f)
\\    &> u_{1}+(\g/2)(\slopeincr/2)(\bub/\g) = u_{1}+\slopeincr\bub/4>u_{1}+L_{1}+l_{1}.
 \end{align*}
The inequality $y_{1}+2L_{1}+\l_{1} < v_{1}$ as well as the continuation of the
proof follows similarly to 
Part~\ref{step:three.x-free.trapc.vert} of the proof of Lemma~\ref{lem:three}.

The relation $p_{n}\leadsto v$ is shown similarly.
 \end{pproof} % ends
 \end{Proof}

\section{Proof of Lemma~\protect\ref{lem:main} (Main)}
\label{sec:main-proof}

Lemma~\ref{lem:main} says that 
if $\m$ is sufficiently large then the sequence $\cM^{k}$ can be
constructed, in such a way that it satisfies all the above conditions and
also  $\sum_{k} \Prob(\cF'_{k}) < 1$.

The construction of $\cM^{k}$ is complete by the definition of $\cM^{1}$ in
Example~\ref{xmp:base} and the scale-up algorithm of
Section~\ref{sec:plan}, after fixing all parameters in Section~\ref{sec:params}.

We will prove, by induction, that every structure $\cM^{k}$ is a mazery.
We already know that the statement is true for $k=1$.
Assuming that it is true for all $i \le k$, we prove it for $k+1$.
The dependency properties in 
Condition~\ref{cond:distr}.\ref{i:distr.indep} are satisfied according to
Lemma~\ref{lem:distr.indep}.
The combinatorial properties in
Condition~\ref{cond:distr}.\ref{i:distr.combinat} have been proved in
Lemmas~\ref{lem:cover} and~\ref{lem:clean}.

The trap probability upper bound in
Condition~\ref{cond:distr}.\ref{i:distr.trap-ub}
has been proved in Lemma~\ref{lem:trap-scale-up}.
The wall probability upper bound in Condition~\ref{cond:distr}.\ref{i:distr.wall-ub}
has been proved in Lemma~\ref{lem:all-wall-ub}.
The cleanness probability lower bounds in
Condition~\ref{cond:distr}.\ref{i:distr.ncln}
have been proved in Lemma~\ref{lem:ncln-ub}.
The hole probability lower bound in
Condition~\ref{cond:distr}.\ref{i:distr.hole-lb} has been proved in 
Lemmas~\ref{lem:emerg-hole-lb}, \ref{lem:all-compound-hole-lb}
and~\ref{lem:heavy-hole-lb}.

The reachability property in Condition~\ref{cond:reachable} is satisfied via 
Lemma~\ref{lem:approx} (the Approximation Lemma).
There are some restrictions on the parameters
$\f, \g, \bub$ used in this lemma.
Of these, condition~\eqref{eq:bub-g-f}
holds if $\R_{0}$ is sufficiently large; the rest
follows from our choice of parameters and Lemma~\ref{lem:ncln-ub}.

Finally, $\sum_{k} \Prob(\cF'_{k}) < 1$ 
follows from Lemma~\ref{lem:0.6}.

The main result, Theorem~\ref{thm:main}, has been proved in
Section~\ref{sec:outline}, using Lemma~\ref{lem:main} and some conditions.
Of these, 
the reachability property in Condition~\ref{cond:reachable-hop} is implied by
Condition~\ref{cond:reachable}.
The property saying that the set of clean points is sufficiently dense,
Condition~\ref{cond:dense}, is implied by
Remark~\ref{rem:distr}.\ref{i:middle-third}.
The property saying that the absence of lower-level bad events near the origin
imply $k$-level cleanness of the origin, 
Condition~\ref{cond:Qbd}, follows immediately from the definition of
cleanness.

\section{Conclusions}\label{sec:concl}

It was pointed out in~\cite{GacsChat04} that the clairvoyant demon does not
really have to look into the infinite future, it is sufficient for it to
look as far ahead as maybe $n^{3}$ when scheduling $X(n), Y(n)$.
This is also true for the present paper.

Another natural question is: how about three independent random walks?
The methods of the present paper make it very likely that three independent
random walks on a very large complete graph can also be synchronized, but
it would be nice to have a very simple, elegant reduction.

It seems possible to give a common generalization of the model of
the paper~\cite{GacsChat04} and the present paper.
Let us also mention that 
we have not used about the independent Markov processes $X,Y$ the fact that
they are homogenous: the transition matrix could depend on $i$.
We only used the fact that for some
small constant $w$, the inequality
$\Pbof{X(i+1) = j \mid X(i)=k} \le w$ holds for all $i,j,k$ (and similarly
for $Y$).

What will strike most readers as the most pressing open question is how to
decrease the number of elements of the smallest graph for which scheduling
is provably possible from super-astronomical to, say, 5.
Doing some obvious optimizations on the present renormalization method is
unlikely to yield impressive improvement: new ideas are needed.

Maybe computational work can find the better probability
thresholds needed for renormalization even on the graph $K_{5}$,
introducing supersteps consisting of several single steps.

\subsection*{Acknowledgement}

I am grateful to Peter Winkler and M\'arton Bal\'azs for valuable comments,
to the referee of paper~\cite{GacsChat04}, to Fred Kochman.
Particular thanks are due
to the referee of an earlier version of the present paper
and to referees of the current version for repeated careful
reading and thorough criticism.
(Of course, the responsibility for the remaining bugs is still mine.)
I also thank J\'anos Koml\'os and Endre Szemer\'edi, who first
challenged me to solve this problem.
G\'abor Tusn\'ady's advices helped readability.

% \bibliographystyle{plain}
% \bibliography{reli,publ}

\newpage

\section*{Glossary}

\subsection*{Concepts}

\begin{flushdescription}

\item[Approximation lemma]
This is the lemma with the longest proof, though the fact that something like it
can be proven is not surprising.
The lemma says that the reachability condition will hold in the higher-order model
$\cM^{k+1}=\cM^{k*}$, and is part of the project to prove that $\cM^{k+1}$ is a
mazery.
The statement and the proof is elementary geometrical:
it can be translated as saying 
that one can pass through a network of walls that are well-separated
and densely perforated by holes, 
with only well-separated traps between them.

The proof is long only because of the awkwardly large number of cases to
consider, though these cases differ from each other only slightly.

\item[Barrier] Barriers are introduced in Section~\ref{sec:mazery}.
Barrier is the same type of object as a wall.
A mazery has a random set of barriers just as a random set of walls.
The set $\cW_{\d}$ of walls is a subset of the set $\cB_{\d}$ of barriers.
The simpler definition of barriers in the renormalization operation of
Section~\ref{sec:plan} allows the computation of probability upper bounds.
The definition of walls guarantees certain combinatorial properties
allowing us to show how to pass through and between them.

\item[Barrier probability upper bound $p(\r)$]  Part of the definition of a mazery,
  this function upper-bounds the probability of a barrier of rank $\r$ arising at a
  specific point.

Moreover, as said in Condition~\ref{cond:distr}.\ref{i:distr.wall-ub},
it will be an upper bound to $\sum_{l}\p(\r,l)$, where $\p(\r,l)$ is
the probability of a barrier of rank $\r$ and length $l$ arising.

The actual size of $\p(\r)$ will depend negative exponentially on $\r$, as
defined in Section~\ref{sec:params}.

\item[Body of a wall]
A right-closed interval $\Body(E)$, part of the definition of a wall or barrier.

\item[Cleanness] Introduced informally in Subsection~\ref{subsec:application}, and
  formally in Section~\ref{sec:mazery}.
There are several kinds of cleanness, but they all express the same idea.
For example, informally, the fact that a 
point $u$ is upper right trap-clean in mazery $\cM^{k}$
stands for the property that for all $i<k$, 
there is no $i$-level trap closer than approximately $\g_{i}$ to $u$ in
the upper right quadrant starting from $u$.
Formally, the random sets $\cS_{\d}$ for $\d=0,1,2$ and $\cC_{\d}$ for $\d=0,1$
describe cleanness.

The combination of several kinds of cleanness can result in some
complex versions, introduced in Definition~\ref{def:H-clean}.
Of these, only H-cleanness will be much used, it is suited for
proving probability lower bounds on holes.

Cleanness must appear with a certain density: this is spelled out in
Conditions~\ref{cond:distr}.\ref{i:distr.clean.1}
and~\ref{cond:distr}.\ref{i:distr.clean.2}.
Actually every point must be clean (in all kinds of sense) with rather large
probability.
This is spelled out in Condition~\ref{cond:distr}.\ref{i:distr.ncln}.

\item[Closed and open points] Defined in Subsection~\ref{subsec:perc}.

\item[Collision] Only used in the introduction.

\item[Compound trap] Introduced informally in Subsection~\ref{subsec:application},
  and formally in the scale-up operation of Section~\ref{sec:plan}.
The event that some traps are too close to each other.
It can be \df{uncorrelated}, when the projections of the traps are disjoint, and
\df{correlated} otherwise.

\item[Compound wall] Introduced informally in Subsection~\ref{subsec:application},
  and formally in the scale-up operation of Section~\ref{sec:plan}.
The event that some walls are too close to each other.
The \df{rank} increases with the rank of its components but decreases with their
distance.

\item[Correlated trap] Introduced 
formally in the scale-up operation of Section~\ref{sec:plan}.
The event that some traps are too close to each other
when the projections of the traps are not disjoint, but still \emph{both
conditional probabilities} are small.
When one of the conditional probabilities is not small, this will give rise to
an emerging wall.

\item[Compatible] Only used in the introduction.

\item[Delay sequence] Only used in the introduction.

\item[Dependencies and monotonicities] The probability estimates also
control which random objects depend on what parts of the original
random process $Z$.
The requirements are formulated in Condition~\ref{cond:distr}.\ref{i:distr.indep}.

\item[Dominant wall] Introduced in Definition~\ref{def:barrier}, it is a wall that
is surrounded by external intervals of size $\ge\bub$.

\item[Emerging barriers and walls] 
Introduced informally in Subsection~\ref{subsec:application},
and formally in Section~\ref{sec:plan}.
Suppose that an event $\cA$ appears that would lead to a
trap of the correlated or missing-hole kind.
If the conditional probability of $\cA$ (with respect to, say, the process $X$)
is not small, then we do not define a trap, but then
this is an improbable event as a function of the process 
$X$ and defines an emerging (vertical) barrier.

Before deciding which of the emerging barriers become walls, we introduce some
cleanness requirements in the notion of \df{pre-wall}, and then apply a
selection process to make walls disjoint.

\item[Exponents] The construction deals with several exponents,
all of which will be fixed, but chosen appropriately to satisfy certain
inequalities.

$\hxp = 0.015$: Informally, if the probability of the occurrence of a certain kind of wall is
upper-bounded by $p$ then the probability, at a given site, to find a hole through the
wall is lower-bounded by $p^{\hxp}$.

The other exponents are defined by their role in the expression of the
parameters already seen elsewhere, as set in Definition~\ref{def:exponential}:
 \begin{align*}
          \T       &= \lg^{R}\text{ with } \lg = 2^{1/2},
\\        \bub     &= \T^{\bubxp},
\quad     \f        = \T^{\fxp},
\quad     \g        = \T^{\gxp},
\quad     \tub      = \T^{-\tubxp}\text{ with }\tubxp=4.5,
\\ \R_{k} &=\R_{0}\txp^{k} \text{ with }\txp = 2-\fxp.
 \end{align*}
The choice of these parameters is made in Lemma~\ref{lem:exponents-choice},
where the necessary inequalities are spelled out.

\item[External interval] Introduced in Definition~\ref{def:barrier}, a
right-closed interval that intersects no walls.

\item[Graph $\cG=(\cV,\cE)$]  This random graph
is the same throughout the proof.
We have $\cV=\bbZ_{+}^{2}$, and $\cE$ consists of
all rightward and upward edges between open lattice points.

\item[Hole] Formally defined in Definition~\ref{def:holes}, a 
vertical hole is defined by
a pair of points on the opposite sides of a horizontal wall,
and an event that it is possible to pass from one to the other.

The hole is called \df{good} if it is lower-left and
upper-right H-clean (see the reference to partial cleanness).
This partial cleanness requirement is useful in the proof of the hole
probability lower bound.
There is also a qualified version of goodness, as ``seen from'' a point $u$.

\item[Hole probability lower bound]
As said in the item on the exponent, the lower bound of the probability of a
hole through a wall of rank $\r$ is essentially $\p(\r)^{\hxp}$.
For technical reasons, it is larger by a polynomial factor, as given 
in~\eqref{eq:h-def}.
In the interests of the inductive proof (in order to handle
compound walls), the actual condition on the hole probability lower bound 
is rather complex, and is spelled out in
Condition~\ref{cond:distr}.\ref{i:distr.hole-lb}.

\item[Hop] A hop, as defined in Definition~\ref{def:hop},
is an inner clean right-closed interval that contains no walls.
The role of a hop is that it is an interval that is manifestly possible to pass
through.

\item[Isolated light wall] This concept is only used in the proof of
  Lemma~\ref{lem:cover}.

\item[Light and heavy walls] Those of rank $<\R^{*}$ and $\ge\R^{*}$.

\item[Mazery $\cM^{k}$] Introduced in Subsection~\ref{subsec:renorm} in an abstract
  way just as a $k$-level ``model''.
Mazeries are defined in Section~\ref{sec:mazery}, and the operation defining
$\cM^{k+1}$ in terms of $\cM^{k}$ is defined in Section~\ref{sec:plan}.
A mazery is a random process containing a random set of traps, random sets of
barriers and walls, and a random assignment of some cleanness properties to points.
These random processes must obey certain conditions, spelled out in
Section~\ref{sec:mazery}.

\item[Number of colors $\m$]  Introduced in the Introduction, it is the size of
the complete graph, and the number of colors in the color-percolation model.
It is connected to the other parameters of the mazeries via~\eqref{eq:m-expr}.

\item[Potential wall] 
Whether a segment is a barrier or not depends only on symbols in the segment
itself.
However, whether a barrier is a wall depends on the environment.
As introduced in Definition~\ref{def:potential-wall},
a segment is a \df{potential wall} when there is an environment making it a wall.

\item[Power law convergence] Only used in the introduction.

\item[Pure sequence] This concept is used only in the proof of
  Lemma~\ref{lem:cover}.

\item[Rank of a wall] This is a value that classifies walls to help in the
  definition of the renormalization operation.
Higher rank will mean lower bound on the probability.
When forming mazery $\cM^{k+1}$ we will not delete all walls of mazery $\cM^{k}$,
only walls of rank less than $\R^{*}$ (the so-called \df{light walls}), as well
as certain heavy walls contained completely in light walls.
Also, in forming compound walls one of the components will be light.

\item[Rank lower bound $\R$]  Introduced in~\ref{sec:mazery}, it is part of the
  definition of a mazery.
In our series of mazeries $\cM^{k}$, we will have $\R_{1}<\R_{2}<\dotsm$.
Actually, $\R_{k}$ grows exponentially with $k$, 
as defined in Section~\ref{sec:params} and above in the description of exponents.

(The rank computation for compound walls
suggests to view rank as analogous to ``free energy'' in statistical 
physics (energy minus entropy), but ignore the analogy if it only confuses you.)

\item[Rank upper bound]
No apriori rank upper bound is given, but the construction guarantees a
that all ranks of a mazery are between lower bound $\R$ and upper bound
$\txpub\R$ for a certain constant $\txpub$.
This will guarantee that each rank value is present in a mazery $\cM^{k}$ for
only a constant number of values of $k$.

\item[Rectangle] Notation for rectangles is introduced in
  Subsection~\ref{subsec:notation}.

\item[Reachability] The notation $\leadsto$ is introduced in
  Subsection~\ref{subsec:notation}.
In a mazery, the reachability condition, Condition~\ref{cond:reachable} says,
essentially, that one can pass from a clean point to another one, if there are
no obstacles (walls, holes) in between, and the slope between the two points is
not too small or too large.

\item[Renormalization, or scale-up] The word is used only informally, as the
  operation $\cM\mapsto \cM^{*}$ that leads from mazery $\cM^{k}$ to mazery
  $\cM^{k+1}$.

\item[Scale parameter $\bub_{k}$] Upper bound on the size of traps and walls in
  the model of level $k$.
The parameter $\bub$ is part of the definition of a mazery.

\item[Scale-up parameters]
The parameters $\g\ll \f$ play a special role in the scale-up operation.
We will have $\bub\ll\g$.
They play several roles, but the most important is this:
We try not to permit walls closer than $\f$ to each other,
and we try not to permit intervals larger than $\g$ without holes.
We will also use $\f$ in the definition of one-dimensional cleanness and $\g$ in
the definition of trap-cleanness.
(The value $\g'=2.2\g$ plays a subsidiary role.)

\item[Sequence of neighbor walls]
A sequence of neighbor walls, as introduced in Definition~\ref{def:neighbor-seq},
is a sequence of walls separated by hops.
Such a sequence is useful since it allows the analysis of
passage through it.
It is a crucial part of the combinatorial construction to establish the
existence of such sequences (under the appropriate conditions) in $\cM^{k+1}$.
The existence requirement is given
in Conditions~\ref{cond:distr}.\ref{i:distr.inner-clean} 
and~\ref{cond:distr}.\ref{i:distr.cover}.
Of course, the proof uses the corresponding property of $\cM^{k}$.

\item[Sequences of random variables]
The basic sequences we consider are the two independent random walks
 $X_{1},X_{2},\dots$ and $Y_{1},Y_{2},\dots$.
But for notational convenience, we sometimes denote $X=Z_{0}$, $Y=Z_{1}$
and $Z=(X,Y)$.
These are the only source of randomness, 
all other random structures introduced later (mazeries $\cM^{k}$)
are functions of these.

\item[Slope constraint] For a pair of points $u,v$ the requirement that the
  slope $s$ of the segment between them be neither too horizontal nor too
  vertical. 
This is expressed by saying $\slb \le s\le 1/\slb$.
The value of the bound $\slb$ belongs to the mazery.
It is $\slb_{k}$ in mazery $\cM^{k}$.
We will have $\slb_{1}<\slb_{2}<\dots<1/2$.
The notation $\minslope(u,v)$ is the smaller of the two slopes of the segment
between points $u$ and $v$, as defined in Subsection~\ref{subsec:notation}.

\item[Strong cleanness] One-dimensional cleanness has a strong version,
  introduced in Definition~\ref{def:1dim-clean}.
The motivation is similar to introduction of barriers, and in fact
in the renormalization operation of Section~\ref{sec:plan}, strong
cleanness will come from the absence of lower-level barriers.

\item[Trap] Introduced informally in Subsection~\ref{subsec:application}, and
  formally in Section~\ref{sec:mazery}, where the set of traps is denoted by
  $\cT$.
On the lowest level, this is where two equal colors collide.
On level $k+1$, it is either the event that two close $k$-level traps occur,
or that some local bad event of a new kind (emerging on this level) occurs.
An example is a trap of the \df{missing-hole} kind: when on a certain wall of
level $k$, there is a long segment without a hole.
The same applies with horizontal and vertical interchanged.

\item[Trap cover]
In the proof of the Approximation Lemma, a strip that covers a trap.
The idea is that under the conditions of that lemma, after covering
all traps with trap covers, all the trap covers and the walls still form a
rather sparse network, which it will be possible to pass through.

\item[Trap probability bound $\tub$] Introduced in Section~\ref{sec:mazery},
it is part of the definition of a mazery.
It bounds the \df{conditional probability} of a trap, as spelled out in
Condition~\ref{cond:distr}.\ref{i:distr.wall-ub}.
The value $w_{k}$ decreases super-exponentially in $k$, 
as defined in Section~\ref{sec:params}.

\item[Trap-cleanness] Introduced in Definition~\ref{def:trap-clean}.
See the item on cleanness.

\item[Ultimate bad event $\cF_{k}$] Introduced in Subsection~\ref{subsec:renorm}.

\item[Wall] Introduced informally in Subsection~\ref{subsec:application}, and
  formally in Section~\ref{sec:mazery} as an object having a \df{body} and a
  \df{rank}.
On level $k+1$, it is either the event that two close $k$-level walls occur
(\df{compound wall}),
or that some bad event of a new kind (emerging on this level) occurs in one of
the projections (\df{emerging wall}).
An example of such an event: 
when the conditional probability that two close traps occur is too high.

\item[Winkler percolation] Only used in the introduction.

\end{flushdescription}

\subsection*{Symbols}

\begin{itemize}

\item[$\cB_{0}$] Set of vertical barriers.

\item[$\cC$] Structure defining the various kinds of one-dimensional
  cleanness. 

\item[$\aux_{1},\aux_{2}$] Constants used in the definition of $\p(\r)$.
$\aux_{1}$ is fixed at the end of the proof of Lemma~\ref{lem:compound-contrib}.
$\aux_{2}$ can be fixed as anything $<1-1/\lg$, as said 
at the end of the proof of Lemma~\ref{lem:pub}.

\item[$\aux_{3}$] Constant used in the definition of $\h(\r)$.
Fixed at the end of the proof of Lemma~\ref{lem:all-compound-hole-lb}.

\item[$\d$] Frequently used to denote an index 0 or 1, with $Z_{0}=X$,
  $Z_{1}=Y$.

\item[$\bub$] Upper bound on the size of walls and traps.
Fixed as $\T^{\bubxp}$ in Definition~\ref{def:exponential}.

\item[$\bubxp$] Exponent used in defining $\bub$.
Chosen $0.15$ in the proof of Lemma~\ref{lem:exponents-choice}.

\item[$\cE$] Set of edges of the random graph $\cG$ in $\bbZ_{+}^{2}$
defined by the processes $X,Y$.

\item[$\f$] Lower bound on the distance between walls.
Fixed as $\T^{\fxp}$ in Definition~\ref{def:exponential}.

\item[$\fxp$] Exponent used in defining $\f$.
Chosen $0.25$ in the proof of Lemma~\ref{lem:exponents-choice}.

\item[$\cG$] The random graph $(\cV,\cE)$ in $\bbZ_{+}^{2}$
defined by the processes $X,Y$.

\item[$\g$] Upper bound on the length of a wall segment without holes
(has also some other roles).
Fixed as $\T^{\gxp}$ in Definition~\ref{def:exponential}.

\item[$\g'$] Plays a subsidiary role in defining correlated traps.
Defined as $2.2\g$ in~\eqref{eq:g-prime}.

\item[$\gxp$] Exponent used in defining $\g$.
Fixed as $0.2$ in the proof of Lemma~\ref{lem:exponents-choice}.

\item[$\h(\r)$] Hole lower bound probability.
Fixed as $\aux_{3}\lg^{-\hxp\r}$ in~\eqref{eq:h-def}.

\item[$\hxp$] Exponent used in lower-bounding hole probabilities.
Fixed as $0.015$ in~\eqref{eq:hxp}.

\item[$H$] Constant defined as $12$ in~\eqref{eq:H-def} and used in the proof of
  the Approximation Lemma.

\item[$k$] Frequently denotes the level of the mazery $\cM^{k}$.

\item[$\l_{j}$] We defined $\l_{1}=7\bub$, $\l_{2} = \g'$ in connection with the
  definition of correlated traps in the scale-up.

\item[$L_{j}$] We defined $L_{1}=4\l_{1}$, $L_{2} = 4\l_{2}$, $L_{3}=\g$.
These parameters 
determine approximately the widths of emerging barriers of three types.

\item[$\cL_{j}$] In the scale-up operation, for 
$j=1,2,3$, the ``bad event'' $\cL_{3}(x,y,I,b)$ triggers
the occurrence of a new trap provided its conditional probability 
(with condition $X=x$) is small.
The event (a function of $x$) that its conditional probability
is large triggers the occurrence of an
emerging vertical barrier of the corresponding type.

\item[$\lg$] The lower base of our double exponents.
Fixed as $2^{1/2}$ in Definition~\ref{def:lg-def}.

\item[$\slopeincr$] Constant used in the definition of $\slb_{k+1}$
and in some bounds.
Fixed as $500$ in Definition~\ref{def:f-g}.

\item[$\m$] Number of elements of the complete graph where the random walk
  takes place, or equivalently, the number of colors in the color percolation
  model.
Lower-bounded by $1/\tub_{1}+1$ in Example~\ref{xmp:base} which serves as the base
case $\cM^{1}$; later $\tub_{1}$ is chosen to make this lower bound exact.

\item[$\cM^{k}$] Mazery of level $k$.

\item[$\p(\r)$] Probability upper bound of a barrier of rank $\r$.
Fixed as $\p(\r) = \aux_{2} \r^{-\aux_{1}} \lg^{-\r}$
in~\eqref{eq:wall-prob}.

\item[$\p(\r,l)$] Supremum of the probabilities of any barrier of
  rank $\r$ and length $l$ starting at a given point.
Depends on the mazery, but the conditions require
$\sum_{l}\p(\r,l)\le\p(\r)$, where the function $\p(\r)$ is fixed.

\item[$\r$] Frequently denotes a rank.

\item[$\R$] Rank lower bound.
Fixed as $\R_{k}=\R_{0}\txp^{k}$ in Definition~\ref{def:ranks}.

\item[$\hat\R$] Rank of emerging walls.
Fixed as $\txp'\R$ in Definition~\ref{def:ranks}.

\item[$\R_{0}$] There are several lemmas that hold when it is chosen sufficiently
  large, and there are no other conditions on it.

\item[$\cS$] Structure describing the one-dimensional kinds of strong
  cleanness and also trap cleanness.

\item[$\slb$] Slope lower bound.
It is defined by $\slb_{1}=0$ in~\eqref{eq:slb1}, and by
$\slb_{k+1} = \slb_{k} + \slopeincr \g_{k} / \f_{k}$ in
Definition~\ref{def:new-slb}.

\item[$\T$] Auxiliary parameter for defining several other parameters.
Fixed as $\lg^{\R}$ in Definition~\ref{def:exponential}.

\item[$\cT$] The set of traps.

\item[$\txp$] Used in the definition of $\R_{k}$.
Fixed as $2-\fxp$ in Definition~\ref{def:exponential}.

\item[$\txp'$] Coefficient used in the definition of $\hat\R$.
Fixed as $2.5$ in the proof of Lemma~\ref{lem:exponents-choice}.

\item[$\cV$] Set of points of the random graph $\cG$ in $\bbZ_{+}^{2}$
defined by the processes $X,Y$.

\item[$\cW_{0}$] Set of vertical walls.

\item[$\tub$] Upper bound on the probability of traps.
Fixed as $\T^{-\tubxp}$ in Definition~\ref{def:exponential}.

\item[$\tubxp$] Exponent in the definition of $\tub$.
Fixed as $4.5$ in Definition~\ref{def:exponential}.
  
\item[$X$] The sequence $X(1),X(2),\dots$ is a random walk over the complete graph $K_{\m}$.

\item[$Y$] The sequence $Y(1),Y(2),\dots$ is a random walk over the complete graph $K_{\m}$.

\item[$Z$] We defined $Z_{0}(i) = X(i)$, $Z_{1}(i)=Y(i)$.

\end{itemize}

\end{document}